\documentclass[opre,nonblindrev]{informs3}

\OneAndAHalfSpacedXI % current default line spacing
\usepackage{endnotes}
\let\footnote=\endnote
\let\enotesize=\normalsize
\def\notesname{Endnotes}%
\def\makeenmark{\hbox to1.275em{\theenmark.\enskip\hss}}
\def\enoteformat{\rightskip0pt\leftskip0pt\parindent=1.275em
  \leavevmode\llap{\makeenmark}}
\def\bfred#1{{\color{red}\bf#1}}
\usepackage{graphicx}
\usepackage{subfigure}
\usepackage{natbib}
 \bibpunct[, ]{(}{)}{,}{a}{}{,}%
 \def\bibfont{\small}%
\usepackage{hyperref}
\hypersetup{hypertex=true,
            colorlinks=true,
            linkcolor=blue,
            anchorcolor=blue,
            citecolor=blue}
\TheoremsNumberedThrough   % Preferred (Theorem 1, Lemma 1, Theorem 2)
\ECRepeatTheorems
\usepackage{url}
\EquationsNumberedThrough    % Default: (1), (2), \cdots
\usepackage{subfigure, epsfig, epstopdf}
\usepackage{graphicx,amssymb}
\usepackage{multirow, multicol}
\usepackage{amsmath,etoolbox}
\usepackage{color}
\usepackage{setspace}
\usepackage{mathrsfs}
\usepackage{dsfont}
\usepackage{bbm}
\usepackage{textcomp}

\usepackage{tikz} 
\usetikzlibrary{shapes.geometric, arrows} 
\usetikzlibrary{snakes}
\usepackage[]{algorithm}
\usepackage{algpseudocode}
\usepackage{array}
\usepackage{threeparttable}
\usepackage{graphicx}
\usepackage{verbatim}

\usepackage{multibib}
\newcites{appendix}{References}

\newcommand{\td}{\mathrm{d}}

\newcommand{\bm}[1]{\boldsymbol #1}
\def\bfred#1{{\color{red}\bf#1}}

\def\ye#1{{\color{red}#1}}

\patchcmd{\subequations}
 {\def\theequation{\theparentequation\alph{equation}}}
 {\def\theequation{\theparentequation.\arabic{equation}}}
 {}{}

\graphicspath{{figure/}}

\begin{document}
%%%%%%%%%%%%%%%%

%\setcounter{page}{1}
%\DOI{10.1287/xxxx.0000.0000}%

\RUNAUTHOR{Sun et al.}

\RUNTITLE{Optimal abort under imperfect condition monitoring}

% \TITLE{Bundling Information Goods of Decreasing Value}
\TITLE{Optimal Abort Policy for Mission-Critical Systems under Imperfect Condition Monitoring}

\ARTICLEAUTHORS{%
	
%\iffalse
%% Enter all authors
\AUTHOR{Qiuzhuang Sun,$^1$ Jiawen Hu,$^2$ and Zhi-Sheng Ye$^3$
%\footnote{Corresponding author.}
}
\AFF{$^1$School of Mathematics and Statistics, University of Sydney, Australia} %, \URL{}}
\AFF{$^2$School of Aeronautics and Astronautics, University of Electronic Science and Technology of China, China} %, \URL{}}
\AFF{$^3$Department of Industrial Systems Engineering and Management, National University of Singapore, Singapore} %, \URL{}}
\AFF{Emails: \texttt{qiuzhuang.sun@sydney.edu.au}; \texttt{hdl@sjtu.edu.cn}; \texttt{yez@nus.edu.sg}}
%\fi
} % end of the block

\ABSTRACT{%
While most on-demand mission-critical systems are engineered to be reliable to support critical tasks, occasional failures may still occur during missions.
To increase system survivability, a common practice is to abort the mission before an imminent failure.
We consider optimal mission abort for a system whose deterioration follows a general three-state (normal, defective, failed) semi-Markov chain.
The failure is assumed self-revealed, while the healthy and defective states have to be {inferred} from imperfect condition monitoring data.
Due to the non-Markovian process dynamics, optimal mission abort for this partially observable system is an intractable stopping problem.
For a tractable solution, we introduce a novel tool of Erlang mixtures to approximate non-exponential sojourn times in the semi-Markov chain.
This allows us to approximate the original process by a surrogate continuous-time Markov chain whose optimal control policy can be solved through a partially observable Markov decision process (POMDP).
We show that the POMDP optimal policies converge almost surely to the optimal abort decision rules when the Erlang rate parameter diverges.
This implies that the expected cost by adopting the POMDP solution converges to the optimal expected cost.
Next, we provide comprehensive structural results on the optimal policy of the surrogate POMDP.
Based on the results, we develop a modified point-based value iteration algorithm to numerically solve the surrogate POMDP.
We further consider mission abort in a multi-task setting where a system executes several tasks consecutively before a thorough inspection.
Through a case study on an unmanned aerial vehicle, we demonstrate the capability of real-time implementation of our model, even when the condition-monitoring signals are generated with high frequency.
}%

%\KEYWORDS{deterministic inventory theory; infinite linear programming duality;
%  existence of optimal policies; semi-Markov decision process; cyclic schedule}


\maketitle
%%%%%%%%%%%%%%%%%%%%%%%%%%%%%%%%%%%%%%%%%%%%%%%%%%%%%%%%%%%%%%%%%%%%%%

% Samples of sectioning (and labeling) in OPRE
% NOTE: (1) \section and \subsection do NOT end with a period
%       (2) \subsubsection and lower need end punctuation
%       (3) capitalization is as shown (title style).
%
%\section{Introduction.}\label{intro} %%1.
%\subsection{Duality and the Classical EOQ Problem.}\label{class-EOQ} %% 1.1.
%\subsection{Outline.}\label{outline1} %% 1.2.
%\subsubsection{Cyclic Schedules for the General Deterministic SMDP.}
%  \label{cyclic-schedules} %% 1.2.1
%\section{Problem Description.}\label{problemdescription} %% 2.

% Text of your paper here

\setcounter{equation}{0}
\setcounter{lemma}{0}
\setcounter{section}{0}
\setcounter{theorem}{0}
\renewcommand*{\theHequation}{\arabic{section}.\arabic{equation}} 
\renewcommand*{\theHlemma}{\arabic{section}.\arabic{lemma}} 
\renewcommand*{\theHsection}{\arabic{section}.\arabic{section}} 
\renewcommand*{\theHtheorem}{\arabic{section}.\arabic{theorem}} 

\newpage

\section{Introduction}\label{sec:intro}
%\subsection{Background and Motivation}\label{subsec:intro_motivation}
{An on-demand mission-critical system is required to operate continuously for a period of time to complete a mission \citep{kim2010contracting}.
Classical examples include commercial jetliner serving routes between cities, emergency response vehicles like ambulance and fire engines, 
and medical devices like a cardiopulmonary bypass machine used in cardiac surgery.
Some recent examples include autonomous underwater vehicles for underwater inspection and maintenance in offshore infrastructure, automated waste collection systems for marine debris removal, and unmanned aerial vehicles (UAVs) used for such tasks as wildfire monitoring, product delivery, combat, and power grids inspection.
% Other examples of mission-critical systems include pulp-and-paper mills manufacturing papers \citep{ranjan2018dataset}, online data-processing systems performing data processing tasks \citep{levitin2021optimal}, and chemical reactors with a cooling system for production tasks \citep{cha2018optimal}.
While these systems are engineered to be reliable to support critical tasks, occasional failures may still occur during missions.
As evidence, a record number of 20 large Air Force drones crashed in major accidents in 2015 according to the Washington Post \citep{whitlock2016more};
electrical failure of the cardiopulmonary bypass device is estimated to be 1 per 1,000 cases \citep{durukan2016electrical};
and the helicopter emergency medical services in the US reported an average 4.5 accidents per 100,000 flying hours, out of which around one third were fatal \citep{bryan2014analysis}.
The system failure causes not only a mission failure, but also a huge system failure cost including possible threaten to human lives. 
To increase system survivability, a common practice is to abort the mission before an imminent failure.
% For example, if we predict that a cooling system of a chemical reactor cannot provide a desired temperature level, the production task using the reactor should be aborted \citep{cha2018optimal}.
For instance, an autonomous vehicle or crewed aircraft diagnosed to be defective could abort its mission and return to the base for crash prevention \citep{yang2019designing};
a malfunctioning ambulance might pull over at a location convenient for the next ambulance to take over the patient, rather than continuing to travel and risking a breakdown on the road;
in the event of an abnormal heart-lung machine, a perfusionist can manage a patient's venous return during weaning from cardiopulmonary bypass and subsequently employ less effective/reliable hand-cranking for systemic perfusion \citep{durukan2016electrical}.
	
Several distinctive features of emergent mission abort can be discerned from the above examples.
Each on-demand mission typically has a short duration when compared to the service life of the system.
Throughout the mission, advanced sensors monitor system health in real time and trigger alarms when a potentially suspicious defect is detected.
Abort, as an option to respond to alarms, has a non-negligible duration compared with the mission duration.
During abort, a system failure is possible, and thus abort may be unnecessary if it takes longer time than completing the mission.
Intuitively, an abort decision balances the system failure risk and the zero-one yield from mission success over a short-term horizon when the system is deemed defective.
These features make abort fundamentally different from condition-based maintenance (CBM),
which monitors system ageing and preventively replaces the system when it reaches wear-out with unacceptably high risk of breakdown.
CBM trades off uptime against system failure costs over an infinite horizon or a horizon significantly longer than a system's service life.
Because of {\it the linear yield in uptime, the long planning horizon, and the commonly assumed negligible maintenance time}, the optimal CBM policy {typically has a simple control limit or a simple control-limit function}. 
See \cite{elwany2011structured,chen2015condition,drent2023real} when system degradation follows a continuous-state stochastic process with random effects,
\cite{kim2016robust} and \cite{zhu2021robust} when system deterioration follows a multi-state Markov chain with model misspecification,
and \cite{khaleghei2021optimal} when system failure follows a three-state semi-Markov process.
In contrast, {\it the zero-one yield, the short planning horizon, and the non-negligible abort duration} rule out {a simple control-limit policy} for mission abort decisions when the state evolution is non-Markovian.
{This is evidenced by Theorem~\ref{thm:struct_policy} ahead, 
where we show that the optimal policy is defined by lower and upper control-limit functions under a spherical coordinate system, when the system deterioration follows a partially observable three-state (normal, defective, failed) semi-Markov chain with phase-type sojourn times.}
% The deep reason is xxx. \sun{[QZ: it seems very difficult to characterize the optimal policy for the original problem; currently I mentioned the policy in Thm~4 is for .]}
%\ye{[one sentence to support the above argument]}
	
For a partially observable three-state system whose health status is inferred from condition-monitoring data and whose failure is self-revealed, existing solutions base the abort decision on a charting statistic that monitors its latent state.
%\ye{[it would be better to link to the literature here. So we can have a thorough review on all rule-based studies here]}. 
The mission is aborted once the charting statistic triggers an alarm.
This is a two-step procedure where the decision rule in the second step is predicated on the alarm generated from the first step.
Some studies considered a more complicated decision rule for mission abort. 
For example,
\cite{qiu2019gamma} classify a system as defective upon receiving a defect signal, leading to a mission abort if the remaining time to complete the mission exceeds an optimized threshold.
\cite{cha2018optimal} confirm the defective status once $\bar m$ warning signals are received in total, and abort the mission if the remaining mission duration is longer than a time threshold $\bar\xi$.
Both $\bar m$ and $\bar\xi$ can be optimized.
\cite{zhao2021multi} compute the expected time for the system to reach a certain degradation level.
The mission is aborted if this time is smaller than a threshold. 
%\ye{[more reviews here]}
% Alternatively, we can predict the system remaining useful life (RUL) based on condition-monitoring signals. Aborting the system occurs when the expected RUL is less than the remaining mission time.  In both cases, we implicitly predetermine a default abort decision rule. However, the structure of abort policies is still predetermined.
These rule-based methods, albeit easy to implement, fall short in optimizing the trade-off between system failure and mission success.
	
Optimality guarantee of an abort policy is possible if the underlying system dynamics governing the system state transition and condition-monitoring data generation are known, but this direction remains unexplored.
For a three-state Markovian deteriorating system with condition-monitoring data, its optimal mission abort control can be attained through a partially observable Markov decision process (POMDP), and our Corollary~\ref{thm:control_limit_CTMC} ahead shows that the optimal policy enjoys a relatively simple structure.
Given that the mission duration is typically much shorter than the system service life, it is reasonable to posit an exponential sojourn time from the healthy state to either the defective or failure state.
However, assuming an exponential sojourn time from defective to failure states may be questionable, as this distribution is contingent on the endogenous physical mechanisms governing the progression from a defect to the ultimate collapse.
This argument is well supported by real-world lifetime datasets from physical systems \citep{mittman2019hierarchical,meeker2022statistical}.
%\ye{[could you help me find a few relevant ref that really explain the phenomenon why transition time from defective to failure should not  be exponential?]}
Violation of the Markovian assumption can severely invalidate solutions from a POMDP.
As an illustration, Section~\ref{sec:num} considers a bimodal sojourn time from defect to failure.
Approximating such a system by a simple Markov process loses critical information about the sojourn time in the defective state, leading to inferior decisions.
Our numerical results show that the cost increment from this simple approximation can be near 20\% compared to the proposed method.
%\ye{[you use this example in the simulation? Then please explicitly tell readers what is the performance loss directly. For example, 30\% cost increment? give exact numbers here.]},  
%\ye{[5min and 3h after defect? or after mission? Is this used in Section 9? Then I should not call it a toy example. Rather, your `extreme example' would be much better]}
	
\iffalse	
\ye{the logic should be like this: in this study, we use semi-Markov for process dynamics. Its opt turns out to be intractable, then link to the optimal stopping formula later in the paper. Then you say that we introduce a phase-type approximation framework to derive an approximation to the optimal policy. then briefly introduce it: we use phase type to appx sojourn times from xx to xx, and from xx to xx. while the tran time from xx to xx is assumed exponential. 
After that, you explain the novelty of the framework: the approx has to be like this otherwise the resulting system cannot be analyzed -- link to Fig 1 and say that detailed discussions are given there. 
Then state your theoretical contribution -- our approx can be made arbitrarily accurate, which is the result in theorem xxx.
as such, using the policy can guarantee optimality.
then a few sentence to explain algorithmic contribution. If this contribution is not strong. then one or two sentences suffice.
Somewhere in the middle of your reasoning, you can point to your table 1 to position our idea in the literature and convince readers that it is novel. ]}
\fi
	
This study considers a general three-state semi-Markov deteriorating system in which we allow general distributions for both transition times from healthy to defective and from defective to failure.
The associated mission abort unfolds as an optimal stopping problem, as formulated in Problem~\eqref{eq:stopping} ahead.	
Directly solving this problem is intractable due to the non-Markovian process dynamics.
For a tractable solution, we introduce a novel tool of Erlang mixtures to approximate non-exponential sojourn times therein.
% With a careful design as illustrated in Figure~\ref{fig:CTMC} ahead, 
The original three-state semi-Markov system can then be approximated with a surrogate continuous-time Markov chain (CTMC) 
%whose absorbing times can be linked to the two Erlang mixture distributions.  The number of states of the surrogate CTMC is larger than three in the original system, and it is 
whose number of states is larger than three in the original semi-Markov system and is determined by the number of components in the two Erlang mixture distributions.
By properly linking the surrogate CTMC to condition-monitoring data, a surrogate POMDP can be constructed.
%, where the state space is the same as that for the surrogate CTMC rather than the three health states in the original semi-Markov system.
The optimal mission abort policy derived from this surrogate POMDP serves as an approximation to optimal stopping of the original problem~\eqref{eq:stopping}.
Consolidating the above idea is not a simple task for reasons as follows.
\begin{itemize}
	\item
	In principle, 
	an infinite array of CTMCs can be devised, each with absorbing time conforming to the same Erlang mixture distribution.
	However, an arbitrary selection of two CTMCs for the two Erlang mixture distributions may make the concatenated CTMC complicated, posing challenges in theoretical analysis and numerical computations of the resulting surrogate POMDP.
	Therefore, the first crucial question to address is the initial design of the two CTMCs.
	\item 
	The surrogate POMDP approximates the original semi-Markov dynamics and retains the same mission abort cost structure. 
	When its solution is used to approximate the optimal solution of the original mission abort problem, i.e., the stopping time problem \eqref{eq:stopping} ahead, it is not immediately clear how the approximation works and what is the quality of the approximation since the two problems are defined in different probability spaces.
	\item 	
	If the optimal policy of the surrogate POMDP well approximates that of Problem~\eqref{eq:stopping}, analyzing its structure would provide deeper insights into the abort problem.
	For example, one may inquire: 
	(a) What is the structure of the optimal policy, and how does it differ from that of a standard CBM when the underlying system dynamics are as simple as a CTMC?
	(b) Does a similar structure persist when the sojourn times become non-exponential?
	(c) If not, can we characterize it? %, possibly by imposing some additional regularity conditions on the sojourn times? 
	The first question aids in understanding the fundamental disparity between mission abort and CBM, while the last two offer insights into the influence of non-exponential sojourn times on the optimal policy.
	However, POMDPs with more than two states are notoriously challenging to analyze \citep{zhang2022analytical}, while our surrogate POMDP is anticipated to have many states for an accurate approximation of the original system.
	\item
	With numerous states in the surrogate POMDP, directly employing standard point-based value iteration (PBVI) algorithms \citep{pineau2003point} for a numerical solution of the optimal policy proves to be exceedingly inefficient. 
	By leveraging structural properties of the optimal policy, it might be possible to tailor an effective algorithm for our mission abort setting.
\end{itemize}

This study aims to address the challenges outlined above. 
The contributions of this work can be summarized as follows.}
	
{Our first contribution lies in the development of a general modeling and solution framework for optimal mission abort with imperfect condition-monitoring data, an imperative consideration for a multitude of on-demand mission-critical systems.
We consider a general three-state semi-Markov deterioration process, resulting in an intractable optimal stopping problem in \eqref{eq:stopping}.
We address this intricacy by introducing an innovative Erlang mixture recipe to approximate the two non-exponential sojourn times of the semi-Markov process, leading to an CTMC approximation of the original semi-Markov process.
The surrogate CTMC is subsequently correlated to surrogate condition-monitoring signals,
{
which are generated by mirroring the signal generation process of the original system.
This leads to a surrogate POMDP for approximating \eqref{eq:stopping}, with the surrogate CTMC as system dynamics, mission abort decisions as actions, and surrogate condition-monitoring signals as observations.
We craft a construction for the CTMC approximation in Section~\ref{subsec:pomdp_phase} that endows this surrogate POMDP with favorable theoretical properties.}

%and makes it easy to optimize numerically. 
%\ye{[does it also facilitate numerical solution?]} 
%\sun{[the speed up of our PBVI algorithm is independent with the chain structure]}
		
Our second contribution stems from asymptotic analysis on the efficacy of the POMDP solution in approximating the optimal mission abort decision rules defined in Problem~\eqref{eq:stopping}.
To comprehend the approximation, we conceptualize the surrogate POMDP as a model defined in a probability space distinct from that of the original abort problem, while the surrogate signals retain the same sample space as the original signals.
As such, the POMDP optimal policy, akin to the optimal abort decision rules from solving the original problem~\eqref{eq:stopping}, is a map from the sample space for the condition-monitoring signals to the binary action space.
Our Theorem~\ref{thm:convergence} shows that the POMDP optimal policies, indexed by a sequence of Erlang rate parameters and treated as measurable functions in the original probability space, converge almost surely to the optimal abort decision rules defined by \eqref{eq:stopping} as the Erlang rate parameter goes to infinity.
% This result leverages the denseness property of Erlang mixtures in the field of positive distributions.
This result implies asymptotic optimality of the POMDP solution for the abort decision in that the expected cost by adopting the POMDP solution converges to the optimal expected cost in \eqref{eq:stopping} almost surely.
		
% Moreover, using mixture of Erlang distributions to approximate a known distribution is closely related to sieve estimations in statistics \citep{chen2007large}. When some distributions are unknown, our method is compatible with a data-driven framework that uses phase-type distributions as a sieve for parameter estimation. Section~\ref{sec:conclusion} discusses this future research in detail.

Our third contribution is to provide comprehensive structural results on the optimal policy of the surrogate POMDPs.
Under a general surrogate POMDP, we demonstrate in Theorem~\ref{thm:threshold_n} that there exists a threshold on the cumulative mission time, beyond which mission abort is always unnecessary. 
This finding offers solid support for the prevailing abort heuristic \citep{cha2018optimal, qiu2019gamma}.
Nevertheless, even when the system dynamics follow the simplest CTMC, our upcoming Corollary~\ref{thm:control_limit_CTMC} reveals that before reaching the time threshold, abort is optimal if and only if the belief state for defect falls in an interval whose upper endpoint is {\it smaller than one}.
This clearly reveals why the existing rule-based methods, which abort upon confirmation of a defect, cannot be optimal.
Moreover, the result is highly indicative of the distinction from CBM where a system is replaced upon confirmation of severe ageing.
% In addition, this result is inspiring in unveiling the structural properties for a general POMDP.
We further show in Theorem~\ref{thm:struct_policy} that the optimal abort policy can be delineated by an interval for the distance from the current belief state vector to the belief vector with unit mass on the worst state right before failure, with the two endpoints dependent on the angle of the current vector in an appropriately designed spherical coordinate system.
In this coordinate system, the region with abort as the optimal action can be shown to be convex.
This important property helps us to develop a modified PBVI algorithm. As demonstrated in Section~\ref{sec:num}, the runtime of our modified algorithm is less than 10\% of the standard PBVI algorithm, while maintaining the same numerical error.}

The remainder of the paper is organized as follows.
{Section~\ref{sec:review} reviews related literature.}
Section~\ref{sec:model} illustrates the mission abort problem. 
Section~\ref{sec:pomdp} introduces Erlang mixture distributions to approximate the system failure process as a CTMC and then formulates a POMDP framework.
Section~\ref{sec:property} derives structural properties of our model and {shows} the optimal abort policy follows a control-limit structure.
To deal with the curse of dimensionality,
Section~\ref{sec:algo} develops an algorithm that approximately solves the optimization problem.
Section~\ref{sec:special} further investigates two special cases of our model which can be exactly solved after discretizing the state space.
{Section~\ref{sec:multiple} extends to a multi-task setting.}
Section~\ref{sec:num} conducts a case study on a UAV, which have become increasingly prevalent in various domains, including asset management, healthcare and supply chain.
Section~\ref{sec:conclusion} concludes this study and discusses future research.
All technical proofs are provided in Appendix~\ref{appen:proof}.

\section{Literature Review}\label{sec:review}
{\cite{myers2009probability} pioneers the mission abort problem by basing the abort decision on the number of failed components, which oversimplifies the system dynamics.
Similar to this work, condition-monitoring signals are not considered in the early literature.	
Some subsequent studies base the abort decision on the number of shocks observed during the mission \citep{levitin2021optimal}.
In practice, however, shocks are hard to define and quantify, and magnitude of shocks often matters.
Both \cite{qiu2019gamma} and \cite{yang2019designing} use condition-monitoring signals to design the mission abort policy, but their models ignore the false and missed alarms.
This leads to suboptimal performances as condition-monitoring data are generally imperfect.
% For the optimal abort decision, we need to specify the relationship between the latent system health state and the observed signals.
We are not aware of an analytically optimal policy for mission abort under condition monitoring, even assuming Markovian system deterioration.
	
Condition-based mission abort shares similarities to CBM, both of which use imperfect signals to predict system failure and make decisions.
In the literature of CBM, Markov decision processes (MDPs) and POMDPs are frequently used.
\cite{elwany2011structured} and \cite{drent2023real} assume observable system degradation and use an MDP to integrate Bayesian learning and maintenance decision-making.
Control-limit policies with simple threshold functions are shown to be optimal in both studies.
When the system health state is unobservable, \cite{kim2013joint} use a POMDP to show the optimality of a two-level control-limit policy for inspection and replacement of a three-state system.
\cite{kim2016robust} considers parameter uncertainty and uses a robust POMDP to show the optimality of a control-limit policy for replacement of a two-state system.
Recently, \cite{zhang2022analytical} develop a dual framework for POMDPs based on \cite{zhang2010partially} to analytically derive the optimal inspection and replacement policy of a three-state system. 
The model introduces a novel perspective for maintenance problems by representing the optimal maintenance policy through six graphs. 
As discussed in Section~\ref{sec:intro}, most existing maintenance models assume a negligible maintenance time compared to the long planning horizon with linear yield in uptime.
In contrast, a distinctive aspect of mission abort problems is that aborting a mission incurs a non-negligible rescue time (see Section~\ref{subsec:model_decision}) and a zero-one mission failure cost.
Hence, existing analysis for maintenance models cannot be directly applied to our problem.
	
Most existing studies on reliability assume a Markovian deterioration system for tractability, which can be challenged in many applications as discussed in Section~\ref{sec:intro}. 
{Departing from this norm,
\cite{krishnamurthy2011bayesian} addresses a change point detection problem where the change time follows a phase-type distribution.
\cite{wang2015multistate} propose a Bayesian control chart that includes a phase-type distributed transition time from in-control to out-of-control states as a special case.
Both studies justify the use of phase-type distributions by highlighting their denseness property within the space of all positive-valued distributions, making them well-suited for approximating non-Markovian systems. 
\cite{khaleghei2021optimal} further explicitly construct a phase-type distribution to approximate a given semi-Markov process.
Based on the corresponding approximating Markov processes, a sequence of POMDPs is formulated for maintenance optimization.}
Their main approximation result, c.f.\ Proposition~2 therein, states that the optimal objective values of the surrogate POMDPs converge to that of the original system.
Their Theorem~3 establishes a control-limit structure for the optimal control policy of the POMDP, and their Figure~3 numerically shows that the control limits of the surrogate POMDPs, represented as a sequence of real numbers, converge to some limit which they regard as the optimal control limit for the original problem.
% Without formal proof, their Figure~3 numerically shows that the control limits of the surrogate POMDPs, represented as a sequence of real numbers, converge to some limit which is regarded as the optimal control limit for the original problem.

In contrast to \cite{khaleghei2021optimal}, the optimal policy for our original mission abort problem under a partially observed semi-Markov system can be highly unstructured.
We approximate the optimal polity, which is a set of functions from the signal space to the binary action space, using the POMDP optimal policy.
Convergence for this functional approximation necessitates a meticulous definition of probability spaces for various stochastic systems and is technically challenging.
In comparison, \cite{khaleghei2021optimal} do not prove convergence of optimal POMDP control limits to that of the original problem, which is convergence of a sequence of real numbers.
Furthermore, they only accommodate a general distribution for the transition time from the defective to failure state.
Our proposed model is more general by allowing the transition time from the healthy to defective state to be general distributed as well.
This generalization leads to a non-trivial and interesting construction of CTMCs as shown in Section~\ref{subsec:pomdp_phase} ahead.
	
To solve a POMDP for optimal stopping, standard methods such as value iteration \citep{piri2022individualized} or continuation-value iteration \citep{zhang2010partially} are effective when the state space is small.
Due to the curse of dimensionality, however, it is computationally infeasible to exactly solve a general large-scale POMDP.
Acceleration algorithms can be devised for solving a POMDP when the optimal policy exhibits certain structural properties.
For example, \cite{khaleghei2021optimal} show that the optimal policy of their surrogate POMDP is defined by a simple control limit on the system reliability;
solving this POMDP thus simplifies to a one-dimensional search for the optimal control limit.
In contrast, the optimal policy of our surrogate POMDP for mission abort lacks such a simple structure.
Popular methods to tackle a general POMDP include online planning methods \citep{silver2010monte} and offline iterative methods such as grid-based methods \citep{sandikcci2013alleviating} and PBVI algorithms \citep{pineau2003point}.
These general-purpose algorithms do not effectively scale to POMDPs with many states.
To attain an accurate solution within a reasonable timeframe, it is desirable to expedite the numerical procedure by capitalizing on structural insights derived from the mission abort problems.}

\section{Problem Statement}\label{sec:model}
Consider a mission-critical system that is required to operate continuously for a constant mission time $H>0$. 
This section formulates the mission abort decision problem given in-situ signals indicating the health condition of the system.
Section~\ref{subsec:model_deg} models the system failure process.
Section~\ref{subsec:model_obs} introduces the observable signals used to infer the system health state.
Section~\ref{subsec:model_decision} describes the decision-making problem.

\subsection{System Failure Process}\label{subsec:model_deg} 
During a mission, suppose the system failure process can be modeled by a stochastic process $\{X(t),~t\ge 0\}$ with a state space $\mathbb{S}=\{1,2,3\}$ defined on a probability space $(\Omega,\mathcal{F},\mathbb{P})$,
where $X(t)=1$, $2$, and $3$ indicate healthy, defective, and failure system states at time $t$, respectively. 
The failure state is self-announcing, but we cannot exactly tell whether the system is healthy ($X(t)=1$) or defective ($X(t)=2$) given that it is working at $t$.
Such a three-state model is commonly used in the existing literature to model failure processes \citep{panagiotidou2010statistical,abbou2019group}
and has been widely validated by real diagnostic data \citep{makis2006application}.
We stipulate $X(0)=1$, meaning a healthy state at the onset.
Based on degradation physics, we require $X(t)$ to be nondecreasing in $t$ \citep{panagiotidou2010statistical}.
As such, the system can either fail directly from a healthy state or go through a two-phase failure process, namely, from healthy to defective and then to the failure state.
Notationally, {let $T_{1j}$ be the underlying time to transition directly from state $1$ to state $j\in\{2,3\}$.}
If $T_{12}\geq T_{13}$, the system fails at time $T_{13}$;
otherwise, the system enters the defective state $2$ at time $T_{12}$, and the sojourn time at this state before transitioning to failure is $T_{23}$.
Based on the above assumptions, the system failure time $\xi\triangleq \inf\{t\geq 0: X(t)=3\}$ is
\begin{equation*}
	\xi=
	\begin{cases}
		T_{13}, & T_{13}\leq T_{12}, \\
		T_{12}+T_{23}, & \text{otherwise}.
	\end{cases}
\end{equation*}

We model $T_{13}$ as an exponential random variable with rate $\zeta$, while $T_{12},T_{23}$ follow a general continuous distribution with probability density functions (PDFs) $g(\cdot),f(\cdot)$ and cumulative distribution functions (CDFs) $G(\cdot),F(\cdot)$, respectively.
Furthermore, we model $T_{12},T_{13},T_{23}$ as independent random variables.
The exponential assumption on $T_{13}$ is a reasonable premise given that the mission time $H$ is typically substantially shorter than the system lifetime.
The failure rate of $T_{13}$ approximates the failure rate of the system at its present age, a value that remains relatively constant over a short mission duration for a long-lived system \citep{khojandi2018dynamic, abbou2019group, tian2021optimal}.
Similarly, an exponential $T_{12}$ is also appropriate, and the manifestation of the defect over the short mission duration should be non-informative on the progression of the defect to eventual failure.
Here, we opt for a general distribution for $T_{12}$ as it does not introduce fundamental complexity into the methodology developed subsequently.
In contrast, a generic distribution for the duration $T_{23}$ from defective to failure states is based on the consideration that $T_{23}$ is usually short and exhibits increasing failure rate when the system is defective.
{An increasing failure rate is a common assumption in reliability modeling, e.g., for production systems \citep{iravani2002integrated} and commercial modular aero-propulsion systems \citep{zhu2021robust}, and has been validated by many real-world datasets \citep{meeker2022statistical}.}
Letting $h(\cdot)$ be the hazard rate of $T_{23}$, this assumption is formally stated as follows.
\begin{assumption}\label{assump:nondecreasing}
The hazard rate $t\mapsto h(t)$ for $T_{23}$ is nondecreasing, and the failure rate at the defective state is larger than that at the healthy state, i.e., {$h(0)>\zeta$}.
\end{assumption}
{Section~\ref{sec:property} analyzes the optimal abort policies with and without Assumption~\ref{assump:nondecreasing}.}

\subsection{Observed Signals}\label{subsec:model_obs}
If the system is functioning at time $t$, it can be either healthy ($X(t)=1$) or defective ($X(t)=2$) for $t<\xi$.
To infer the underlying health state, we consider condition-monitoring signals collected with an inter-sampling interval $\delta$.
The interval $\delta$ is exogenously determined by sensor properties and is typically short in practice that entails real-time abort decision-making.
The condition-monitoring signal $Y_n$ at time $n\delta$ takes values in $[K]\triangleq\{1,\ldots,K\}$ for $n\in[N]\triangleq\{1,\ldots,N\}$, where we assume $N\triangleq H/\delta$ is an integer.
This assumption is not restrictive because if the condition-monitoring data are multi-dimensional, many state-of-the-art prognostic algorithms have been developed to fuse them into a discrete categorical signal \citep{hong2018big,liu2022machine}.
% , e.g., green and red light signals when $K=2$. Therefore, it is not restrictive to assume $Y_n\in[K]$. As in \cite{zhu2021robust}, fusing multisensor signals into a scalar value can be regarded as a variance reduction approach and also significantly facilitates optimization.
For completeness, we let $Y_n=0$ with probability one if the system has failed at time $n\delta$ for $n\in[N]$.

Letting $X_n \triangleq X(n\delta)$, we model its stochastic dependence with $Y_n$
as $\mathbb{P}(Y_n=k\mid X_n=i)=d_{ik}$ for $i\in\{1,2\}$ and $k\in[K]$.
The state-observation matrix is then denoted by $\mathbf{D}=(d_{ik})_{2\times K}$. 
In practice, this matrix can be estimated based on historical experimental data \citep{piri2022individualized}. 
The following assumption is needed when we establish the monotonicity of value functions in Section~\ref{sec:property}.
\begin{assumption}\label{assump:TP2}
	The matrix $\mathbf{D}$ is totally positive of order $2$ (TP2), i.e., $d_{ik}d_{i'k'}\geq d_{ik'}d_{i'k}$ for all $i\leq i'$ and $k\leq k'$.
\end{assumption}
In words, Assumption~\ref{assump:TP2} requires that a larger $Y_n$ indicates the system is more likely to be at the defective state.
Similar assumptions are commonly imposed on the state-observation matrix on disease severity in the healthcare literature \citep{sandikcci2013alleviating}.

\subsection{Mission Abort Decision Making}\label{subsec:model_decision}
Upon receiving a new signal $Y_n$ at time $n\delta$, 
our available data are $Y_{1:n}\triangleq (Y_1,\ldots,Y_n)$, and we have two action choices:
(i) to abort the mission and initiate a rescue procedure, 
or (ii) to continue the mission until observing $Y_{n+1}$.
The system can fail during either rescue or mission, incurring a \textit{system} failure cost $C_{\text{s}}$.
Moreover, a \textit{mission} failure cost $C_{\text{m}}$ is incurred if the mission is not completed.
When the mission is aborted at time $n\delta$, we let $w_n$, $n=0,\ldots,N-1$, be the rescue time with $w_0=0$.
When the mission is finished at time $H\triangleq N\delta$, we assume there is an additional time $w_N\geq 0$ to stop the system.
For example, a UAV takes time $w_N$ to fly back to the base upon finishing an inspection mission.
% On the other hand, when a paper mill completes the manufacturing mission, we can immediately stop the system so that $w_N=0$.
The following assumptions on the rescue procedure are in force.
\begin{assumption}\label{assump:rescue}
	(i) The rescue time $w_n$ is nondecreasing with $n$.
	(ii) The rescue procedure does not affect the failure process $X(t)$.
\end{assumption}
Assumption~\ref{assump:rescue} is common in literature \citep{yang2019designing,levitin2021optimal}.
We illustrate it using the UAV example again.
Consider a UAV that leaves the ground base to conduct a mission at a certain altitude. 
During the mission, the altitude of the UAV increases over time before it reaches the required altitude. 
The monotone assumption on $w_n$ is valid since a UAV at a higher altitude needs more time to land on the ground base.
Moreover, the rescue procedure does not change the system failure rate so the failure process $X(t)$ is not affected.
As a consequence of Assumption~\ref{assump:rescue}(i), 
we have $w_n\leq H-n\delta+w_N$ for all $n=0,\ldots,N-1$.
This inequality indicates that the rescue time at time $n\delta$ is shorter than the time to complete the mission and to stop the system;
otherwise, there is no need to abort the mission.

Based on the above setting, our objective is to find the optimal mission abort policy during the mission time $[0,H]$ that minimizes the expected total loss for a mission.

\section{POMDP Formulation}\label{sec:pomdp}
%Let $\mathcal{F}_n\triangleq \sigma(Y_{1:n},\mathbbm{1}{\{\xi\geq n\delta\}})$ be the $\sigma$-field generated by the condition-monitoring signals of a working system at time $n\delta$.
%Let $\mathcal{F}\triangleq(\mathcal{F}_n)_{n\in[N]}$ be the natural filtration.
Our decision-making problem is to determine a set of decision rules $a_n:\{0,\ldots,K\}^n\mapsto\{0,1\}$ for $n\in[N-1]$, such that given $Y_{1:n}=y_{1:n}\in[K]^n$, the action $a_n(y_{1:n})=1$ means aborting the mission at time $n\delta$ and $a_n(y_{1:n})=0$ means doing nothing until the next decision epoch.
On the other hand, $Y_n=0$ indicates the system has failed at time $n\delta$, and we let $a_n(Y_{1:n})=1$ if $Y_n=0$ for completeness.
%denote $\mathcal{T}(\mathcal{F})$ as the collection of all the $\mathcal{F}$-stopping times.
%Any feasible abort policy corresponds to a stopping time for mission abort, which has support $\{n \delta:n\in[N]\}$ and is denoted by $\mathfrak{T}\in\mathcal{T}(\mathcal{F})$.
Given $(a_n)_{n\in[N-1]}$, we denote $\mathfrak{T}\triangleq \min \{n\delta: a_n(Y_{1:n})=1,~n\in[N-1]\}\cup \{N\delta\}$ as the random time when the mission is aborted ($\mathfrak{T}<H$) or finished ($\mathfrak{T}=H$).
Let $C(X(\mathfrak{T}+w_{\mathfrak{T}/\delta}),\mathfrak{T})$ be the total operational cost for the mission, which is a function of the random time $\mathfrak{T}$ and the state $X(\mathfrak{T}+w_{\mathfrak{T}/\delta})$ when the system is fully stopped. 
Then $C(X(\mathfrak{T}+w_{\mathfrak{T}/\delta}),\mathfrak{T})$ is given by
$$
C(X(\mathfrak{T}+w_{\mathfrak{T}/\delta}),\mathfrak{T})
=(C_\text{m}+C_\text{s})\mathbbm{1}{\{X(\mathfrak{T}+w_{\mathfrak{T}/\delta})=3\}}
+C_\text{m}\mathbbm{1}{\{X(\mathfrak{T}+w_{\mathfrak{T}/\delta})<3,~\mathfrak{T}<H\}},
$$
where $\mathbbm{1}\{\cdot\}$ is the indicator function.
Our abort decision-making problem can be written as:
%\bfred{[Do we need to add the constraint $a_n(y_{1:n})\leq a_n(y_{1:n+1})$?]}
\begin{equation}\label{eq:stopping}
	\begin{aligned}
		\min_{\{a_n: n\in[N-1]\}} & \mathbb{E}[C(X(\mathfrak{T}+w_{\mathfrak{T}/\delta}),\mathfrak{T})] \\
		\text{s.t.~} \quad & \mathfrak{T}= \min\{n\delta: a_n(Y_{1:n})=1,~n\in[N-1]\}\cup \{N\delta\}; \\
		& a_n(y_{1:n})=1 \text{ if } y_n=0,~ n\in[N-1],
	\end{aligned}
\end{equation}
where the expectation $\mathbb{E}$ is taken with respect to the stochastic process $\{X(t),~0\leq t\leq H+w_N\}$ and condition-monitoring signals $Y_{1:N}$.
Let $\{a_n^*:n\in[N-1]\}$ be the set of optimal actions, and $\mathfrak{T}^*\triangleq \min \{n\delta: a_n^*(Y_{1:n})=1,~n\in[N-1]\}\cup \{N\delta\}$ be the optimal stopping time.
We assume the set of optimal actions is unique.
In case there are multiple sets of optimal actions, we stipulate the set with the smallest $\mathfrak{T}^*$ as the optimal solution.

\subsection{Erlang Approximation}\label{subsec:pomdp_phase}

Compared with most existing optimal stopping problems with Markovian system dynamics \citep{abbou2019group,tian2021optimal},
{the non-exponential sojourn times $T_{12}$ from healthy to defect and $T_{23}$ from defect to failure make Problem~\eqref{eq:stopping} challenging to solve.}
We approximate them using two sequences of Erlang mixture distributions.
This approximation will induce a sequence of surrogate POMDPs whose optimal policies approximate that for Problem~\eqref{eq:stopping}.

We start with approximating the distribution of the nonnegative $T_{23}$ by an Erlang mixture (also known as the hyper-Erlang distribution).
Let $E_m^{(\lambda)}(\cdot)$ be the CDF of an Erlang distribution with rate $\lambda>0$ and shape $m\in\mathbb{N}_+$.
% It is the distribution of the absorption time of a $\ell$-state Markov chain where t 
%\cite{schassberger1973warteschlangen} show that the CDF $F(\cdot)$ of any positive random variable can be approximated by a mixture of Erlang distributions given by
%\begin{equation*}
%F^{(\lambda)}^\infty(t)=\sum_{i=0}^{\infty}\bigg[F\left(\dfrac{i}{\lambda}\right)-F\left(\dfrac{i-1}{\lambda}\right)\bigg]E_\lambda^{(i)}(t),\quad t\ge0,
%\end{equation*} 
Corollary 4.1 of \citet{wiens1979distributions} shows that for any integer-valued function $m_2\triangleq m_2(\lambda)$ of $\lambda$ satisfying $\lim_{\lambda\to\infty}{m_2}/{\lambda} \to \infty$, the random variable $T_{23}^{(\lambda)}\geq 0$ with CDF 
\begin{equation}\label{eq:herlang}
	F^{(\lambda)}(t)=\sum_{i=1}^{m_2-1}\left[F\left(\dfrac{i}{\lambda}\right)-F\left(\dfrac{i-1}{\lambda}\right)\right]E_i^{(\lambda)}(t)+\left[1-F\left(\frac{m_2-1}{\lambda}\right)\right]E^{(\lambda)}_{m_2}(t),\quad t\ge0,
\end{equation}
converges weakly to $T_{23}$ as $\lambda\to\infty$.
In view of this property, we approximate the distribution of $T_{23}$ by that of $T_{23}^{(\lambda)}$ for some fixed $\lambda>0$ and $m_2\triangleq m_2(\lambda)>0$.
{
We choose $\lambda>\zeta$ and let
\begin{equation}\label{eq:herlang2}
	G^{(\lambda)}(t)=\sum_{i=1}^{m_1-1}\left[G\left(\dfrac{i}{\lambda-\zeta}\right)-G\left(\dfrac{i-1}{\lambda-\zeta}\right)\right]E_i^{(\lambda-\zeta)}(t)+\left[1-G\left(\frac{m_1-1}{\lambda-\zeta}\right)\right]E^{(\lambda-\zeta)}_{m_1}(t),\quad t\ge0,
\end{equation}	
for some $m_1\triangleq m_1(\lambda)$ satisfying $\lim_{\lambda\to\infty}m_1/\lambda\to\infty$.
Then $T_{12}^{(\lambda)}\sim G^{(\lambda)}$ converges weakly to $T_{12}$ as $\lambda\to\infty$ and is used to approximate $T_{12}$.
The choice of $\lambda-\zeta$ in \eqref{eq:herlang2} is a pivotal maneuver.
It allows us to construct a surrogate CTMC with identical transition rates out of any transient state, and the resulting surrogate POMDP has nice structural properties, as shown in Section~\ref{sec:property}.

\begin{figure}[]
	\centering
	\subfigure[Start from state $i$ with probability $\pi_i$]{
		\begin{tikzpicture}[scale=1.8, auto,swap]
			\node [circle, draw] (zero) {1};
			\node [circle, draw] (one) [right of=zero, xshift=0.22cm] {2};
			\node [circle, draw] (two) [right of=one, xshift=0.22cm] {3};
			\node [text width=0.4cm] (dots) [right of=two, xshift=0.22cm] {$\cdots$};
			\node [circle, draw] (last) [right of=dots,xshift=0.22cm] {$m$};
			\node [circle, draw] (absorb) [right of=last,xshift=0.22cm] {$A$};
			\node [text width=0.44cm] (pi1) [above of=zero] {$\pi_{1}$};
			\node [text width=0.44cm] (pi2) [above of=one] {$\pi_{2}$};
			\node [text width=0.44cm] (pi3) [above of=two] {$\pi_{3}$};
			\node [text width=0.44cm] (pi-last) [above of=last] {$\pi_{m}$};
			\path[->] (zero) edge (one);
			\path[->] (one) edge (two);
			\path[->] (two) edge (dots);
			\path[->] (dots) edge (last);
			\path[->] (last) edge (absorb);
			\path[->] (pi1) edge (zero);
			\path[->] (pi2) edge (one);
			\path[->] (pi3) edge (two);
			\path[->] (pi-last) edge (last);
		\end{tikzpicture}}
	\hspace{0.66cm}
	\subfigure[Start from state $1$ with probability $1$]{
		\begin{tikzpicture}[scale=1.8, auto,swap]
			\node [circle, draw] (zero) {1};
			\node [circle, draw] (one) [right of=zero, xshift=0.22cm] {2};
			\node [circle, draw] (two) [right of=one, xshift=0.22cm] {3};
			\node [text width=0.4cm] (dots) [right of=two, xshift=0.22cm] {$\cdots$};
			\node [circle, draw] (last) [right of=dots,xshift=0.22cm] {$m$};
			\node [circle, draw] (absorb) [right of=last,xshift=0.22cm] {$A$};
			\node [text width=0.44cm] (pi1) [above of=zero] {w.p.~$1$};
			\path[->] (zero) edge (one);
			\path[->] (one) edge (two);
			\path[->] (two) edge (dots);
			\path[->] (dots) edge (last);
			\path[->] (last) edge (absorb);
			\path[->] (pi1) edge (zero);
			\path[->] (zero) edge [bend left] (absorb);
			\path[->] (one) edge [bend left] (absorb);
			\path[->] (two) edge [bend left] (absorb);
		\end{tikzpicture}}
	\caption{Two CTMC schemes with $m$ transient states and an absorbing state $A$.
	With proper choice of the initial probabilities $(\pi_i)_{i\in[m]}$ and the transition rates, their absorbing times can follow the same distribution.
	% corresponding to either state $m_1+1$ or $m_1+m_2+1$ in our mission abort setting.
	The two schemes are used to approximate the transitions from $\mathcal{X}_1$ to $\mathcal{X}_2$ and from $\mathcal{X}_2$ to $\mathcal{X}_3$, respectively.
	%\ye{[could you add a circle with number $m$ to denote $m$ transient states? This will echo the Erlang mixture with $m$ components as I have added in the main text.]}
	%\ye{
	%[Deep reason: our observe data include two part: $Y$ and system is still functioning.
	%The 2nd piece of information rules out the transition from state $m_1+i$ to $m_1+m_2+1$.
	%But state $m_1+1$ is unobservable.] 
	%}
}\label{fig:CTMC}
\end{figure}
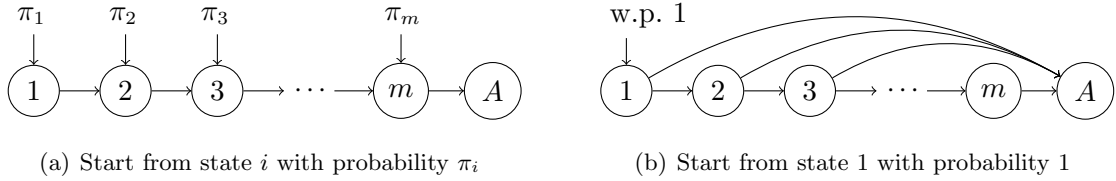

After the Erlang mixture approximations,}
two CTMCs with $m_1$ and $m_2$ transient states can be constructed with absorbing times following $G^{(\lambda)}$ and $F^{(\lambda)}$, respectively.
Concatenation of the two CTMCs yields a surrogate CTMC with $m_1+m_2$ transient states and one absorbing state (labeled as $m_1+m_2+1$) that approximates the original non-Markovian dynamics $X(\cdot)$.
However, the CTMC construction is non-unique \citep{buchholz2014input}.
{For illustration, Figure~\ref{fig:CTMC} shows two simple schemes to construct a CTMC based on an Erlang mixture distribution with $m$ components.
In Figure~\ref{fig:CTMC}(a), the CTMC can start from any transient state, after which the state transition sequence is deterministic. 
In Figure~\ref{fig:CTMC}(b), the CTMC starts from the first state. 
In each transition, it moves either to the next transient state or directly to the absorbing state.
In both CTMCs, we can select appropriate initial probabilities or transition rates to make their absorbing times follow a given $m$-component mixture of Erlang distribution. 
In principle, there are infinitely many CTMC schemes for $G^{(\lambda)}$ and $F^{(\lambda)}$ in \eqref{eq:herlang} and \eqref{eq:herlang2}.
It is extremely difficult, if not unattainable, to establish structural results for the resulting surrogate model if we casually choose two CTMCs.

After thorough experimentation, we find that using the scheme in Figure~\ref{fig:CTMC}(a) for $G^{(\lambda)}$ and that in Figure~\ref{fig:CTMC}(b) for $F^{(\lambda)}$ notably simplifies the ensuing theoretical analysis. 
The rationale behind this choice is as follows. 
In Figure~\ref{fig:CTMC}(a) for $G^{(\lambda)}$, we obtain a simple CTMC with deterministic transitions at a cost of a random start captured by an initial distribution. 
Its absorbing state $A$ corresponds to the collection of states $m_1+1,m_2+2,\ldots,m_1+m_2$ in the concatenated surrogate CTMC.
This explains the importance of starting the second chain, i.e., Figure~\ref{fig:CTMC}(b) for $F^{(\lambda)}$, at the first state that corresponds to state $m_1+1$ in the surrogate CTMC after concatenation.
As such, state $m_1$ transitions to state $m_1+1$ with probability one in the surrogate CTMC.
Initiating the second CTMC deterministically comes with a cost of a complicated transition rate matrix.
Nonetheless, this complexity does not pose an obstacle for our mission abort problem, because the absorbing state in Figure~\ref{fig:CTMC}(b) for $F^{(\lambda)}$ corresponds to the self-revealed failure state $m_1+m_2+1$ in the concatenated surrogate CTMC.
The acquisition of condition-monitoring data signifies the system's operational status, conditional on which we can logically exclude any transition from states $i\in[m_1+m_2]$ to the failure state $m_1+m_2+1$ in the CTMC analysis.
With the above construction, any transient state in our constructed surrogate CTMC can only transition to the next state when condition-monitoring data are received. 
This greatly simplifies the Kolmogorov backward equation for the surrogate CTMC.}

With Figure~\ref{fig:CTMC}(a) for $G^{(\lambda)}$ and Figure~\ref{fig:CTMC}(b) for $F^{(\lambda)}$, we construct a surrogate CTMC $\{X^{(\lambda)}(t):t\geq 0\}$ and the associated condition-monitoring signals $\{Y^{(\lambda)}_n:n\in[N]\}$ on some probability space $(\Omega^{(\lambda)},\mathcal{F}^{(\lambda)},P^{(\lambda)})$, which is different from $(\Omega,\mathcal{F},\mathbb P)$ for $\{X(t),~t\geq 0\}$ and $Y_{1:N}$.
Let the state space of ${X}^{(\lambda)}(t)$ be $\mathcal{X}\triangleq\{1,\ldots,m_1+m_2+1\}$ with state $m_1+m_2+1$ absorbing and others transient.
The states in $\mathcal{X}_1\triangleq\{1,\ldots,m_1\}$ correspond to the healthy state $X(t)=1$, states in $\mathcal{X}_2\triangleq\{m_1+1,\ldots,m_1+m_2\}$ correspond to the defective state $X(t)=2$, and the  state in $\mathcal{X}_3\triangleq\{m_1+m_2+1\}$ corresponds to the failure state $X(t)=3$.
We let $P^{(\lambda)}({X}^{(\lambda)}(0)=i)=\pi_{0i}^{(\lambda)}$, where
\begin{align}\label{eq:init_prob}
\begin{aligned}
\pi_{0i}^{(\lambda)}=
\begin{cases}
	1-G((m_1-1)/(\lambda-\zeta)), & i=1, \\
	G((m_1+1-i)/(\lambda-\zeta))-G((m_1-i)/(\lambda-\zeta)), & i=2,\ldots,m_1, \\
	0, & i=m_1+1,\ldots,m_1+m_2+1.
\end{cases}
\end{aligned}
\end{align}
Denote $q_{ij}\triangleq \lim_{\epsilon\downarrow 0}\frac{1}{\epsilon}P^{(\lambda)}(X^{(\lambda)}({t+\epsilon})=j\mid X^{(\lambda)}(t)=i)$, $i\neq j$, as the instantaneous transition rate from state $i\in\mathcal{X}$ to state $j\in\mathcal{X}$.
We construct the CTMC ${X}^{(\lambda)}(t)$ such that its instantaneous transition rate matrix $\mathbf{Q}=(q_{ij})_{i,j\in\mathcal{X}}$ satisfies
\begin{equation}\label{eq:trans_prob}
	{
		q_{ij}=
		\begin{cases}
			\lambda-\zeta, & j=i+1,~i\in\mathcal{X}_1, \\
			\zeta, & i\in\mathcal{X}_1,~ j=m_1+m_2+1,\\
			\lambda[1-p_{i-m_1}(\lambda)], & j=i+1,~i,j\in\mathcal{X}_2,\\
			\lambda p_{i-m_1}(\lambda), & i\in\mathcal{X}_2\backslash\{m_1+m_2\},~ j=m_1+m_2+1, \\
			\lambda, & i=m_1+m_2,~j=m_1+m_2+1,
	\end{cases}}
\end{equation}
$q_{ij}=0$ for other $i\neq j$, and $q_{ii}=-\sum_{j\in\mathcal{X}:j\neq i}q_{ij}$, where 
\begin{equation}\label{eq:jump_prob}
	p_i(\lambda)\triangleq\dfrac{F(i/\lambda)-F((i-1)/\lambda)}{1-F((i-1)/\lambda)}, \quad i\in [m_2-1].
\end{equation}
{The above construction leads to $q_{ii}=-\lambda$ for all $i\in[m_1+m_2]$, facilitating the subsequent analysis in Section~\ref{sec:property}.}
% \sun{[because $|q_{ii}|=\lambda-\zeta+\zeta=\lambda$ for $i\in\mathcal{X}_1$; $|q_{ii}|=[1-p_{i-m_1}]\lambda+p_{i-m_1}\lambda=\lambda$ for $i\in\mathcal{X}_2\backslash\{m_1+m_2\}$; $|q_{ii}|=\lambda$ when $i=m_1+m_2$]}
Consistent with $X(t)$, the states in $\mathcal{X}_1$ and $\mathcal{X}_2$ are unobservable, but the failure state $m_1+m_2+1$ of ${X}^{(\lambda)}(t)$ is observable.
We associate the surrogate signal $Y^{(\lambda)}_n$ with $X^{(\lambda)}(n\delta)$ as $P^{(\lambda)}(Y^{(\lambda)}_n=k \mid X^{(\lambda)}(n\delta)\in \mathcal{X}_i)=d_{ik}$ for $i=1,2$ and $k\in[K]$.
When $X^{(\lambda)}(n\delta)=m_1+m_2+1$, we let $Y^{(\lambda)}_n=0$ with $P^{(\lambda)}$-probability one.
%\bfred{[u need to be more formal in introducing the distribution of $Y^{(\lambda)}_n$]}

With the above construction, the surrogate CTMC ${X}^{(\lambda)}(\cdot)$ can be linked to the independent random variables {$T_{12}^{(\lambda)}$, $T_{13}$, and $T_{23}^{(\lambda)}$} as follows.

\begin{proposition}\label{prop:absorption}
	Consider the CTMC ${X}^{(\lambda)}(\cdot)$ with initial probabilities in \eqref{eq:init_prob} and  transition rate matrix $\mathbf{Q}$ in \eqref{eq:trans_prob}.
	Then the first hitting time of ${X}^{(\lambda)}(\cdot)$ to the states in $\mathcal{X}_2$ has the same distribution as {$T_{12}^{(\lambda)}$},
	and the absorption time of ${X}^{(\lambda)}(\cdot)$ has the same distribution as $T_{13}\mathbbm{1}{\{T_{13}\leq T_{12}^{(\lambda)}\}}+(T_{12}^{(\lambda)}+T_{23}^{(\lambda)})\mathbbm{1}{\{T_{13}> T_{12}^{(\lambda)}\}}$, where $T_{12}^{(\lambda)},T_{23}^{(\lambda)}$ and $T_{13}$ are independent random variables.
\end{proposition}
Recall the absorption time of $\{X(t):t\geq0\}$ is $\xi\triangleq T_{13}\mathbbm{1}{\{T_{13}\leq T_{12}\}}+(T_{12}+T_{23})\mathbbm{1}{\{T_{13}> T_{12}\}}$.
Since {$(T_{12}^{(\lambda)},T_{23}^{(\lambda)})\rightsquigarrow (T_{12},T_{23})$} as $\lambda\to\infty$, 
the above proposition, together with the continuous mapping theorem, implies that the absorption time of $\{X^{(\lambda)}(t):t\geq0\}$ converges weakly to the absorption time of $\{X(t):t\geq0\}$.

In the original problem, the decision at time $n\delta$ is based on the observed signal $Y_{1:n}$ with $Y_{1:n}\triangleq (Y_1,\ldots,Y_n)$, and the decision rule is a map $a_n:\{0,1,\ldots,K\}^n\mapsto \{0,1\}$.
We are interested in the sequence of the optimal decision rules $(a_n^*)_{n\in[N-1]}$ determined by Problem~\eqref{eq:stopping}, which is hard to solve.
In view of the distributional similarity in Proposition~\ref{prop:absorption},
we can consider a surrogate problem where the state-signal process $\{{X} (\cdot),Y_{1:N}\}$ in the original problem is replaced by $\{{X}^{(\lambda)}(\cdot),Y_{1:N}^{(\lambda)}\}$, while the cost structure remain unchanged.
The optimal decision rule $a_{\lambda,n}^{*}:\{0,1,\ldots,K\}^n\mapsto \{0,1\}$ for this surrogate problem is expected to be a good approximation of the map $a_n^*$, $n\in[N-1]$.
To consolidate this intuition, we can obtain $a_{\lambda,n}^*$ from solving the following surrogate problem:
\begin{equation}\label{eq:stopping_CTMC}
	\begin{aligned}
		\min_{\{a_{\lambda,n}: n\in[N-1]\}} & 
		\mathbb{E}^{(\lambda)}[C(X^{(\lambda)}(\mathfrak{T}+w_{\mathfrak{T}/\delta}),\mathfrak{T})] \\
		\text{s.t.~} \quad & \mathfrak{T}= \min\{n\delta: a_{\lambda,n}(Y_{1:n}^{(\lambda)})=1,~n\in[N]\}\cup\{N\delta\}, \\
		& a_{\lambda,n}(y_{1:n})=1, \quad \text{if } y_{n}=0,~ n\in[N-1],
	\end{aligned}
\end{equation}
where the expectation $\mathbb{E}^{(\lambda)}$ is taken {with} respect to $\{{X}^{(\lambda)}(t),~0\leq t\leq H+w_N\}$ and $Y^{(\lambda)}_{1:N}$.
Let $\mathfrak{T}^{(\lambda*)} \triangleq \min \{n\delta: a_{\lambda,n}^*(Y_{1:n}^{(\lambda)})=1,~n\in[N-1]\}\cup \{N\delta\}$
be the optimal stopping time of the surrogate problem.
Again, if there are multiple sets of optimal actions $\{a_{\lambda,n}^*:n\in[N-1]\}$, we stipulate the set with the smallest $\mathfrak{T}^{(\lambda*)}$ as the optimal solution.
This optimal solution is solved under the $P^{(\lambda)}$-law, which is used to approximate $\{a_{n}^*:n\in[N-1]\}$ obtained under the $\mathbb{P}$-law.
To assess its performance in the original problem, we define $\mathfrak{T}^{(\lambda)}\triangleq \min \{n\delta: a_{\lambda,n}^*(Y_{1:n})=1,~n\in[N-1]\}\cup \{N\delta\}$ as the stopping time using the approximate solution in the original problem.
The following theorem shows the accuracy of the approximation under the $\mathbb P$-law.
%\bfred{[i think it is about use of $A^{(\lambda)*}(t)$ to approximate $A^*(t)$.]}
\begin{theorem}\label{thm:convergence}
	Consider the sequence of optimal actions $\{a_{\lambda,n}^*:n\in[N-1]\}$ defined in the surrogate problem~\eqref{eq:stopping_CTMC}, and $\{a_{n}^*:n\in[N-1]\}$ defined in the original problem~\eqref{eq:stopping}. 
	For any $n\in[N-1]$, we have $a_{\lambda,n}^*(Y_{1:n})\to a_n^*(Y_{1:n})$ $\mathbb{P}$-almost surely as $\lambda\to\infty$. 
	As a result, $C(X(\mathfrak{T}^{(\lambda)}+w_{\mathfrak{T}^{(\lambda)}/\delta}),\mathfrak{T}^{(\lambda)})\rightarrow C(X(\mathfrak{T}^*+w_{\mathfrak{T}^*/\delta}),\mathfrak{T}^*)$ $\mathbb{P}$-almost surely as $\lambda\to\infty$, and
	\begin{equation}\label{eq:conv_exp_cost}
		\begin{aligned}
			\mathbb{E}[ C(X(\mathfrak{T}^*+w_{\mathfrak{T}^*/\delta}),\mathfrak{T}^*)]
			=&\lim_{\lambda\to\infty}\mathbb E[C(X(\mathfrak{T}^{(\lambda)}+w_{\mathfrak{T}^{(\lambda)}/\delta}),\mathfrak{T}^{(\lambda)})]\\
			=&{\lim_{\lambda\to\infty}\mathbb E^{(\lambda)}[C(X^{(\lambda)}(\mathfrak{T}^{(\lambda*)}+w_{\mathfrak{T}^{(\lambda*)}/\delta}),\mathfrak{T}^{(\lambda*)})]} .
		\end{aligned}
	\end{equation}
	%\begin{equation}\label{eq:convergence}
	%	\min_{\mathfrak{T}\in\mathcal{T}(\mathcal{F})}\mathbb{E}[C(X(\mathfrak{T}+w_{\mathfrak{T}/\delta}),\mathfrak{T})]=\lim_{\lambda\to\infty}\mathbb{E}[C(X(\mathfrak{T}^*_\lambda+w_{\mathfrak{T}^*_\lambda/\delta}),\mathfrak{T}^*_\lambda)].
	%\end{equation}
\end{theorem}

The first part of Theorem~\ref{thm:convergence} indicates that we can solve the surrogate problem~\eqref{eq:stopping_CTMC} with sufficiently large $\lambda$ and use the resulting sequence of decision rules $(a_{\lambda,n}^*)_{n\in[N-1]}$ for online decision-making in the original problem.
{The first equality of \eqref{eq:conv_exp_cost} reveals that this sequence of decision rules is asymptotically optimal in that the expected costs of the original mission abort problem, when applying these approximate policies, converge to the optimal expected costs of Problem~\eqref{eq:stopping}.
In comparison, Proposition~2 of \cite{khaleghei2021optimal} shows that the optimal objective values of their surrogate problems converge to that of the original semi-Markov maintenance problem, which corresponds to the second equality in \eqref{eq:conv_exp_cost}.
Such convergence result might not be interesting enough, as the two objective values are obtained from two different stochastic systems.
In contrast, our convergence results in Theorem~\ref{thm:convergence}, except the second equality in \eqref{eq:conv_exp_cost}, are under the $\mathbb{P}$-law for the original non-Markovian system.
The result directly shows the expected cost of applying the actions obtained from the surrogate problem to the original problem.
Comparing this cost with the optimal cost provides a straightforward assessment of solution quality.}

% We propose selecting finite $\lambda$ to ensure the system failure rate is nondecreasing over time. Specifically, 
The above theorem is an asymptotic result. 
To gain some insight on an appropriate choice of the $\lambda$ value, we look at the system failure rate when the CTMC $X^{(\lambda)}(\cdot)$ starts from the defective state, 
i.e.,
$h^{(\lambda)}(t)\triangleq\lim_{\Delta t\downarrow 0}{P}^{(\lambda)}(X^{(\lambda)}(t+\Delta t)=m_1+m_2+1\mid X^{(\lambda)}(0)=m_1+1,~X^{(\lambda)}(t)\neq m_1+m_2+1)/\Delta t$.
Under Assumption~\ref{assump:nondecreasing}, the time $T_{23}$ from defective to failure has a nondecreasing hazard rate $h(\cdot)$, and $h(0)>\zeta$.
If $T_{23}^{(\lambda)}$ well approximates $T_{23}$ in distribution, $h^{(\lambda)}(t)$ should also have this property.
% \bfred{[the failure rate is for either state 1 or 2 to transition to state 3. If you don't have a reasonable assumption on the failure rate of $T{12}$, intuitively, your failure rate result below might not be valid. Please double check]}
We formalize these statements as follows.

\begin{lemma}\label{lemma:increase_rate}
	Under Assumptions~\ref{assump:nondecreasing}, there exists a constant $\bar\lambda>0$ such that $h^{(\lambda)}(0)>\zeta$ and $h^{(\lambda)}(t)$ is nondecreasing in $t>0$ for any $\lambda>\bar\lambda$. 
\end{lemma}
In view of the above lemma, a necessary condition for good approximation of the surrogate problem is that $h^{(\lambda)}(t)$ is nondecreasing in $t\geq 0$ and $h^{(\lambda)}(0)>\zeta$ when Assumption~\ref{assump:nondecreasing} holds.
As such, we need to carefully monitor this property when selecting $\lambda$.
The properties in Lemma~\ref{lemma:increase_rate} will also be useful in Section~\ref{sec:property} when we establish structural properties and derive the optimal abort policy for the surrogate problem under Assumption~\ref{assump:nondecreasing}.
When Assumption~\ref{assump:nondecreasing} does not hold, $\lambda$ can be chosen through a visual comparison of the CDFs $F^{(\lambda)},G^{(\lambda)}$ and the original $F,G${, as illustrated in Appendix~\ref{appen:num_setting}}.
% Numerically results in Section~\ref{sec:num} show that moderate $\lambda$ has brought significant performance improvement compared with some commonly used benchmarks in abort decision-making.

\subsection{Bellman Recursion}\label{subsec:pomdp_model}

For fixed $\lambda>0$, {$m_1\triangleq m_1(\lambda)$}, and $m_2\triangleq m_2(\lambda)$, we now formulate the surrogate problem~\eqref{eq:stopping_CTMC} as a POMDP where $\mathcal{X}_1\cup\mathcal{X}_2$ are treated as hidden states.
%For notational convenience, we suppress the dependence of $Y_{1:n}^{(\lambda)}$ on $\lambda$ and write $Y_{1:n}$ when there is no confusion.
Let $p_{ij}(\delta)\triangleq P^{(\lambda)}({X}^{(\lambda)}(t+\delta)=j\mid {X}^{(\lambda)}(t)=i)$, $i,j\in\mathcal{X}$, be the transition probability from state $i$ to state $j$ during a sampling interval $\delta$, which can be computed by solving the Kolmogorov backward equations \citep{liu2021optimal}.
Then $\mathbf{P}(\delta)=(p_{ij}(\delta))_{i,j\in\mathcal{X}}$ is the one-step transition probability matrix.
For $n\in[N]$ and $i\in\mathcal{X}_1\cup\mathcal{X}_2$, define $\Pi_{ni}\triangleq P^{(\lambda)}({X}^{(\lambda)}(n\delta)=i\mid Y_{1:n}^{(\lambda)})$ as the conditional probability that the CTMC ${X}^{(\lambda)}(\cdot)$ is at hidden state $i$ at time $n\delta$ given the signals $Y_{1:n}^{(\lambda)}$.
After the value of $Y_{1:n}^{(\lambda)}$ is observed, the system \textit{belief state} at time $n\delta$ is $\bm\pi_n\triangleq [\pi_{n1},\ldots,\pi_{n,m_1+m_2}]'$, where $\pi_{ni}$ is a realization of $\Pi_{ni}$.
The belief state space is then given by the probability simplex $\mathcal{S}\triangleq\{\bm\pi=(\pi_i)_{i\in[m_1+m_2]}\in\mathbb{R}^{m_1+m_2}_+:\sum_{i=1}^{m_1+m_2}\pi_{i}=1\}$.
As with the existing POMDP literature (e.g., \citealp{liu2021optimal,wang2022optimal}), we introduce a value function $V^{(\lambda)}(n,\bm\pi_n)$, $n=0,\ldots,N$, %\ye{[$V^{(\lambda)}(n,\bm\pi_n)$?]}
that denotes the minimum cost-to-go at time $n \delta$ when the corresponding belief state is $\bm \pi_n \in \mathcal{S}$.
The Bellman equations are given by
\begin{equation}\label{eq:bellman_base}
	V^{(\lambda)}(n,\bm\pi_n) =\text{min}\{V_{\text{ab}}(n,\bm\pi_n),V_{\text{c}}(n,\bm\pi_n)\}, \quad n=0,\ldots,N-1, 
\end{equation}
where $V_{\text{ab}}(n,\bm\pi_n)$ and $V_{\text{c}}(n,\bm\pi_n)$ are the respective minimum cost-to-go when the action at time $n\delta$ is to \textit{abort} and to \textit{continue} the mission (we omit their dependence on $\lambda$ for notational simplicity).

We first derive the explicit expression of $V_{\text{ab}}(n,\bm\pi_n)$.
Assume that the system is working with a belief state $\bm\pi$.
Let $\kappa(t,\bm\pi)$ be the probability that the system's remaining useful life (RUL) is less than $t$, which is given by
% \ye{[is it better to use ${\bm p}_{\cdot,m_1+m_2+1}(t)$?]}
\begin{equation}\label{3}
	\kappa(t,\bm\pi)=\bm\pi' {\bm p}_{\cdot,m_1+m_2+1}(t)
\end{equation}
with ${\bm p}_{\cdot,m_1+m_2+1}(t)\triangleq [p_{1,m_1+m_2+1}(t),\ldots,p_{m_1+m_2,m_1+m_2+1}(t)]'$.
The expected cost of aborting the mission and initiating a rescue procedure at time $n\delta$ can then be expressed as 
\begin{equation}\label{eq:v_ab}
	V_{\text{ab}}(n,\bm\pi_n)=C_{\text{m}}+C_{\text{s}}\kappa(w_n,\bm\pi_n), \quad n=0,\ldots,N,
\end{equation}
where the first term is the mission failure cost due to abort and the second term is the expected cost of system failure during the rescue procedure.

Next, we compute the expected cost $V_{\text{c}}(n,\bm\pi_n)$, which is given by
\begin{equation}\label{eq:v_c}
	V_{\text{c}}(n,\bm\pi_n)=(C_{\text{s}}+C_{\text{m}})\kappa(\delta,\bm\pi_n)+\mathbb{E}^{(\lambda)}\left[V^{(\lambda)}(n+1,\bm\Pi_{n+1})\mid \bm\Pi_n=\bm\pi_n\right], \quad n=0,\ldots,N-1,
\end{equation}
with $\bm\Pi_n\triangleq(\Pi_{ni})_{i\in[m_1+m_2]}$.
The first term in \eqref{eq:v_c} is the expected cost for system and mission failure before the next decision epoch $(n+1)\delta$.
The second term is the cost-to-go at the $(n+1)$th decision epoch and can be evaluated as
\begin{equation}\label{eq:U_function}
	\mathbb{E}^{(\lambda)}\left[V^{(\lambda)}(n+1,\bm\Pi_{n+1})\mid \bm\Pi_n=\bm\pi_n\right]=\sum_{k=1}^K
	P^{(\lambda)}(Y_{n+1}^{(\lambda)}=k\mid\bm\Pi_n=\bm\pi_n)V^{(\lambda)}(n+1,\bm\pi_{n+1}(\bm\pi_n,k)), \nonumber
\end{equation}
where $\bm\pi_{n+1}(\bm\pi_n,k)$ denotes the belief state at time $(n+1)\delta$ if the belief state at time $n\delta$ is $\bm\pi_n\in\mathcal{S}$ and the observed signal at time $(n+1)\delta$ is $Y_{n+1}=k\in[K]$.
By Bayes' theorem, $\bm\pi_{n+1}(\bm\pi_n,k)$ is given by
\begin{equation}\label{eq:bayes_update}
	\bm\pi_{n+1}'(\bm\pi_n,k)
	=\dfrac{\bm \pi_n'\tilde{\mathbf{P}}(\delta)\mathrm{diag}(\tilde{\bm d}_k)}
	{\boldsymbol{\pi}_n'\tilde{\mathbf{P}}(\delta)\tilde{\boldsymbol{d}}_k},
\end{equation}
where $\tilde{\mathbf{P}}(\delta)$ is the square submatrix of $\boldsymbol{\mathbf P}(\delta)$ containing the first $m_1+m_2$ rows and columns of $\boldsymbol{\mathbf P}(\delta)$, $\tilde{\bm d}_k\triangleq[d_{1k}\bm 1_{m_1}',d_{2k}\bm 1_{m_2}']'$, and $\bm 1_r$ is the $r$-dimensional column vector of ones.
Moreover, we have $P^{(\lambda)}(Y_{n+1}^{(\lambda)}=k\mid\bm\Pi_n=\bm\pi_n)=\sum_{i=1}^{m_1+m_2}\sum_{j=1}^{m_1+m_2}\pi_{ni}p_{ij}(\delta)\tilde d_{jk}$.

Finally, the terminal cost of the Bellman equation is given by
\begin{equation}\label{eq:bellman}
	V^{(\lambda)}(N,\bm\pi_N)\triangleq V_\text{c}(N,\bm\pi_N)=(C_{\text{s}}+C_{\mathrm{m}})\kappa(w_N,\bm\pi_N), 
\end{equation}
which is the expected loss for stopping the system (e.g., let the UAV fly back to the base) immediately after it completes the mission at time $H$.

{Due to our CTMC construction, the system starts from state $i$ with $P^{(\lambda)}$-probability $\pi_{0i}^{(\lambda)}$.
The belief state at time $0$ is $\bm\pi_0^{(\lambda)}=(\pi_{0i}^{(\lambda)})_{i\in[m_1+m_2]}$.}
Based on the above finite-horizon POMDP framework, Problem~\eqref{eq:stopping_CTMC} is equivalent to finding the optimal mission abort policy that minimizes {$V^{(\lambda)}(0,\bm \pi_0^{(\lambda)})$}, where the policy at time $n\delta$ maps any belief state $\bm\pi_n\in\mathcal{S}$ to an action, i.e., to abort or continue the mission.
Since $\bm\pi_n$ encapsulates all historical information from $Y_{1:n}$ \citep{bertsekas2012dynamic,wang2022optimal}, basing abort on $\bm\pi_n$ is more reliable than on a single signal $Y_n$.
When the system is healthy, an unusual $Y_n$ signal cannot significantly change $\bm\pi_n$.
On the other hand, successive alarms would change the belief state towards $\mathcal{X}_2$ and trigger mission abort.
Our structural results in the next section formalize the above reasoning.
% Moreover, numerical experiments in Section~\ref{sec:num} show a substantial cost saving of making decisions using $\bm\pi_n\in\mathcal{S}$ compared with a two-step benchmark that makes abort decisions based on alarms.

\section{Structural Properties}\label{sec:property}

For notational simplicity, we drop the superscript in $V^{(\lambda)}(n,\bm\pi)$ and write it as $V(n,\bm\pi)$ hereafter.
{Different from most existing literature on POMDPs that assumes all the state are unobservable (e.g., \citealp{smallwood1973optimal}),}
the POMDP in Section~\ref{sec:pomdp} is based on a CTMC with both unobservable ($X^{(\lambda)}(t)\in\mathcal{X}_1\cup\mathcal{X}_2$) and observable ($X^{(\lambda)}(t)=m_1+m_2+1$) states.
Nevertheless, we can still show that the value functions of our POMDP exhibit the same piecewise linearity and concavity as those in POMDPs where all states are unobservable.

\begin{lemma}\label{lemma:piecewise_and_concave}
	% \ye{[please check if you need the condition $\lambda>\bar\lambda$, which essentially means two things: (a) you need to make assumption 1, and (b) Lemma 1 holds.]}	
	The value function $V(n,\bm\pi)$ is piecewise linear and concave in $\bm\pi$ for all $n=0,\ldots,N$.
	Moreover, both $V_{\mathrm{ab}}(n,\bm\pi)$ and $V_{\mathrm{c}}(n,\bm\pi)$ are concave in $\bm\pi$ for all $n=0,\ldots,N-1$.
\end{lemma}

To fully characterize the optimal abort policy, we need structural properties in addition to the convexity in Lemma~\ref{lemma:piecewise_and_concave}.
{We first study the optimal abort policy when Assumption~\ref{assump:nondecreasing} holds.}
In this case, we derive an upper bound of $V_\text{c}(n,\bm\pi_n)$ and compare it with $V_\text{ab}(n,\bm\pi_n)$ to obtain a sufficient condition under which to continue is the optimal action.
Let 
\begin{equation}\label{eq:upper_V}
	\overline{V}_\text{c}(n,\bm\pi_n)\triangleq (C_{\text{s}}+C_{\mathrm{m}})\kappa((N-n)\delta + w_N,\bm\pi_n), \quad n=0,\ldots,N,
\end{equation}
be the cost-to-go of continuing the mission until the end of the mission, while ${V}_\text{c}(n,\bm\pi_n)$ is the optimal cost-to-go if we continue the mission at the current decision epoch but use the optimal actions in the sequel. 
The cost-to-go of the latter case is obviously lower.
%\ye{[this result is too obvious. i am not sure if it deserves a lemma. i would suggest that we delete this part highlighted in red] The above statement is then formalized as follows. 
%\begin{lemma}\label{lemma:upper_bound}
%	The function $\overline{V}_\textup{c}(n,\bm\pi_n)$ defined in \eqref{eq:upper_V} satisfies $V_{\mathrm{c}}(n,\bm\pi_n)\le \overline{V}_{\mathrm{c}}(n,\bm\pi_n)$ for all $n=0,\ldots,N-1$ and $\bm\pi_n\in\mathcal{S}$. 
%\end{lemma}}

We want to compare the functions $\overline{V}_\text{c}(n,\bm\pi)$ and $V_\text{ab}(n,\bm\pi)$.
For any two vectors $\bm\pi^{(1)},\bm\pi^{(2)}\in\mathcal{S}$, let $\pi_{j}^{(i)}$ be the $j$th entry of $\bm\pi^{(i)}$, $i=1,2$.
Then $\bm\pi^{(1)}$ is greater than $\bm\pi^{(2)}$ with respect to the monotone likelihood ratio (MLR) ordering, denoted by $\bm\pi^{(1)}\succeq_\textnormal{LR}\bm\pi^{(2)}$, if $\pi_{i}^{(1)}/\pi_{j}^{(1)}\geq\pi_{i}^{(2)}/\pi_{j}^{(2)}$ for all $i>j$ \citep{shaked2007stochastic}.
A function $f:\mathcal{S}\mapsto\mathbb{R}$ is MLR nondecreasing (resp.\ nonincreasing) in $\bm\pi\in\mathcal{S}$ if $f(\bm\pi^{(1)})\geq f(\bm\pi^{(2)})$ (resp.\ $f(\bm\pi^{(1)})\leq f(\bm\pi^{(2)})$) for any $\bm\pi^{(1)}\succeq_\textnormal{LR}\bm\pi^{(2)}$.
We then show the monotonicity of $V_\textnormal{ab}(n,\bm\pi)$ and the upper bound $\overline{V}_\textnormal{c}(n,\bm\pi)$.

%\ye{[Lemma 4 is again very obvious. Since it shares the same conditions as lemma 5, i suggest combine Lemma 4 and 5.]}
\begin{lemma}\label{lemma:monotonous_pi_n}
	%\bfred{[check if you need any assumption/condition]}
	Suppose that Assumptions~\ref{assump:nondecreasing} and \ref{assump:rescue} hold and $\lambda>\bar\lambda$ {for the $\bar\lambda$ in Lemma~\ref{lemma:increase_rate}}. 
	{We have $p_{i,m_1+m_2+1}(t)$ is nondecreasing in $i\in[m_1+m_2+1]$ for any $t>0$.}
	Moreover, the cost-to-go $V_{\mathrm{ab}}(n,\bm\pi)$ for the abort decision is MLR nondecreasing in $\bm\pi\in\mathcal{S}$ and nondecreasing in $n$,
	and $\overline{V}_{\mathrm{c}}(n,\bm\pi)$ is MLR nondecreasing in $\bm\pi\in\mathcal{S}$ and nonincreasing in $n$. 
\end{lemma}

Based on Lemma~\ref{lemma:monotonous_pi_n}, we can derive the following results.
\begin{theorem}\label{thm:never_abort}
	%\bfred{[check if you need any assumption/condition]}	
	Suppose that Assumptions~\ref{assump:nondecreasing} and \ref{assump:rescue} hold and $\lambda>\bar\lambda$, and consider the surrogate problem~\eqref{eq:stopping_CTMC}.
	The optimal actions satisfy $a_{\lambda,n}^*(y_{1:n})=0$ for all $y_{1:n}\in[K]^n$ and $n=0,\ldots,N-1$ if \begin{equation}\label{eq:cost_ratio}
		\dfrac{ C_{\mathrm{m}}}{C_{\mathrm{s}}}
		\ge
		\exp[\lambda(H+w_N)]-1.
	\end{equation}
\end{theorem}
Intuitively, Theorem \ref{thm:never_abort} indicates that aborting a mission is never optimal during the planning horizon $[0,\,T]$ if the mission failure cost $C_\text{m}$ is significantly larger than the system failure cost $C_\text{s}$.
In this case, the priority of completing a mission is much higher than avoiding a system failure.
On the other hand, we derive the following result when \eqref{eq:cost_ratio} is violated.
\begin{theorem}\label{thm:threshold_n}
	%\bfred{[check if you need any assumption/condition]}	
	Suppose that Assumptions~\ref{assump:nondecreasing} and \ref{assump:rescue} hold and $\lambda>\bar\lambda$ {for the $\bar\lambda$ in Lemma~\ref{lemma:increase_rate}}.
	If the condition in \eqref{eq:cost_ratio} does not hold, let
	%\ye{[in your $\hat n$ def below, it is over $n=0,1,\cdots,N$. But after the equation, you say for $n=\hat{n},\ldots,N-1$. some thoughts need on the range]}
	\begin{equation}\label{eq:hat_n}
		\hat n\triangleq\min\left\{n=0,\ldots,N:\max_{\bm\pi\in\mathcal{S}}\left(\overline V_\textnormal{c}(n,\bm\pi)-V_\textnormal{ab}(n,\bm\pi)\right)\leq 0\right\}.
	\end{equation}
	If $\hat n<N$, then the optimal actions satisfy $a_{\lambda,n}^*(y_{1:n})=0$ for any $y_{1:n}\in[K]^n$ and $n=\hat{n},\ldots,N-1$.
	Moreover, the value function satisfies $V(n,\bm\pi)=\overline{V}_\textnormal{c}(n,\bm\pi)$ for $n=\hat n,\ldots,N-1$.
	% In addition, $\hat n= \min\left\{n=1,\ldots,N: \overline{V}_{\mathrm{c}}(n,\bm e_{m_1+m_2})\le V_{\mathrm{ab}}(n,\bm e_{m_1+m_2})\right\}.$
\end{theorem}
Theorem \ref{thm:threshold_n} suggests that continuing the mission is the optimal action after a time threshold $\hat n\delta$, intuitively because the mission is nearly finished after this time.
\cite{yang2019designing} show a similar result when the state $X(t)$ is observable, under which the abort decision is straightforward.
In comparison, the unobservable $X(t)$ makes it difficult to characterize the optimal abort policy for $n=0,\ldots,\hat n - 1$.
To help with the analysis, we provide the following structural properties of the value function $V(n,\bm\pi)$ in \eqref{eq:bellman_base} and the belief state vector $\bm\pi_{n+1}(\bm\pi,k)$ in \eqref{eq:bayes_update}.

\begin{lemma}\label{lemma:MLR_montonicity}
	Under Assumptions~\ref{assump:nondecreasing} and \ref{assump:TP2}, the vector $\bm\pi_{n+1}(\bm\pi,k)$ satisfies the following properties for all $n=0,\ldots,N-1$ if $\lambda>\bar\lambda$.
	\begin{itemize}
		\item For any $\bm\pi^{(1)},\bm\pi^{(2)}\in\mathcal{S}$ such that $\bm\pi^{(1)}\succeq_\textnormal{LR}\bm\pi^{(2)}$, we have $\bm\pi_{n+1}(\bm\pi^{(1)},k)\succeq_\textnormal{LR}\bm\pi_{n+1}(\bm\pi^{(2)},k)$ for all $k\in[K]$. 
		\item For any $k_1,k_2\in[K]$ such that $k_1\geq k_2$, we have $\bm\pi_{n+1}(\bm\pi,k_1)\succeq_\textnormal{LR}\bm\pi_{n+1}(\bm\pi,k_2)$ for all $\bm\pi\in\mathcal{S}$. 
	\end{itemize}
\end{lemma}

\begin{lemma}\label{lemma:monotone}
	Under Assumptions~\ref{assump:nondecreasing} and \ref{assump:TP2}, the value function $V(n,\bm\pi)$ is MLR nondecreasing in $\bm\pi\in\mathcal{S}$ for any fixed $n=1,\ldots,N$ if $\lambda>\bar\lambda$.
\end{lemma}

Lemma~\ref{lemma:MLR_montonicity} orders the belief states at the next decision epoch based on the current belief $\bm\pi\in\mathcal{S}$ and observed signal $k\in[K]$.
It is established by showing that 
%\ye
{the matrix $\tilde{\mathbf{P}}(\delta)$ is TP2, c.f.\ Lemma~\ref{lemma:TP2} in Appendix~\ref{proof:MLR montonicity}.}
Lemma~\ref{lemma:monotone} indicates that the cost-to-go is lower when the system is in a ``healthier'' state (in the sense of the MLR ordering).
Let $\bm e_i$ denotes the unit vector with $1$ in the $i$th component.
We need an additional structural property to fully characterize the optimal abort policy. 

\begin{lemma}\label{lemma:ab_opt}
	%If $V_\textnormal{c}(n,\bm e_{m_1+m_2})>V_\textnormal{ab}(n,\bm e_{m_1+m_2})$, then $V_\textnormal{c}(n-1,\bm e_{m_1+m_2})>V_\textnormal{ab}(n-1,\bm e_{m_1+m_2})$.
	Assume Assumption~\ref{assump:rescue} holds.
	If the program $\max\{n=0,1,\ldots,N-1:V_\textnormal{c}(n,\bm e_{m_1+m_2})>V_\textnormal{ab}(n,\bm e_{m_1+m_2})\}$ has an optimal feasible solution denoted as $\tilde{n}$, then $V_\textnormal{c}(n,\bm e_{m_1+m_2})>V_\textnormal{ab}(n,\bm e_{m_1+m_2})$ for all $n\leq\tilde n$.
	%\ye{[is it possible to add the following conclusion?: In addition, if $V_\textnormal{c}(\check{n},\bm e_{m_1+m_2})>V_\textnormal{ab}(\check{n},\bm e_{m_1+m_2})$ for some $\check n\in [n-1]$, then $V_\textnormal{c}(n,\bm e_{m_1+m_2})>V_\textnormal{ab}(n,\bm e_{m_1+m_2})$ for all $n\leq \check{n}$. Is the above is possible, then we can define $\tilde n$ as the largest $n$ st $V_\textnormal{c}(n,\bm e_{m_1+m_2})>V_\textnormal{ab}(n,\bm e_{m_1+m_2})$ ]}
\end{lemma}
In the sense of MLR ordering, the smallest and largest vectors in $\mathcal{S}$ are $\bm e_1$ and $\bm e_{m_1+m_2}$, respectively.
Lemma~\ref{lemma:ab_opt} then specifies a range of periods where aborting the mission is optimal when the system is in the most ``unhealthy'' state $\bm e_{m_1+m_2}$.
Theorem~\ref{thm:struct_policy} ahead shows that the optimal mission abort policies before and after time $\tilde n\delta$ exhibit different structures.
% \ye{[it would be better to add one to two sentences to discuss why this range is important/useful]}

The above analysis, together with the two thresholds $\hat n$ in Theorem~\ref{thm:threshold_n} and $\tilde n$ in Lemma~\ref{lemma:ab_opt}, sets the basis to derive the optimal abort policy.
The main challenge in the analysis is the multivariate belief state $\bm\pi_n$ in an $(m_1+m_2)$-dimensional probability simplex $\mathcal{S}$, under which it is difficult to order an arbitrary pair of vectors $\bm\pi^{(1)},\bm\pi^{(2)}\in\mathcal{S}$.
%Even in the sense of MLR ordering, there exist $\bm\pi^{(1)},\bm\pi^{(2)}\in\mathcal{S}$ such that neither $\bm\pi^{(1)}\succeq_\text{LR}\bm \pi^{(2)}$ nor $\bm\pi^{(2)}\succeq_\text{LR}\bm \pi^{(1)}$.
To circumvent this difficulty, we follow \cite{wang2015multistate} and transform $\bm\pi$ from a Cartesian coordinate system to a spherical coordinate system with the vertex $\bm e_{m_1+m_2}$ as the origin.
Specifically, a vector $\bm\pi\in\mathcal{S}\backslash\{\bm e_{m_1+m_2}\}$ can be represented by $(r,\bm\Phi)\triangleq (r(\bm\pi),\bm\Phi(\bm\pi))$ in this spherical coordinate system by the transformation in Table~\ref{tab:spherical}.
Here, $r$ is the spherical radius and $\bm\Phi\triangleq[\phi_1,\ldots,\phi_{m_1+m_2-1}]'$ is the angle vector.
The transformations in Table~\ref{tab:spherical} ensure $\phi_j\in[0,\pi/2]$ for all $j\in[m_1+m_2-2]$ and $\phi_{m_1+m_2-1}\in(\pi/2,\pi]$.
Under the spherical coordinate system, we can easily order two vectors with the same angle vector $\bm\Phi$ as shown below.
\begin{proposition}\label{prop:order}
	%\bfred{[check if you need any assumption/condition]}	
	Consider the original space $\mathcal{S}=\{\bm\pi=(\pi_i)_{i\in[m_1+m_2]}\in\mathbb{R}^{m_1+m_2}_+:\sum_{i=1}^{m_1+m_2}\pi_{i}=1\}$ and the corresponding spherical coordinate system using the transformation in Table~\ref{tab:spherical}.
	For any two vectors in $\mathcal{S}$ with spherical coordinates $(r_1,\bm\Phi)$ and $(r_2,\bm\Phi)$,
	their Cartesian coordinates $\bm\pi^{(1)}\triangleq\bm\pi^{(1)}(r_1,\bm\Phi)$ and $\bm\pi^{(2)}\triangleq\bm\pi^{(2)}(r_2,\bm\Phi)$ satisfy $\bm\pi^{(1)}\succeq_\textnormal{LR}\bm\pi^{(2)}$ if $r_1<r_2$.
\end{proposition}

For a fixed angle vector $\bm\Phi$, the radius $r$ can be understood as the distance of the belief state $(r,\bm\Phi)$ to the most ``unhealthy'' state $\bm e_{m_1+m_2}$.
Proposition~\ref{prop:order} then implies that the system failure risk decreases in $r$.
Based on this observation, we show that when the angle vector is fixed, the optimal abort policy is defined by a lower and an upper control limit with respect to the radius $r$.

\begin{theorem}\label{thm:struct_policy}
	%\bfred{[check if you need any assumption/condition]}	
	Suppose that Assumptions~\ref{assump:nondecreasing}--\ref{assump:rescue} hold and $\lambda>\bar\lambda$ {for the $\bar\lambda$ in Lemma~\ref{lemma:increase_rate}}.
	Let $(r_n,\bm\Phi_n)$ be the system belief state at time $n\delta$ under the spherical coordinate system in Table~\ref{tab:spherical}.
	When the condition in \eqref{eq:cost_ratio} does not hold, the optimal mission abort policy has the following structure. 
	\begin{itemize}
		\item[(i)] For $n=0,\ldots,\tilde n$, 
		there exists an upper control-limit function $\bar r_n(\bm\Phi_n):[0,\pi/2]^{m_1+m_2-2}\times(\pi/2,\pi] \mapsto \mathbb{R}_+$ such that it is optimal to abort the mission if and only if $r_n<\bar r_n(\bm\Phi_n)$. 
		\item[(ii)] For $n=\tilde n+1,\ldots,\hat n-1$, there exist lower and upper functions $\underline r_n,\bar r_n:[0,\pi/2]^{m_1+m_2-2}\times(\pi/2,\pi] \mapsto \mathbb{R}_+$ such that mission abort is the optimal action if and only if $r_n\in [\underline r_n(\bm\Phi_n),\bar r_n(\bm\Phi_n)]$, where the interval is allowed to be empty. 
		% \ye{[since your $\hat{n}$ is not the smallest $n$ after which abort is never optimal, within your $n$ range here, some $n$ might have never abort as optimal action, and some $n$ might have an optimal abort interval. I would suggest this way to state the optimal policy rather than forcing $\underline r_n(\bm\Phi_n)\leq r_n\leq\bar r_n(\bm\Phi_n)$, which means it must be non empty.]}
		\item[(iii)] The set of states $\{\bm\pi\in\mathcal{S}:V_\textnormal{ab}(n,\bm\pi)\leq V_\textnormal{c}(n,\bm\pi) \}$ is convex for all $n=0,\ldots,N-1$. 
	\end{itemize}
\end{theorem}

\begin{table}[]
	\centering
	\caption{The spherical transform and inverse transform of the Cartesian belief vector ($\bm \pi\neq \bm e_{m_1+m_2}$).}
	\label{tab:spherical}
	\small
	\begin{tabular}{cc}
		\hline
		Cartesian $\bm\pi$ to spherical $(r,\bm\Phi)\triangleq(r(\bm\pi),\bm\Phi(\bm\pi))$ 
		& Spherical $(r,\bm\Phi)$ to Cartesian $\bm\pi\triangleq \bm\pi(r,\bm\Phi)$ \\
		\hline
		$r=\sqrt{\sum_{i=1}^{m_1+m_2-1}\pi_i^2+(\pi_{m_1+m_2}-1)^2}$
		& $\pi_2=r\cos\phi_1$ \\
		$\phi_j=\arccos\frac{\pi_{j+1}}{\sqrt{\pi_{j+1}^2+\cdots+(\pi_{m_1+m_2}-1)^2+\pi_1^2}}$, 
		& 
		$\pi_j=r\cos \phi_{j-1}\prod_{k=1}^{j-2}\sin\phi_ k$, \\
		\qquad \qquad \qquad \qquad \qquad \qquad \qquad
		$j\in[m_1+m_2-2]$ 
		& \qquad \qquad \qquad  $j=3,\ldots,m_1+m_2-1$ \\
		$\phi_{m_1+m_2-1}=\arccos\frac{\pi_{m_1+m_2}-1}{\sqrt{\pi_1^2+(\pi_{m_1+m_2}-1)^2}}$ 
		& $\pi_{m_1+m_2}=1+r\cos \phi_{m_1+m_2-1}\prod_{k=1}^{m_1+m_2-2}\sin\phi_k$ \\
		& 
		$\pi_1=-r\left(\cos\phi_1+\sum_{j=2}^{m_1+m_2-1}\cos\phi_j\prod_{k=1}^{j-1}\sin\phi_k\right)$ \\
		\hline
	\end{tabular}
\end{table}%
The above theorem naturally divides the mission into three stages.
By Proposition~\ref{prop:order}, a smaller $r_n$ indicates that the system is more likely in a defective state when $\bm\Phi_n$ is fixed.
Intuitively, Theorem~\ref{thm:struct_policy}(i) indicates that the mission should be aborted when it is very likely to be defective at the initial stage of the mission.
Nevertheless, Theorem~\ref{thm:struct_policy}(ii) presents a somewhat counter-intuitive result: during the middle stage of the mission, mission termination is optimal only if the system is in certain intermediate states.
During such stages, it is optimal to continue the mission if the system is healthy.
Interestingly, the existence of the lower control limit $\underline r_n(\bm\Phi_n)$ implies that continuation is also optimal if the system is close to failure.
This lower bound results from the rescue time $w_n$, which is assumed nondecreasing in $n$.
If the system is very ``unhealthy'' upon abort at the middle stage, it is still highly likely to fail during the rescue procedure.
Consequently, mission abort cannot effectively mitigate the failure risk, and it is worth continuing the mission even though the failure probability is large.
%{A similar counter-intuitive optimal policy is found in \cite{ross1971quality} for quality control problems.}
Leveraging the piecewise-linear concavity of $V(n,\bm\pi)$ in Lemma~\ref{lemma:piecewise_and_concave}, the existence of the threshold $\tilde n$ in Lemma~\ref{lemma:ab_opt}, and the convexity of the set of belief states wherein an abort is optimal in Theorem~\ref{thm:struct_policy}(iii), 
%\ye{[which structure exactly? Please be specific here] the structure of the optimal abort policy}, 
Section~\ref{sec:algo} develops an efficient algorithm to solve the POMDP.

The following proposition provides two sufficient conditions under which $\{\tilde n+1,\ldots,\hat n-1\}$ is an empty set, rendering the middle stage in Theorem~\ref{thm:struct_policy}(ii) nonexistent.
In this case, we need only the upper bound $\bar r_n(\bm\Phi_n)$ to characterize the mission abort region for each $n=0,\ldots,\hat n-1$.

\begin{proposition}\label{prop:no_intermediate}
	%\bfred{[check if you need any assumption/condition]}	
	Suppose that Assumptions~\ref{assump:nondecreasing} and \ref{assump:rescue} hold and $\lambda>\bar\lambda$.
	We have $\tilde{n}\geq \hat n-1$ if either conditions below is met:
	\begin{itemize}
		\item [(i)] The optimal solution to the maximization problem:
		\begin{equation}\label{eq:remove_counter}
			\max_{i\in[m_1+m_2]} \left(1+\frac{C_\textnormal{m}}{C_\textnormal{s}}\right)p_{i,m_1+m_2+1}((N-\hat n+1)\delta+w_N)-p_{i,m_1+m_2+1}(w_{\hat n-1})
		\end{equation}
		is $i^*=m_1+m_2$.
		\item [(ii)] The rescue time satisfies $w_n=0$ for all $n=0,\ldots,N$.
	\end{itemize}
\end{proposition}

To understand condition (i) in Proposition~\ref{prop:no_intermediate}, recall that $p_{i,m_1+m_2+1}((N-\hat n+1)\delta+w_N)$ and $p_{i,m_1+m_2+1}(w_{\hat n-1})$ are the system failure probabilities of never terminating since time $(\hat n-1)\delta$ and 
%\ye{[replace as: "never terminating since"?] continuing} 
terminating at $(\hat n-1)\delta$, respectively.
The objective function in \eqref{eq:remove_counter} can then be understood as the gain of aborting the mission instead of continuing till the end, weighted by the costs $C_\textnormal{m}$ and $C_\textnormal{s}$.
If the gain attains its maximum when the system is in the most ``unhealthy'' state $\bm e_{m_1+m_2}$, we only need the upper control-limit function $\bar r_n$ in each period to make the abort decision.
Moreover, when the rescue time keeps zero during the mission, we also have $\tilde n\geq \hat n-1$.
This scenario could happen for paper mills, where we can halt the system in a negligible time to prevent imminent failures such as the breaking of bolts \citep{ranjan2018dataset}.
Interestingly, when $w_n$ remains zero throughout the mission, our mission abort problem resembles a CBM problem, which typically assumes negligible maintenance time as discussed in Section~\ref{sec:intro}.
	A potential violation of condition (ii) in Proposition~\ref{prop:no_intermediate} underscores a distinction between mission abort and CBM problems.
%Numerically, we find that the optimal abort policy only exhibits a one-threshold control limit structure under most parameter settings, even when the sufficient conditions in Proposition~\ref{prop:no_intermediate} are violated.
%This may partially explain why the two-threshold control limit policy in Theorem~\ref{thm:struct_policy}(ii) is not common in practice.

When Assumption~\ref{assump:nondecreasing} does not hold, Theorem~\ref{thm:threshold_n} may not hold, i.e., we cannot guarantee the existence of a time threshold $\hat n\delta$ after which the mission is never aborted.
Nevertheless, the structure of the optimal mission abort policy described in Theorem~\ref{thm:struct_policy} keeps largely the same.
We only need to replace $\hat n$ in Theorem~\ref{thm:struct_policy}(ii) with $N$ as detailed below.
	
\begin{corollary}\label{coro:opt_policy_wo_a1}
	%\ye{[if your assumption 1 does not hold, then what is your $\bar{\lambda}$?]}
	Under Assumptions~\ref{assump:TP2} and \ref{assump:rescue},
	the optimal mission abort policy for $n=0,\ldots,\tilde n$ has the same structure as that in Theorem~\ref{thm:struct_policy}(i), and the optimal policy for $n=\tilde n+1,\ldots,N-1$ has the same structure as that in Theorem~\ref{thm:struct_policy}(ii).
	Moreover, Theorem~\ref{thm:struct_policy}(iii) still holds.
\end{corollary}

\section{Computational Solution Approach}\label{sec:algo}
The belief state space is an $(m_1+m_2)$-dimensional probability simplex $\mathcal{S}$.
Due to the curse of dimensionality,
it is computationally intractable to solve the Bellman equations \eqref{eq:bellman_base} and \eqref{eq:bellman} to optimality when $m_1$ or $m_2$ is large \citep{bertsekas2012dynamic}.
Nevertheless, we can significantly reduce the computational burden of our problem by using the structural properties derived in Section~\ref{sec:property}.

First,
Theorem~\ref{thm:threshold_n} indicates the existence of a time threshold $\hat n$ under Assumption~\ref{assump:nondecreasing}, such that it is optimal to continue the mission for $n=\hat n,\ldots,N$ no matter of the observed signal.
To find $\hat n$, we need to solve the maximization problem: 
%\ye{
$\max_{\bm\pi\in\mathcal{S}}\{\overline V_\textnormal{c}(n,\bm\pi)-V_\textnormal{ab}(n,\bm\pi)\}$ 
%[do you need to add $\{\}$ for the objective?]} 
as shown in \eqref{eq:hat_n} for $n=0,\ldots,N$.
%\ye{[i shall say $\hat{n}$ is not unique. you are searching for an upper bound for it if you solve this program.]}
%\sun{[$\hat n$ is unique because we have min in (16)? The previous $\tilde n$ may not be unique, and I have corrected Lemma 6]}
Since $\mathcal{S}$ is a probability simplex, the above maximization problem is a linear program, and the optimal solution is one of the vertices of $\mathcal{S}$ \citep{bertsekas2012dynamic}, i.e., one of $\bm e_1,\ldots,\bm e_{m_1+m_2}$. 
For fixed $n$, we only need to enumerate $\bm\pi\in\{\bm e_1,\ldots,\bm e_{m_1+m_2}\}$ to check if $\max_{\bm\pi\in\mathcal{S}}\{\overline V_\textnormal{c}(n,\bm\pi)-V_\textnormal{ab}(n,\bm\pi)\}\leq 0$, based on which we can know whether $\hat n \leq n$ or $\hat n > n$.
Then a binary search can be adopted to find $\hat n$, as detailed in Algorithm~\ref{algo:binary} in Appendix~\ref{appen:implement}.
The binary search takes $\mathcal{O}(\log(N))$ steps to find the time threshold $\hat n$, where each step enumerates $\bm\pi\in\{\bm e_1,\ldots,\bm e_{m_1+m_2}\}$ to solve the linear program $\max_{\bm\pi\in\mathcal{S}}\{\overline V_\textnormal{c}(n,\bm\pi)-V_\textnormal{ab}(n,\bm\pi)\}$.

For $n<\hat{n}$, exactly computing the value function $V(n,\bm\pi)$ is computationally intractable.
We approximate $V(n,\bm\pi)$ by using the structural properties derived in Section~\ref{sec:property}.
Specifically, Lemma~\ref{lemma:piecewise_and_concave} shows that $\bm\pi\mapsto V(n,\bm\pi)$ is piecewise linear and concave for all $n=0,\ldots,N$.
This means that for each $n$, there exists a finite set $\mathcal{A}_n\subseteq\mathbb{R}^{m_1+m_2}$ such that $V(n,\bm\pi)$ can be written as
\begin{equation}\label{eq:min-linear}
	V(n,\bm\pi)=\min_{\bm\alpha_n\in\mathcal{A}_n} \bm\alpha_n'\bm\pi, 
	\quad n=0,\ldots,N.
\end{equation} 
%The size of $\mathcal{A}_n$ can be very large, 
Solving the Bellman equations in \eqref{eq:bellman_base} and \eqref{eq:bellman} amounts to obtaining the finite sets $\mathcal{A}_n$, $n=0,\ldots,N$.
This motivates us to approximate $\mathcal{A}_n$ by a sequence of sets $\{\widehat{\mathcal{A}}_n^{(\tau)},~\tau\in\mathbb{N}_+\}$ for $n=0,\ldots,N$, where each $\widehat{\mathcal{A}}_n^{(\tau)}$ is a subset of $\mathcal{A}_n$.
With $\widehat{\mathcal{A}}_n^{(\tau)}$, we can approximate the value function $V(n,\bm\pi)$ as 
\begin{equation}
	\widehat V^{(\tau)}(n,\bm\pi)=\min_{\bm\alpha_n\in\widehat{\mathcal{A}}_n^{(\tau)}}\bm\alpha_n' \bm\pi,
	\quad n=0,\ldots,N, \nonumber
\end{equation}
where $V(n,\bm\pi)\leq\widehat{V}^{(\tau)}(n,\bm\pi)$ for all $n=0,\ldots,N$, $\bm\pi\in\mathcal{S}$, and $\tau\in\mathbb{N}_+$ because $\widehat{\mathcal{A}}_n^{(\tau)}\subseteq\mathcal{A}_n$.
To obtain $\{\widehat{\mathcal{A}}_n^{(\tau)},~\tau\in\mathbb{N}_+\}$, we propose a modified PBVI algorithm to iteratively construct $\widehat{\mathcal{A}}_n^{(\tau)}$.

The classical PBVI algorithm \citep{pineau2003point} iterates between a backup step and an expansion step.
The backup step uses a finite set of belief states $\{\bm\pi^{(1)}_n,\ldots,\bm\pi^{(L_\tau)}_n\}$ to construct $\widehat{\mathcal{A}}_n^{(\tau)}$ for $n\in[\hat n-1]$, 
where $L_\tau$ is the size of such a set at the $\tau$th PBVI iteration. 
%\ye{[does $L_\tau$ depend on $n$? or does it depend on the iteration number only?]} 
%\sun{[it can depend on $n$ but we keep it the same across different $n$ and denote it to be $L_\tau$ for simplicity]}
The construction of $\widehat{\mathcal{A}}_{n}^{(\tau)}$ is done by approximating the cost-to-go for continuing the mission $V_\text{c}(n,\bm\pi_n)$ in \eqref{eq:bellman_base} as $\widehat V^{(\tau)}_\textnormal{c}(n,\bm\pi_n)$, where $V^{(\tau)}_\textnormal{c}(\cdot,\cdot)$ is defined analogously to \eqref{eq:v_c} but with $\widehat{V}^{(\tau)}(\cdot,\cdot)$ replacing $V(\cdot,\cdot)$ therein.
This approximation enables us to approximately solve the Bellman equation to obtain $\widehat{A}_{n}^{(\tau)}$ for each $n\in[\hat n-1]$ in a backward manner (c.f.\ Line 8 of Algorithm~\ref{algorithm:PBVI}).
Then the expansion step augments $\{\bm\pi_n^{(l)}\}_{l\in [L_\tau]}$ by sampling new elements $\bm\pi_n^{(l)}$, $l=L_\tau+1,\ldots,L_{\tau+1}$, for each $n\in[\hat n-1]$ for constructing $\widehat{\mathcal{A}}_n^{(\tau+1)}$ in the next iteration.
We initialize the algorithm with a finite set of belief states $\{\bm\pi^{(1)}_n,\ldots,\bm\pi^{(L_1)}_n\}$ for each $n\in[\hat n-1]$ by simulation, as detailed in Algorithm~\ref{algorithm:PBVI} ahead.
The algorithm terminates if the approximate value functions between two consecutive iterations are close enough or after a predetermined number of iterations.
%\ye{[I am confused here. your  $\widehat{A}_n^{(0)}$ should consist of the coefficient vector $\bm\alpha_n$. But now you initiate the algorithm using a set of belief states rather than a set of $\bm\alpha_n$? Belief state and $\bm\alpha_n$ can have the same dimension, but $\bm\alpha_n$ does not need to be in the probability simplex. It would be good to add one more sentence here to link (or explain how to link) the initialized belief states to $\bm\alpha_n$ vectors.]}
%\ye{[i cannot understand the idea of this algorithm based on your description. 
%You have already obtained $\widehat{\mathcal{A}}_n^{(\tau+1)}$ in the backup step? Then why do you need to sample again in the expansion step?
%Something should be wrong here. please relook and rewrite.]} 

Generally, the classical PBVI algorithm does not use any structural properties of $V(n,\bm\pi)$ in the iteration.
We develop a modified PBVI algorithm by using the structural properties of the optimal mission abort policy to further reduce the computational cost.
% \ye{[i would say this is novel idea and not a straightfoward extension of the original PBVI algorithm. Original PBVI, according to your description above, aims to find $\bm\alpha_n$, the coefficients to approximate the value function. But we are looking for the optimal policy. Plesae think deeper into this. To make this thinking coherent, you need to explain why this candidate belief set can be linked to the coefficient $\bm\alpha_n$.]}
Specifically, Theorem~\ref{thm:struct_policy} indicates that 
%\ye{[in the optimal policy, the mission abort region?] the set} 
the region $\{\bm\pi\in\mathcal{S}:V_{\textnormal{ab}}(n,\bm\pi)\leq V_\textnormal{c}(n,\bm\pi)\}$ with mission abort as the optimal action is convex for all $n=0,\ldots,N-1$.
For a fixed $n$ and a large $\tau$, the above set can be well approximated by the convex hull
% \ye{[it seems to me that $\widehat V^{(\tau)}_\textnormal{c}$ is undefined so far]}
\begin{equation}\label{eq:ab_opt_hull}
	\mathrm{conv}\left\{\bm\pi_n^{(l)},~l\in[L_\tau]:V_{\textnormal{ab}}(n,\bm\pi_n^{(l)})\leq \widehat V^{(\tau)}_\textnormal{c}(n,\bm\pi_n^{(l)})\right\}.
\end{equation}
In words, mission abort is very likely to be the optimal action for the belief states in the convex hull in \eqref{eq:ab_opt_hull}.
If aborting the mission is optimal for any belief states $\bm\pi$ in period $n$, then $V(n,\bm\pi)=V_\text{ab}(n,\bm\pi)$ is linear in $\bm\pi$ as shown in \eqref{eq:v_ab}, which satisfies $V(n,\bm\pi)=\bm\alpha_n'\bm\pi$ for some $\bm\alpha_n\in\widehat{\mathcal{A}}_n^{(\tau)}$.
As such, all such belief states share the same coefficient $\bm\alpha_n$.
To gain more information in later iterations, %\ye{[this argument depends heavily on the link between $\bm\alpha_n$ and the belief states. you might need deeper reasoning to support this argument.] the expansion step should expand the set $\{\bm\pi_n^{(l)}\}_{l\in[L_\tau]}$ by adding elements falling outside of the above convex hull to explore the belief states with less known information.}
%To this end, 
a sampled state inside \eqref{eq:ab_opt_hull} is discarded during the expansion step when $\tau$ is large.
%\ye{[why do you make this assumption? Please justify here]}
Moreover, in view of Lemma~\ref{lemma:ab_opt}, if $V_\text{ab}(\tilde n,\bm e_{m_1+m_2})\leq \widehat{V}_\text{c}^{(\tau)}(\tilde n,\bm e_{m_1+m_2})$ for a large $\tau$, it is likely that $V_\text{ab}(n,\bm e_{m_1+m_2})\leq \widehat{V}_\text{c}^{(\tau)}(n,\bm e_{m_1+m_2})$ for all $n\leq \tilde n$.
When we find such $\tilde n$ during the expansion step, we can include the state $\bm e_{m_1+m_2}$ to \eqref{eq:ab_opt_hull} and update the convex hull for all $n\leq \tilde n$. 
Other steps remain the same as the classical PBVI algorithm.
%This modified PBVI algorithm is motivated by \cite{liu2022machine} that works on a one-dimensional belief space.
%Our problem is more complex since the state space $\mathcal{S}$ is an $(m_1+m_2)$-dimensional probability simplex, entailing the above proposed procedure to accelerate the computation.
This modified PBVI algorithm is summarized in Algorithm~\ref{algorithm:PBVI}, and details on the backup and expansion steps are elaborated in Appendix~\ref{appen:implement}.
With the approximate value function obtained by the PBVI algorithm, the optimal action at each belief state can be readily obtained.
Numerical experiments in Section~\ref{sec:num} show that the modified PBVI algorithm achieves a significant reduction in computational time compared with the classical PBVI algorithm.

\begin{algorithm}[t]
	\caption{The modified PBVI algorithm to approximately calculate the value function.}
	\begin{algorithmic}[1]
		{\renewcommand\baselinestretch{1}\selectfont
		\Require
		Transition rate matrix $\mathbf{Q}$, state-observation matrix $\mathbf{D}$, costs $C_{\text{m}},C_{\text{s}}$, inspection interval $\delta$, number of periods $N$, the set of $\bm\alpha$-vectors $(\mathcal{A}_n)_{n=\hat n}^N$, numbers of iterations $Z_1,Z_2$, initial number of states $L_1$, error $\varepsilon$, number of samples to simulate $W$
		\Ensure
		The approximate value function
		\State Init $\widehat{\mathcal{A}}^{(\tau)}_n\leftarrow \mathcal{A}_n$ for all $n=\hat n,\ldots,N$ and $\tau\geq 0$, and init $\widehat{\mathcal{A}}_n^{(0)}\leftarrow\emptyset$ for all $n\in[\hat n-1]$
		\State Running $L_1$ sets of simulations from state $\bm\pi_0$ to the end of the decision horizon
		\State Compute the sets of belief states $\{\bm \pi_n^{(l)}\}_{l\in L_1}$, $n\in[\hat n-1]$, by Equation \eqref{eq:bayes_update} based on the $L_1$ trajectories of simulated signals
		\For{$\tau\leftarrow 1,\ldots,Z_1$}
		\For{$n\leftarrow \hat n-1,\ldots,1$}
		% \State Set $\widehat{\mathcal{A}}^{(\tau)}_n\leftarrow \widehat{\mathcal{A}}^{(\tau-1)}_n$
		\For{$l\leftarrow 1,\ldots,L_\tau$}
		\State Compute $V_\textnormal{ab}(n,\bm\pi_n^{(l)})$ by \eqref{eq:v_ab} and $\widehat{V}^{(\tau)}_\textnormal{c}( n,\bm\pi_n^{(l)})$ as 
		$$
		\kappa(\bm\pi_n^{(l)},\delta)(C_\textnormal{s}+C_\textnormal{m})+\sum_{k=1}^K\min_{\bm\alpha_{n+1}\in\widehat{\mathcal{A}}^{(\tau)}_{n+1}}\sum_{i=1}^{m_1+m_2}\sum_{j=1}^{m_1+m_2}\alpha_{n+1,j}\pi_{ni}^{(l)} p_{ij}(\delta) \tilde d_{jk}
		$$
		\State \parbox[t]{\dimexpr 14.1cm}{%
			Set $\widehat{\mathcal{A}}_n^{(\tau)}$ be the collection of coefficients corresponding to the belief vector
			$\bm{\pi}_n^{(l)}$ in the function $\min\{V_{\text{ab}}(n,\bm{\pi}_n^{(l)}),\widehat V^{(\tau)}_{\text{c}}(n,\bm{\pi}_n^{(l)})\}$
		}
		\State \parbox[t]{\dimexpr 14.1cm}{%
			Sample $W$ underlying states from $\bm\pi_n^{(l)}$ and simulate a signal for each sampled state to compute the belief state at time $(n+1)\delta$
		}
		\If{$\tau>Z_2$}
		\State Discard all belief states at time $(n+1)\delta$ in the convex hull in \eqref{eq:ab_opt_hull}
		\EndIf
		\State \parbox[t]{\dimexpr 14.1cm}{%
			Compute the minimum Euclidean distance between each of the above (kept) belief states and each state in $\{\bm\pi_{n+1}^{(l)}\}_{l\in[L_\tau]}$
		}
		\State Set $\bm\pi^{(L_\tau+l)}_{n+1}$ to be the belief state that has the largest distance
		\EndFor
		\EndFor
		\If{$\max_{\bm\pi\in\mathcal{S},n=0,\ldots,N}|\widehat{V}_n^{(\tau)}(\bm\pi)-\widehat{V}_n^{(\tau-1)}(\bm\pi)|<\varepsilon$}
		\State Break the for-loop
		\EndIf
		\EndFor
		\par}
	\end{algorithmic}
	\label{algorithm:PBVI}
\end{algorithm}

Recall that the value function $V(n,\bm\pi)$ is MLR nondecreasing in $\bm\pi$ by Lemma~\ref{lemma:MLR_montonicity} when Assumption~\ref{assump:nondecreasing} holds.
According to \cite{tian2022adaptive}, approximate dynamic programming (ADP) algorithms using a finite set of sampled states, such as the PBVI algorithms, typically yield accurate approximations of $V(n,\bm\pi)$ for $\bm\pi$ not visited by the algorithm ($\bm\pi\notin \{\bm\pi_n^{(i)}:i=1,2,\ldots\}$), provided that the approximation maintains monotonicity of the value function. 
%\ye{[i cannot understand this. please make your statement clearer] 
% in the states that the algorithm does not visit.}
The following proposition illustrates that our PBVI algorithm can preserve such structural properties of $V(n,\bm\pi)$.
\begin{proposition}\label{lemma:preservation}
%\bfred{[check if you need any assumption/condition]}	
	Suppose that Assumptions~\ref{assump:nondecreasing}--\ref{assump:rescue} hold and $\lambda>\bar\lambda$ {for the $\bar\lambda$ in Lemma~\ref{lemma:increase_rate}}.
	The approximate value function $\widehat V^{(\tau)}(n,\bm\pi)$ is piecewise linear, concave, and MLR nondecreasing in $\bm\pi\in\mathcal{S}$ for all $n=0,\ldots,N$ and $\tau\in\mathbb{N}_+$.
\end{proposition}

In view of the preservation of the structural properties of $V(n,\bm\pi)$, our PBVI algorithm is more appealing than other ADP methods, 
%\ye{[by saying `other ADP methods'' here, do you mean PBVI is also an ADP? if yes, please make this statement clearly right after you introduce the name `AD'.]}
e.g., using a neural network to approximate $V(n,\bm\pi)$, that generally do not have the preservation property in Proposition~\ref{lemma:preservation}.
Moreover, the PBVI algorithm is an offline method as it optimizes a policy before realization of observed signals.
After optimization, the optimal policy can be implemented online to yield the optimal action in real time.
In contrast, online planning methods such as Monte Carlo tree search \citep{silver2010monte} may not be suitable for real-time abort decision-making when the inter-sampling interval $\delta$ is very small.

\section{Insights into Specialized Models}\label{sec:special}
This section derives some additional structural results for two special cases of our model.
The first case in Section~\ref{subsec:special_CTMC} investigates the decision problem when all $T_{12}$, $T_{13}$, and $T_{23}$ are exponential.
The second case in Section~\ref{subsec:special_phase} assumes an exponential distribution for $T_{23}$ while $T_{12}$ remains non-exponential.
% Our results provide insights on how to reduce the computational cost when some transition times are exponentially distributed.
In both cases, the POMDP can be exactly solved after discretizing the state space.

\subsection{Markovian Failure Process}\label{subsec:special_CTMC}
We first study the special case where $T_{12}$, $T_{13}$, and $T_{23}$ are all exponential, under which the underlying system dynamics is a CTMC.
Its optimal mission abort decision making is a POMDP, and our POMDP analysis above applies to this problem directly with $m_1=m_2=1$.
This simple setting lends itself to a straightforward structure for the optimal policy.
% \ye{[defective state is 1 or 2?] 
Let $\Pi_n=\mathbb{P}(X_n=2\mid Y_{1:n})$ be the conditional probability that the system is at the defective state given $Y_{1:n}$.
After observing the values of $Y_{1:n}$, the system state at time $n\delta$ is realized as a scalar $\pi_n\in[0,1]$.
The probability of the system being at the healthy state is then $1-\pi_n$.

The general result in Theorem~\ref{thm:threshold_n} applies to this simple setting.
That is, there exists $\hat n$ such that, for all $n\geq\hat n$, continuation is optimal regardless of the system state $\Pi_n$.
Similarly, the threshold $\tilde n$ in Lemma~\ref{lemma:ab_opt} becomes $\tilde n=\max\{n=0,\ldots,N:V_\textnormal{ab}(n,1)< V_\textnormal{c}(n,1)\}$.
The univariate belief state $\Pi_n$ obviates the need for introducing the spherical coordinate system to depict the optimal policy when $n<\hat n$.
The optimal policy can be shown to be a control-limit policy with respect to $\pi_n$.
\begin{corollary}\label{thm:control_limit_CTMC}
%\bfred{[check if you need any assumption/condition]}	
	Suppose $X(t)$ is a CTMC.
	Under Assumptions~\ref{assump:nondecreasing}--\ref{assump:rescue}, the optimal mission abort policy has the following structure.
	\begin{itemize}
		\item For $n=0,\ldots,\tilde n$, there exists a threshold $\underline{\pi}_n$ such that it is optimal to abort the mission at time $n\delta$ if and only if the system belief state $\pi_n$ satisfies $\pi_n>\underline{\pi}_n$ and to continue the mission otherwise.
		% \ye{[it would be better to use $\underline{\pi}_n$ to replace your $\underline{\pi}_n$ here. If adopt, please also revise the supp]}
		\item For $n=\tilde n+1,\ldots,\hat n -1$, there exists an interval $[\underline\pi_{n},\overline\pi_{n}]\subseteq[0,1]$ such that it is optimal to abort the mission if and only if $\underline\pi_{n}\leq\pi_n\leq\overline\pi_{n}$, where the interval is allowed to be empty. 
		% \ye{this is very strange. In case (i), the middle stage does not exist. So you don't have the optimal policy. Therefore, i would suggest that you remove case (i).]}
		\item 
		For $n\geq \hat n$, the optimal policy is to continue the mission.
	\end{itemize}
\end{corollary}

{
Compared with optimal CBM for a system with three states (e.g., \citealp{kim2010contracting}), 
% \ye{[ref here]}
the fundamental difference comes from the middle stage where continuation of the mission is optimal upon confirmation of the defective state and the last stage where abort is never optimal.
From Corollary~\ref{thm:control_limit_CTMC}, it is interesting to see that a POMDP based on a three-state CTMC makes mission abort decisions only using the defective probability $\pi_n$.
For a non-Markovian system, the sojourn time on the defective state is informative for decision-making.
Theorem~\ref{thm:struct_policy} implicitly considers this information using a spherical coordinate system.
As shown by our numerical experiments in Section~\ref{sec:num}, loss of such critical information can lead to a significant cost increase if we approximate a semi-Markov system by a simple three-state CTMC.}

Since the belief state space can be represented by a one-dimensional interval $[0,1]$, it is unnecessary to use the PBVI algorithm in Section~\ref{sec:algo} to approximately solve the POMDP.
Instead, we can discretize the interval $[0,1]$ and use standard backward induction to exactly solve the stochastic dynamic program (see Algorithm~\ref{algo:backward} in Appendix~\ref{appen:implement}).

\subsection{One-Phase Transition}\label{subsec:special_phase}
\iffalse
As illustrated in Section~\ref{subsec:model_deg}, the transition rate from healthy to defective state may hardly change during a mission due to the short duration of a mission compared with the system lifetime.
In this case, $T_{12}$ can be modeled by an exponential distribution, which means $m_1=1$.
The time $T_{23}$ from defect to failure still follows a general distribution $F(\cdot)$, and is approximated as an Erlang mixture random variable $T_{23}^{(\lambda)}$ with $m_2>1$ phases as before.
On the other hand, for completeness, we also consider the case that $T_{23}$ follows an exponential distribution ($m_2=1$) but we keep $m_1>1$.
Although less common in practice, this case is easier to analyze and facilitates the exposition of our dimension reduction idea introduced ahead.
\fi

Next, consider the surrogate POMDP with $m_1>1$ and $m_2=1$.
% Such a system may be susceptible to early failures, often stemming from inefficient burn-in procedures.
Motivated by \cite{wang2015multistate}, we can derive a dimension reduction property in this case that facilitates solving the POMDP.
Specifically, let 
$$
\varphi_j(\bm\Phi)\triangleq(p_{2j}(\delta)-p_{1j}(\delta))\cos \phi_1+\sum_{i=3}^{m_1+1}(p_{ij}(\delta)-p_{1j}(\delta))\cos\phi_{i-1}\prod_{k=1}^{i-2}\sin\phi_k
$$
for $j\in[m_1+1]$.
We arrive at the following proposition.
\begin{proposition}\label{prop:dim_red}
%\bfred{[check if you need any assumption/condition]}	
	Consider $m_1>1$ and $m_2=1$.
	Suppose that given $Y_1,\ldots,Y_{n-1}$, we have computed the system belief state $(r_{n-1},\bm\Phi_{n-1})$ under the spherical coordinate system defined in Table~\ref{tab:spherical}.
Then after observing $Y_n$ at time $n\delta$, $\bm\Phi_n=(\phi_{ni})_{i\in[m_1]}$ can be updated as
	\begin{equation}\label{eq:angle}
		\begin{aligned}
			\phi_{ni}&=\arccos
			\frac{\varphi_{i+1}(\bm\Phi_{n-1})}
			{\sqrt{\varphi_1^2(\bm\Phi_{n-1})+\sum_{j=i+1}^{m_1} \varphi_j^2(\bm\Phi_{n-1})+(\sum_{j=1}^{m_1}\varphi_j(\bm\Phi_{n-1}))^2}}, \quad i\in[m_1-1]; \\
			\phi_{n,m_1} & = \arccos \frac{-\sum_{i=1}^{m_1}\varphi_i(\bm\Phi_{n-1})}{\sqrt{\varphi_1^2(\bm\Phi_{n-1})+(\sum_{j=1}^{m_1}\varphi_i(\bm\Phi_{n-1}))^2}},
		\end{aligned}
	\end{equation}
while $r_n$ can be updated as a function of $r_{n-1}$ and the observed $Y_n$.
This function is complicated and is provided in Appendix~\ref{appen:proof_dim_red}.
\end{proposition}

Proposition~\ref{prop:dim_red} reveals that the angle vector $\bm\Phi_n$ at each decision epoch is \textit{irrelevant} to the signals $Y_1,\ldots,Y_n$ and can be deterministically computed by~\eqref{eq:angle} at time $0$.
In contrast, the radius $r_n$ depends on the realizations of $Y_1,\ldots,Y_n$.
This means that the signals affect the cost-to-go only through the scalar $r_n$.
Instead of using PBVI, we can discretize the domain of $r_n$ and then use standard backward induction to exactly solve the Bellman equation. 
Routine algebra shows that to meet the constraint $\bm\pi(r_n,\bm\Phi_n)\in\mathcal{S}$, $r_n$ should be within $[0, -(\cos \phi_{n,m_1}\prod_{k=1}^{m_1-1}\sin \phi_{nk})^{-1}]$.
%\ye{[help me check, is there a negative sign on the upper end point?]}.
With this range for $r_n$, the detailed optimization procedure is similar to that in Section~\ref{subsec:special_CTMC} and is relegated to Appendix~\ref{appen:implement}.

\iffalse
\ye{[The claim here is wrong, right? If wrong, then probably we can delete completely]} We next consider the case that $m_1=1$ and $m_2>1$.
Unfortunately, we cannot obtain a dimension reduction result as in Proposition~\ref{prop:dim_red} if we use Table~\ref{tab:spherical} to transform $\bm\pi$ to a vector in the spherical coordinate measured with respect to $\bm e_{m_1+m_2}$.
Nevertheless, if we construct a new spherical coordinate measured with respect to vertex $\bm e_1$ (rather than with respect to $\bm e_{m_1+m_2}$ as in Table~\ref{tab:spherical}), the computation of the angle vector $\bm\Phi_n$ (under the new spherical coordinate system) at each decision epoch can also be decoupled from the computation of $r_n$.
The derivation of this result largely follows the same procedure as that in the case of $m_2=1$, so we provide the details and some additional structural results under the new spherical coordinate system in Appendix~\ref{appen:m1=1}.
Since the dimension reduction property when $m_1=1$ is derived under a different spherical coordinate system, the structural results in Theorem~\ref{thm:struct_policy} no longer hold.
Nevertheless, the corresponding POMDP is a one-dimensional problem on the radius $r_n$, so it can also be solved quickly and exactly by standard backward induction after discretizing the state space (see Appendix~\ref{appen:implement}).
\fi

\section{Extension to Multiple Tasks}\label{sec:multiple}
%{\color{red}
%[why you consider multiple tasks throughout the lifetime? This make things much more complicated.]
%}
%\sun{[the ``lifetime planning'' was raised by R1. So I wrote in this way to tell R1 that we only need to consider the decision-making during consecutive thorough inspections]}
{In applications, a thorough system inspection may only be scheduled once the system has completed several scheduled tasks. 
Examples include a delivery drone tasked with multiple deliveries in a single day, and a UAV performing day-long line inspections in a power grid where the UAV may need to return for charging between flights.
Another scenario is that a big task can be divided into a few phases and each phase can be treated as a mission.
For example, a mission of tracking, telemetry, and command systems typically consists of two phases: data transformation and telecommand \citep{yu2020extended}.
% A graphical illustration of a possible multi-task setting is provided in Figure~\ref{fig:multi_insp}.
As such, there is no assurance that the system is in a perfectly healthy state at the start of a mission.
%\ye{[any example? After that, relate to Fig 2.]}

Consider a system that completes missions $1,\ldots,L$ in sequence, and each mission needs $N_l$ periods to complete.
A system failure incurs a cost $C_{\text{s}}$.
To prevent failures, we can abort the mission any time during the operation, 
upon which the current and subsequent missions fail, and failure to complete mission $l\in[L]$ incurs a loss $C_{\text{m}}^{(l)}$.
Upon completion of all $L$ missions or rescue procedure, a comprehensive inspection is carried out to ensure that the system is at the healthy state.
If the inspection reveals system defects, a repair cost $C_\text{r}$ is incurred.
After inspection, the system starts the next set of $L$ missions with a healthy state.
In this context, it would be advantageous to maximize the abort decision for all $L$ missions jointly.
Then the total number of decision periods is $N\triangleq \sum_{l=1}^L N_l$.

In practice, $L$ is usually a small number, and thus the assumptions on the sojourn times of the underlying semi-Markov system dynamics are still valid. 
% This model can be further extended to encompass long-term system deterioration. Every time the system restarts to finish the next $L$ missions, we can adjust the distributions $G$ and $F$ and re-construct our CTMC for decision-making. For illustration, we only use a pair of $(G,F)$ in our numerical studies in Section~\ref{subsec:num_multi}.
As such, we adhere to the model setting in Section~\ref{sec:model} for the system failure process, observed signals, and rescue procedures. 
{However, Assumption~\ref{assump:rescue} may not hold.
For example, consider a UAV that consecutively executes $L=3$ tasks by departing from a depot, visiting three sites in sequence, and returning to the depot.
If site 3 is closer to the depot than site 2, the rescue time is not necessarily monotone nondecreasing in $n\in[N]$.}
%There are only two slight differences.
%First, due to the repair cost $C_\text{r}$, we now consider a salvage value depending on the systems state upon completion of the mission or rescue procedure.
%Second, as shown in Figure~*, the rescue time $w_n$ may not be nondecreasing with $n$.
Nevertheless, the solution methodology is largely the same as before.
We first use two Erlang mixtures to approximate the non-exponential sojourn times $T_{12}$ and $T_{23}$ in the original semi-Markov process and formulate a CTMC.
The surrogate POMDP is then constructed by closely following that outlined in Section~\ref{sec:pomdp} with two minor changes.
First, the terminal cost needs to consider the possible repair cost after thorough inspection.
Second, the mission failure cost now depends on the number of unfinished missions.
Accordingly, the Bellman recursion in Eqs.~\eqref{eq:v_ab}--\eqref{eq:bellman} is changed to 
\begin{align}
	V_{\text{ab}}(n,\bm\pi_n)&=\sum_{l=1}^LC_{\text{m}}^{(l)}\mathbbm{1}\left\{n<\sum_{l'=1}^l N_{l'}\right\}+C_{\text{s}}\kappa(w_n,\bm\pi_n)+C_\text{r}\bm\pi_n' \tilde{\mathbf{P}}(w_n)\cdot(\bm 0_{m_1}',\bm 1_{m_2}')'; \label{eq:multi_ab} \\
	V_{\text{c}}(n,\bm\pi_n)&=\left(C_{\text{s}}+\sum_{l=1}^LC_{\text{m}}^{(l)}\mathbbm{1}\left\{n<\sum_{l'=1}^l N_{l'}\right\}\right)\kappa(\delta,\bm\pi_n)+\mathbb{E}^{(\lambda)}\left[V^{(\lambda)}(n+1,\bm\Pi_{n+1})\mid \bm\Pi_n=\bm\pi_n\right] \label{eq:multi_c}
\end{align}
for $n=0,\ldots,N-1$, where $\bm 0_r$ is the $r$-dimensional column vector of zeros.
The terminal cost is
\begin{equation}\label{eq:multi_terminal}
	V^{(\lambda)}(N,\bm\pi_N)\triangleq V_\text{c}(N,\bm\pi_N)=(C_{\text{s}}+C_{\mathrm{m}}^{(L)})\kappa(w_N,\bm\pi_N)+C_\text{r}\bm\pi_N' \tilde{\mathbf{P}}(w_N)\cdot(\bm 0_{m_1}',\bm 1_{m_2}')'.
\end{equation}

The last terms on the right hand side of \eqref{eq:multi_ab} and \eqref{eq:multi_terminal} represent the expected repair cost when system defects are revealed during the thorough inspection.
Because this cost is linear in the belief state, it is readily shown that $\bm\pi\mapsto V^{(\lambda)}(n,\bm\pi)$ remains concave for all $n\in[N]$.
% \ye{[which repair cost? Introduce it before pointing to it. For examle, you can say that the first term on the right hand side of $V_{xxx}$ represents the expected failure cost. Since it is a linear function of xxx]} 
Then the optimal mission abort policy in the multi-task setting shares a similar structure to those in Theorem~\ref{thm:struct_policy}.
Using Table~\ref{tab:spherical} to transform a belief vector $\bm\pi_n$ at period $n$ to $(r_n,\bm\Phi_n)$ in the spherical coordinate system, the optimal mission abort policy can be characterized as below.

%\ye{[it would be better to call it a theorem because it is a new setting, not a special case like the CTMC above]}

\begin{theorem}\label{thm:multi_opt}
	The optimal mission abort policy in the multi-task setting has the following structure.
	\begin{itemize}
		\item[(i)] The optimal action is always to continue the mission if
		\begin{equation}\label{eq:no_abort_multi}
			C_\textup{m}^{(L)}\frac{\exp(-\lambda(T+w_N))}{1-\exp(-\lambda(T+w_N))}\geq C_\textup{s}+\sum_{l=1}^{L-1}C_\textup{m}^{(l)}
			\frac{1-\exp(-\lambda\delta\sum_{l'=1}^l N_{l'})}
			{1-\exp(-\lambda(T+w_N))}.
		\end{equation}
		Otherwise, the set of states $\{\bm\pi\in\mathcal{S}:V_\textnormal{ab}(n,\bm\pi)\leq V_\textnormal{c}(n,\bm\pi) \}$ is convex for all $n=0,\ldots,N-1$.
		Moreover, aborting the missions is never optimal after completion of mission $l\in[L-1]$ if \eqref{eq:no_abort_after_l} in Appendix~\ref{appen:proof_thm:multi_opt} is satisfied.
		\item[(ii)] If $V_\textup{ab}(n,\bm e_{m_1+m_2})\leq V_\textup{c}(n,\bm e_{m_1+m_2})$, then there exists a function $\bar r_n(\bm\Phi_n):[0,\pi/2]^{m_1+m_2-2}\times(\pi/2,\pi] \mapsto \mathbb{R}_+$ such that it is optimal to abort the mission if and only if $r_n<\bar r_n(\bm\Phi_n)$. 
		% \ye{[because this is a quite different setting, it would be good to restate the optimal policy again.]}
		\item[(iii)] If $V_\textup{ab}(n,\bm e_{m_1+m_2})> V_\textup{c}(n,\bm e_{m_1+m_2})$, then there exist two functions $\underline r_n,\bar r_n:[0,\pi/2]^{m_1+m_2-2}\times(\pi/2,\pi] \mapsto \mathbb{R}_+$ such that the mission abort is optimal if and only if $r_n\in [\underline r_n(\bm\Phi_n), \bar r_n(\bm\Phi_n)]$, where the interval is allowed to be empty.
		% \ye{[because this is a quite different setting, it would be good to restate the optimal policy again.]}
	\end{itemize}
\end{theorem}
% \ye{[this sentence comes from nowhere. Where do you need to compare? Why you need to remark? Why is it useful/important? Please ALWAYS set stages before you perform.]}
Characterizing the optimal mission abort policy needs to compare $V_\text{ab}(n,\bm e_{m_1+m_2})$ and $V_\text{c}(n,\bm e_{m_1+m_2})$ as stated in Theorem~\ref{thm:multi_opt}(ii) and (iii).
Computing these functions for comparison is easy and does not require solving the Bellman equations in \eqref{eq:multi_ab}--\eqref{eq:multi_terminal}.
This is because the system can only transition to a worse state.
As such, we can obtain values of $V_\text{ab}(n,\bm e_{m_1+m_2})$ and $V_\text{c}(n,\bm e_{m_1+m_2})$ in a backward induction way by only using the values of $\{V(i,\bm e_{m_1+m_2})\}_{i=n+1}^N$.}

\section{Case Study}\label{sec:num}
We consider a fixed-wing UAV that executes an inspection task on high-voltage power grids located at a mountainous area.
This UAV application is widely adopted for line inspection of power and electronic industry in recent years, due to its advantage of low cost and high effectiveness.
Our application aligns with the growing utilization of Internet-of-Things technologies to enhance operational safety \citep{olsen2020industry,li2022after}.
The UAV line inspection mission consists of three phases, i.e., (i) the UAV leaves its base to the predetermined inspection task location; (ii) the UAV executes the inspection task; (iii) the UAV flies back to the base.
When the mission is aborted, the UAV halts the inspection task and directly heads back to the base.
Hence, we only need to focus on the optimal mission abort policy for the first two phases of the mission.
We conduct comprehensive simulations to assess the performance of our method.

\subsection{Parameter Setting}\label{subsec:num_setting}

In the simulation, the distance between the UAV base and power grid is 8 km.
During the mission, the UAV departs from the base and flies to the power grid at a constant speed 19.2 km/h \citep{liu2022energy}.
Upon arrival at the power grid, it executes the inspection task for 135 min.
Hence, the total duration for the first two phases of the mission is $H=60 \times 8/19.2 +135=160$ min, which is typical for a power-grid inspection task \citep{qiu2019gamma}. 
We set the inspection interval as $\delta=1$ min, so the number of decision epochs is $N=H/\delta=160$.
Based on this setting, we obtain the time duration (in min) of a rescue procedure after aborting a mission at the $n$th decision epoch:
\begin{equation*}\label{13}
	w_n =
	\begin{cases}
		n, & n\leq 25, \\
		25, & 25 < n \leq 160.
	\end{cases}
\end{equation*} 
The system failure cost is set as $C_{\text{s}}=2,000$ as with \cite{yang2019designing}.
The mission failure cost is more subjective.
We set $C_\text{m}=2,000$ and conduct sensitivity analysis on $C_{\text{m}}$ in Appendix~\ref{appen:sen}.

% The time distributions between the transition of states are set as follows.
In practice, the UAV failure during a mission is a rare event because the times $T_{12}$ and $T_{13}$ from healthy to defective state and to failure state are expected to be long.
If we directly use their distributions in simulations, the majority of missions would have no failure record, making it difficult to compare our method with benchmarks.
Hence, we modify the distributions of $T_{12}$ and $T_{13}$ to render more failures during simulations.
% This setting may correspond to the real scenario in which we use an aging UAV that have worked for several years to execute the mission.
Following \cite{khaleghei2021optimal}, we use an Erlang distribution for $T_{12}$.
Since an Erlang distribution corresponds to a CTMC, we do not need an Erlang mixture distribution for approximation.
This leads to a CTMC with a similar structure to our constructed CTMC in Section~\ref{subsec:pomdp_phase}, and the structural results hold for this CTMC; see Appendix~\ref{appen:num_setting} for details.
Next,
we set $\zeta=10^{-3}$ for $T_{13}$ such that $\mathbb{E}[T_{13}]/\mathbb{E}[T_{12}]$ is around three to four according to \cite{panagiotidou2010statistical}.
{We consider two different distributions of $T_{23}$.
In the first setting,
%To estimate the time distribution from defective to failure state, one can manually introduce some defect to a UAV in laboratory and then conduct a flying test to obtain the failure time data.
we use a Weibull distribution for $T_{23}$ with the following PDF:
\begin{align}\label{eq:weibull}
	f(t)=\frac{2.3}{108.8}\left(\frac{t}{108.8}\right)^{1.3}\exp\left(-(t/108.8)^{2.3}\right), \quad t\geq 0.
\end{align}
In the second setting, we assume $F$ follows a mixture distribution, which is commonly used in reliability \citep{krivtsov2007practical,yao2017shelf} and is given by}
\begin{align}\label{eq:mixture}
	{f(t)=\frac{1.3}{180.8}\left(\frac{t}{180.8}\right)^{1.6}\exp\left(-\left(\frac{t}{180.8}\right)^{2.6}\right)+\frac{1.15}{36.3}\left(\frac{t}{36.3}\right)^{1.3}\exp\left(-\left(\frac{t}{36.3}\right)^{2.3}\right).}
\end{align}
for $t\geq 0.$
Last, we determine the state-observation matrix $\mathbf{D}$.
According to \cite{wang2019multivariate}, the UAV flight control system typically monitors five physical states as shown in Table~\ref{table:parameters}.
During flight, sensors on the UAV measure the in-situ values of these physical states, and an alarm system computes the differences between the measured and desired values.
The differences for all the five physical states should be close to zero when the UAV is at the healthy state and deviate from zero at the defective state.
We obtain a dataset about the differences of these physical states when the UAV is in a healthy and defective state, respectively.
The differences of the five physical states lead to a five-dimensional signal.
To use our proposed model, we fit a logistic regression model by treating the five-dimensional signal as predictors and the system state (healthy or defective) as label.
A green (red) light condition-monitoring signal is observed when the predicted probability of being defective is smaller (larger) than $0.5$, as commonly adopted in practice.
The state-observation matrix $\mathbf{D}$ is then estimated as the empirical probability of observing green and red light signals under each latent state, given by
$$
\mathbf{D}=\left(
\begin{array}{cc}
	0.737 & 0.263\\
	0.101 & 0.899\\
\end{array}
\right).
$$

\begin{table}
	\centering
	\caption{Five key physical states of a UAV under condition monitoring \citep{wang2019multivariate}.}
	\label{table:parameters}
	\begin{tabular}{cc}
		\hline
		Description & Units \\
		\hline
		X-axis velocity & m/s\\
		Y-axis velocity &  m/s \\
		Z-axis velocity & m/s \\
		Pitch rate & degree/s\\
		Yaw rate & degree/s \\
		\hline
	\end{tabular}%
\end{table}%

%we follow \cite{zhu2021robust} and adopt principal component analysis to compress the multivariate signals.
%The first principle component accounting for the largest variability of data is used for condition monitoring.

%\ye{[this is strange. $K$-means is unsupervised. If you have the healthy/defective label for each signal, then you should train a linear classifier and then determine the decision boundary. This decision boundary can be used to determine which signals are green and which are red.] We aggregate the univariate signals under both the healthy and defective states, and then discretize the signals into $K=2$ levels by $K$-means.}

\subsection{Performance of the Model}\label{subsec:num_result}
\begin{figure}[htbp]
	\centering
	\subfigure{\includegraphics[scale=0.4,clip]{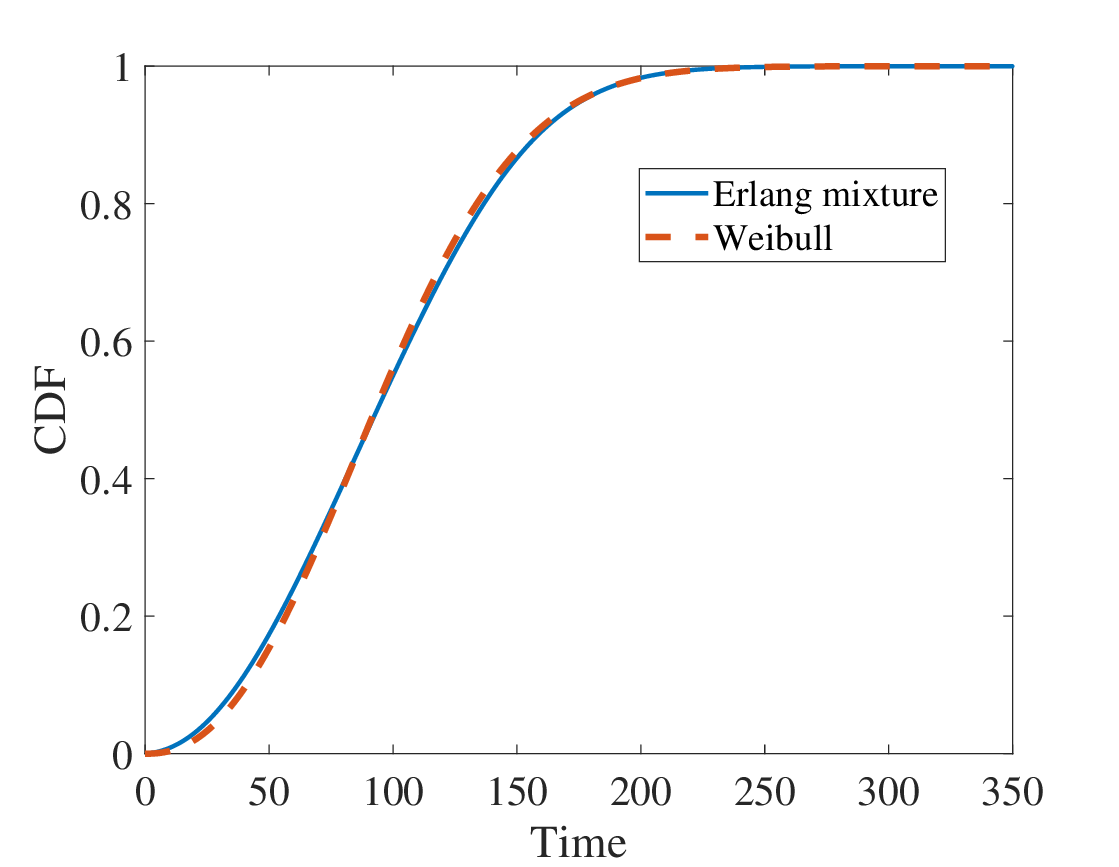}}
	\subfigure{\includegraphics[scale=0.4,clip]{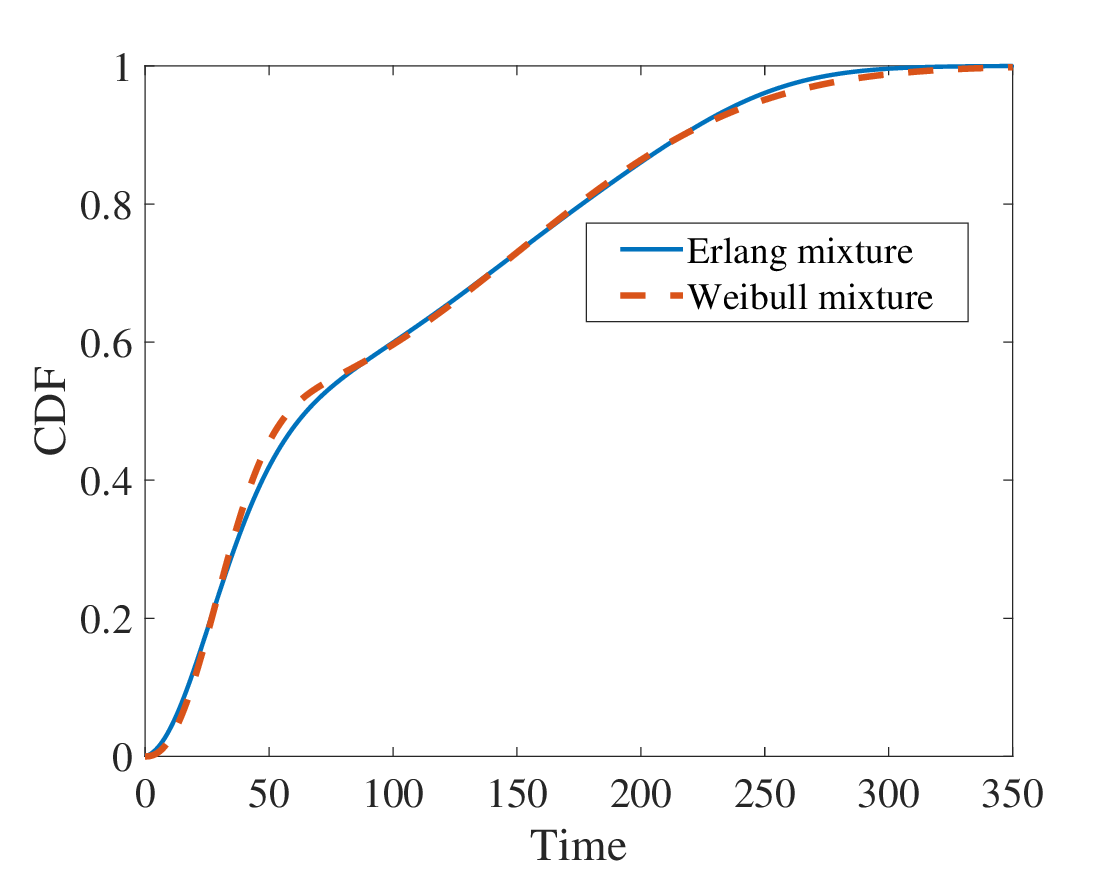}}
	\caption{Approximating the distribution of $T_{23}$ by the Erlang mixture distribution $F^{(\lambda)}(\cdot)$: The left panel uses $m_2= 20$ and $\lambda=0.134$ to approximate the Weibull distribution in \eqref{eq:weibull}.
	The right panel uses $m_2=50$ and $\lambda=0.209$ to approximate the mixture distribution in \eqref{eq:mixture}.} 
	\label{fig:cdf-compare}
\end{figure}

We simulate the system dynamics $\{X(t),~t\geq 0\}$ and $Y_{1:n}$ based on the parameter setting in Section~\ref{subsec:num_setting} and then use our method for abort decision-making.
To approximate the distribution $F(\cdot)$ by an Erlang mixture in \eqref{eq:herlang}. we first fix $m_2$ and use the commonly used moment-matching method to determine $\lambda$ such that the resulting Erlang mixture distribution has the same mean as the original distributions.
We enumerate $m_2$ and compute the approximation error defined by $\max_{0\leq t\leq H+w_N} |F(t)-F^{(\lambda)}(t)|$.
By the elbow method that seeks the point where the rate of decrease sharply changes, $m_2= 20$ with corresponding $\lambda=0.134$ is chosen when $T_{23}$ is Weibull in \eqref{eq:weibull}, and $m_2=50$ is chosen when $T_{23}$ follows the mixture distribution in \eqref{eq:mixture} with corresponding $\lambda=0.209$.
Figure~\ref{fig:cdf-compare} displays a comparison between the CDFs of $T_{23}$ and $T_{23}^{(\lambda)}$ {in the two settings}, where we see $F^{(\lambda)}(\cdot)$ well approximates $F(\cdot)$.
For the Weibull distribution, we verify that the constructed CTMC satisfies the monotone property in Lemma~\ref{lemma:increase_rate}, supporting our choice of $m_2$.
{Appendix~\ref{appen:sen} further examines the sensitivity of the optimal policy with respect to $(m_2,\lambda)$ and demonstrates how to select $(m_2,\lambda)$.}
The abort policy is optimized offline by Algorithm~\ref{algorithm:PBVI}, with hyper-parameters $Z_1=10^5$, $Z_2=5\times 10^5$, $L_1=10^4$, $\varepsilon=10^{-3}$, and $W=5,000$.

For comparison, we consider {four} benchmarks for mission abort decision-making.
%The first benchmark does not make any abort decision during mission.
%This benchmark corresponds to the case in which condition monitoring is unavailable.
The first two are {{rule}-based heuristics.
The first heuristic (referenced as ``$\mathcal{C}$-policy'') uses a control-chart scheme that deems the system as defective when receiving $\check{m}$ warning signals in $\check{N}$ consecutive periods.
The mission is aborted once the system is deemed defective.
We optimize both $\check{m}$ and $\check{N}$ by minimizing the expected operational cost.
This rule-based heuristic generalizes the method in \cite{yang2019designing}.
The second heuristic (referenced as ``$\mathcal{R}$-policy'') uses the system {RUL} to make decisions, and the mission is aborted when the $p$th percentile of the predicted RUL is smaller than the remaining time to complete the mission.
We optimize $p\in(0,100)$ for minimizing the expected operational cost.}
The third benchmark (referenced as ``$\mathcal{M}$-policy'') uses a three-state CTMC to approximate the failure process.
The transition rate matrix of this CTMC is determined such that the transition times between healthy, defective, and failure states have the same mean values as the original semi-Markov failure process.
A similar approximation is used in \cite{iravani2002integrated} for maintenance optimization.
Since this approximation uses a three-state CTMC, Corollary~\ref{thm:control_limit_CTMC} implies that the optimal $\mathcal{M}$-policy has a simple control-limit structure.
{The fourth benchmark keeps $T_{12}$ as $m_1$-phase Erlang distributed but sets $m_2=1$ for approximating $T_{23}$.
This benchmark degenerates to $\mathcal{M}$-policy when we enforce $m_1=1$.
Our dimension reduction result in Proposition~\ref{prop:dim_red} can be used to exactly solve the POMDP for the one-phase approximation.
%\bfred{[To discuss: if we need to compare our method with an EWMA chart? In healthcare papers, people generally compare their model with a benchmark commonly used in industry to gain insights. The benchmark commonly used in industry for our setting might be a control chart.]}
}
Details on optimizing these benchmarks are relegated to Appendix~\ref{appen:implement}.

{We first compare the model performance when $F$ is Weibull.}
Table~\ref{tab:cost_compare} summarizes results for our model and the four benchmarks under 10,000 simulation replications.
We can see that our method is consistently better than the benchmarks in term of smaller costs, larger mission success probability, and smaller system failure probability.
{Compared with the proposed method, the relative increase of the operational cost for the optimal $\mathcal{C}$-policy, $\mathcal{R}$-policy, $\mathcal{M}$-policy, and one-phase approximation is $4.89\%$, $7.52\%$, $4.93\%$, and $4.74\%$, respectively.
We additionally conduct a sensitivity analysis in Appendix~\ref{appen:sen} by varying $C_\text{m}$, $\mathbf{D}$, and $m_2$.
%\ye{[xxx]}.
It reveals that our method consistently outperforms the benchmarks, and the benchmarks incur a cost increase by around 2\%--13\% compared to our model.}
We shall highlight that since the optimal policy of our model is solved offline, it can be implemented online in real time. 

{We next consider the mixture distribution in \eqref{eq:mixture}.
The results averaged across 10,000 simulation replications are summarized in Table~\ref{tab:cost_compare_bimodal}. 
Our method again outperforms the benchmarks in terms of costs, mission success probability, and system survival probability.
The relative increase of the operational cost for the $\mathcal{C}$-policy, $\mathcal{R}$-policy, $\mathcal{M}$-policy, and one-phase approximation is $15.83\%$, $11.84\%$, $16.02\%$, and $15.96\%$, respectively.
Compared with the results when $F$ is Weibull, the gaps between the benchmarks and our method are significantly larger.
This outcome is unsurprising since exponential approximation to a bimodal density typically lacks accuracy. 
In contrast, a mixture of Erlang distribution offers the flexibility to approximate any $F$, resulting in significantly enhanced policies through subsequent POMDP approximation.
The sensitivity analysis in Appendix~\ref{appen:sen} further shows that cost increases of the benchmarks generally range from 5\% to 18\% compared to our method under various $C_\text{m}$ and $\mathbf{D}$.}

\begin{table}
	\centering
	\caption{A comparison of the cost, the mission success probability and the system failure probability between the proposed model and the four benchmarks when $T_{23}$ is Weibull in \eqref{eq:weibull}. The results are computed based on 10,000 Monte Carlo replications.}
	\label{tab:cost_compare}
	\begin{tabular}{cccc}
		\hline	
		Policy & Cost per mission & Mission success prob. & System failure prob. \\
		\hline
		$\mathcal{C}$-policy & 1063.0 & 0.668 & 0.198  \\
		$\mathcal{R}$-policy & 1089.6 & 0.662 & 0.207  \\
		%\midrule
		$\mathcal{M}$-policy & 1063.4 & 0.666 & 0.198 \\
		One-phase approximation & 1061.4 & 0.670 & 0.201 \\
		%\midrule
		Proposed model & 1013.4 & 0.681 & 0.188 \\
		\hline
	\end{tabular}%
\end{table}%

\begin{table}
	\centering
	\caption{
	A comparison of the cost, the mission success probability and the system failure probability between the proposed model and the four benchmarks when $T_{23}$ is a mixture distribution in \eqref{eq:mixture}. 
	The results are computed based on 10,000 Monte Carlo replications.}
	\label{tab:cost_compare_bimodal}
	\begin{tabular}{cccc}
		\hline	
		Policy & Cost per mission & Mission success prob. & System failure prob. \\
		\hline
		$\mathcal{C}$-policy & 1293.1 & 0.621 & 0.2675 \\
		$\mathcal{R}$-policy & 1248.6 & 0.6697 & 0.294 \\
		%\midrule
		$\mathcal{M}$-policy & 1295.2 & 0.6242 & 0.2718 \\
		One-phase approximation & 1294.6 & 0.6241 & 0.2714 \\
		%\midrule
		Proposed model & 1116.4 & 0.6881 & 0.2463 \\
		\hline
	\end{tabular}%
\end{table}%

{We see an interesting and intuitive result by comparing Tables~\ref{tab:cost_compare} and \ref{tab:cost_compare_bimodal}.
When $F$ is Weibull, an exponential distribution provides a relatively accurate approximation of $F$.
Hence, the one-phase approximation outperforms all other benchmarks.
However, when $T_{23}$ is a mixture distribution, both the $\mathcal{M}$-policy and one-phase approximation exhibit significant deterioration in performance, as an exponential approximation to a bimodal density is inaccurate. 
In this scenario, the $\mathcal{R}$-policy excels among all benchmarks, as it effectively captures the bimodality of $T_{23}$ in the RUL prediction.
There is no benchmark that consistently outperforms others under the two distributions in \eqref{eq:weibull} and \eqref{eq:mixture}.
In contrast, our model uses the mixture of Erlang distribution to enhance advantages of both rule-based methods and POMDPs.
It performs well even when the optimal mission abort policy is highly unstructured due to the bimodality of time from the defective to failure state.}
%\ye{[the one-phase has cost 907 and the M-policy has cost 908. how do you arrive at this conclusion? I would suggest deleting this part highlighted in red] However, a counter-intuitive result may be that the optimal $\mathcal{M}$-policy performs better than the one-phase approximation when $F$ is a mixture distribution.
%This phenomenon is not observed when $F$ is Weibull.
%Our results suggest that increasing the model complexity does not necessarily lead to a better performance when some parts of the system dynamics are oversimplified.}

\subsection{Performance of the Modified PBVI}
We compare our modified PBVI algorithm with the classical PBVI algorithm. 
{} 
Due to the curse of dimensionality, the value function {for the problem instance in Section~\ref{subsec:num_setting} cannot be solved exactly.}
For comparison, we first run the classical PBVI algorithm until hitting its termination condition.
We then run our modified PBVI algorithm until the corresponding approximate value function is close to that obtained from the classical PBVI algorithm (see Appendix~\ref{appen:num_setting} for the detailed rule to terminate the modified PBVI algorithm). % \ye{[see this for what? }
We replicate the above procedure 10 times and compare the running times.
The mean times to finish the modified and classical PBVI algorithms are 5.75 and 63.45 hours, respectively. 
The runtime of our modified PBVI algorithm is less than 10\% of that for the classical PBVI algorithm.

{
To get a comprehensive insight into the runtime and accuracy of both algorithms, we create a small problem instance with $m_1=1$ and $m_2=2$, leading to a belief state $\{\bm\pi\in\mathbb{R}^3_+:\pi_1+\pi_2+\pi_3=1\}$ before system failure.
Other parameters are kept the same as those in Section~\ref{subsec:num_setting}.
% \ye{[how many states? what is the problem? it deserves a one-sentence intro before you throw everything to the appendix]}. 
This instance represents possibly the maximum size that can be solved through belief space discretization and backward induction within a reasonable timeframe.
The detailed parameter setting is given in Appendix~\ref{appen:num_setting}.
The availability of the exact value function allows us to compute the optimality gaps of the two PBVI algorithms.
The running times (averaged by 20 simulations) for the two PBVI algorithms to achieve solutions with different optimality gaps are summarized in Table~\ref{tab:run_time}.
The two PBVI algorithms implement the same procedures at early stages, so their performances in this stage are comparable.
In contrast, the modified PBVI algorithm becomes much faster in the presence of a small optimality gap. 
When the optimality gap is 0.3\%, the running time of the modified PBVI algorithm is only around 3.5\% of that for the classical PBVI algorithm.

\begin{table}[htbp]
	\centering
	\caption{Runtime comparison (in minutes) for the classical PBVI and the proposed modified PBVI algorithms to achieve solutions with different optimality gaps.}
	\label{tab:run_time}
	\begin{tabular}{ccccccc}
		\hline
		Optimality gap & 15\% & 10\% & 5\% & 1\% & 0.5\% & 0.3\% \\
		\hline
		Classical PBVI & 0.37 & 0.66 & 0.86 & 161.86 & 1712.08 & 9803.31 \\
		Modified PBVI & 0.34 & 0.62 & 0.81 & 92.28 & 229.30 & 344.00 \\
		%\midrule
		\hline
	\end{tabular}%
\end{table}%
}

\subsection{Mission Abort in a Multi-Task Setting}\label{subsec:num_multi}
{We conduct numerical experiments in a multi-task setting as described in Section~\ref{sec:multiple}.
We set $L=3$, and the $L=3$ tasks need $N_1=35$, $N_2=50$, and $N_3=50$ periods to complete, respectively.
% The inter-sampling interval is $\delta=1$.
The mission failure costs for each task are $(C_\text{m}^{(1)},C_\text{m}^{(2)},C_\text{m}^{(3)})=(500,300,200)$.
The repair cost is $C_\text{r}=1,000$.
The rescue times $(w_n)_{n\in[N]}$ are specified in Appendix~\ref{appen:num_setting}.
For illustration, we only conduct experiments when the distribution $F$ follows \eqref{eq:mixture}.
Other settings such as the inter-sampling interval, system failure cost, and observation matrix are the same as those in Section~\ref{subsec:num_setting}.
% We compare our model with the benchmarks in Section~\ref{subsec:num_result}.
Using 1,000 simulation replications, the operational costs of the proposed method, $\mathcal{C}$-policy, $\mathcal{R}$-policy, $\mathcal{M}$-policy, and one-phase approximation are 893.6, 1002.5, 960.2, 1028.7, and 1026.7, respectively.
Comparing these numbers with our model, we find that the $\mathcal{C}$-policy, $\mathcal{R}$-policy, $\mathcal{M}$-policy, and one-phase approximation increase the operational cost by 12.19\%, 7.45\%, 15.12\%, and 14.89\%, respectively.
This reaffirms the cost-effectiveness of our method.}

\iffalse
\begin{table}[htbp]
	\centering
	\caption{A comparison of the cost, the mission success probability and the system failure probability between the proposed model and the four benchmarks when $T_{23}$ follows \eqref{eq:mixture}, under the multi-task setting. 
	The results are computed based on 10,000 Monte Carlo replications.}
	\label{tab:cost_compare_multi}
	\begin{tabular}{ccc}
		\hline	
		Policy & Cost per mission & Relative increase in cost \\
		\hline
		$\mathcal{C}$-policy & 1002.5 & 0.2357 \\
		$\mathcal{R}$-policy & 960.2 & 0.267 \\
		%\midrule
		$\mathcal{M}$-policy & 1028.7 & 0.1806 \\
		One-phase approximation & 1026.7 & 0.1738 \\
		%\midrule
		Proposed model & 893.6 & 0 \\
		\hline
	\end{tabular}%
\end{table}%
\fi

\section{Conclusions}\label{sec:conclusion}

This study has developed a general surrogate POMDP framework for mission abort decisions of a three-state mission-critical system with imperfect condition-monitoring signals.
Although the original problem \eqref{eq:stopping} is intractable and the optimal mission abort policy is expected to be highly unstructured, the optimal solution of the surrogate POMDP, with proper design in Figure~\ref{fig:CTMC}, enjoys nice structural properties under a spherical coordinate system.
In particular, mission abort is never optimal at the late stage of the mission;
before this stage and given the belief state $(r_n,\bm\Phi_n)$ under the spherical coordinate system at the $n$th decision epoch, the optimal abort policy has a simple upper-lower control-limit structure with respect to the radius $r_n$, but the control limits depend on the angle $\bm\Phi_n$.
These properties enabled us to design a customized algorithm that first finds a cut-point for the late stage and then develop a modified PBVI algorithm to deal with the curse of dimensionality.
When $w_n=0$ for all $n\in[N]$, mission abort resembles the CBM setting \citep{zhu2021multicomponent,sun2023robust}, which implies that our framework is applicable to CBM.
The case study on a UAV reveals around 5\%--20\% reduction in the operational cost of the proposed method compared with rule-based policies and simple CTMC approximations under various parameter settings.

Beyond mission abort, the proposed surrogate POMDP framework applies to other discrete-time optimal stopping problems with bounded costs and .
For example, it can be used to handle non-Markovian transition of the patient health state in treatment planning \citep{skandari2021patient}.
The framework is also applicable to some sequential decision-making problems, such as inventory control with inter-restocking intervals following a phase-type distribution \citep{wang2023inventory}.
% Similar to the false and missed alarms in mission abort, the false positive and negative diagnoses are common in healthcare literature \citep{tian2021optimal}. Our framework that factors in measurement errors has potential in improving reliability of sequential medical tests for management of pandemic like COVID-19. Last, our modified PBVI algorithm leverages the convex structure of the optimal abort policy. This strategy can be applied to other ADP problems if the corresponding optimal policy has a similar structural property.
Some natural extensions are of interest for further investigation.
{
If we specialize in mission abort scenarios for UAVs, range constraints due to battery capacity can be a potential issue.  
Consequently, location designs for rescue sites, where a UAV can safely land for rescue missions, becomes crucial.
In a general mission abort problem, uncertainty in the transition and state-observation matrices can be considered.
Promising approaches include Bayesian methods and robust POMDPs.
Structural properties of our problem may be leveraged to speed up computation in these more complex cases.
% In addition to the binary abort decision, we can also control the mission execution process. For example, we can dynamically adjust the concentration of chemicals in a chemical reactor to control the production process according to the cooling system performance. Joint optimization of the pre-abort control and mission abort can be an interesting topic.
When historical data are available to infer these parameters, a prominent focus in recent literature has been the development of an end-to-end data-driven methodology that directly maps the available data to a control policy.
Our approximation scheme provides a viable framework for this purpose with the Erlang mixture distribution as a key intermediary. 
Our Theorem~\ref{thm:convergence} has shown that given the best approximating Erlang mixture distribution with a given phase, the optimal solution closely approximates the true optimal policy.
The next step is to devise a statistical procedure to estimate this best approximating Erlang mixture distribution.
This estimation is closely related to the method of sieve in statistics \citep{chen2007large}.
With a proper data-driven choice of the sieve, the estimated Erlang mixture distribution and the estimated policy are expected to be asymptotically consistent, which deserve future research.}

\section*{Acknowledgment}
This work was supported in part by the National Science Foundation of China (72071071, 72171037, 72371161, 72471144), Singapore MOE AcRF Tier 2 grant (A-8001052-00-00, A-8002472-00-00), and the Future Resilient Systems project supported by the National Research Foundation Singapore under its CREATE programme.
The authors would like to thank the editor, an associate editor, and two anonymous referees for their valuable comments and constructive suggestions, which have led to a significant improvement in the quality and presentation of this work.

\section*{Author Biographies}
\begin{itemize}
	\item Qiuzhuang Sun is a Lecturer in the School of Mathematics and Statistics at the University of Sydney. His research focuses on data-driven decision-making, applied data science, and system reliability.
	\item Jiawen Hu is an Associate Professor with the School of Aeronautics and Astronautics, University of Electronic Science and Technology of China. His current research interests include maintenance optimization and degradation modeling.
	\item Zhi-Sheng Ye is an Associate Professor and the Dean's Chair in the Department of Industrial Systems Engineering \& Management at the National University of Singapore. His research focuses on data-driven operations management, mathematical and industrial statistics, and system resilience.
\end{itemize}

\bibliographystyle{informs2014} % outcomment this and next line in Case 1
\bibliography{ref} % if more than one, comma separated

\newpage
\begin{APPENDICES}

\noindent
\begin{center}
	{\Large{Electronic Companion to ``Optimal Abort Policy for Mission-Critical Systems under Imperfect Condition Monitoring''}}
\end{center}

\vspace{1em}

{\centering
Qiuzhuang Sun,$^1$ Jiawen Hu,$^2$ and Zhi-Sheng Ye$^3$ 

{\scriptsize
$^1$School of Mathematics and Statistics, University of Sydney, Australia

$^2$School of Aeronautics and Astronautics, University of Electronic Science and Technology of China, China

$^3$Department of Industrial Systems Engineering and Management, National University of Singapore, Singapore

Emails: \texttt{qiuzhuang.sun@sydney.edu.au}; \texttt{hdl@sjtu.edu.cn}; \texttt{yez@nus.edu.sg}}

}

% \vspace{1em}

\setcounter{figure}{0}  
\renewcommand{\thefigure}{EC.\arabic{figure}}
\setcounter{table}{0}  
\renewcommand{\thetable}{EC.\arabic{table}}
\setcounter{equation}{0}  
\renewcommand{\theequation}{EC.\arabic{equation}}
\setcounter{lemma}{0}  
\renewcommand{\thelemma}{EC.\arabic{lemma}}
\setcounter{theorem}{0}  
\renewcommand{\thetheorem}{EC.\arabic{theorem}}
\setcounter{proposition}{0}  
\renewcommand{\theproposition}{EC.\arabic{proposition}}
\setcounter{algorithm}{0}  
\renewcommand{\thealgorithm}{EC.\arabic{algorithm}}
\setcounter{page}{1} 
\renewcommand{\thepage}{EC.\arabic{page}}

\section{Technical Proofs}\label{appen:proof}

\subsection{Proof of Proposition~\ref{prop:absorption}}
Consider the continuous-time Markov chain (CTMC) $\{X^{(\lambda)}(t):t\geq0\}$ with initial probabilities in \eqref{eq:init_prob} and  transition rate matrix $\mathbf{Q}$ in \eqref{eq:trans_prob}.
If ${X}^{(\lambda)}$ starts from state $i$, $i\in[m_1]$, 
the transition time $T_{12}^{(\lambda)}$ from this state to state $m_1+1$ follows an Erlang distribution with rate $\lambda-\zeta$ and shape $m_1-i+1$.
Since ${X}^{(\lambda)}$ starts start from $m_1$ with probability $\pi_{01}^{(\lambda)} = 1-G((m_1-1)/(\lambda-\zeta))$ and from $i\in\{2,\ldots,m_1\}$ with probability $\pi_{0i}^{(\lambda)} = G((m_1+1-i)/(\lambda-\zeta))-G((m_1-i)/(\lambda-\zeta))$,
$T_{12}^{(\lambda)}$ is a mixed Erlang random variable unconditionally.
Its distribution can be specified from the initial probability $\pi_{0i}^{(\lambda)}$ and the Erlang distributions above, which is exactly $G^{(\lambda)}$ in \eqref{eq:herlang2}.

%An Erlang random variable with shape $m_1$ and rate $\nu$ has the same distribution as the sum of $m_1$ independent exponential random variables with rate $\nu$.
%Based on this fact and the memoryless property of exponential distributions, it is easy to see that the hitting time of ${X}^{(\lambda)}(t)$ to the states in $\mathcal{X}_2$ has the same distribution as $T_{12}$.

A careful inspection of $q_{ij}$ shows that for any state $i\in\mathcal{X}_1$, there is a  transition rate $q_{i,m_1+m_2+1}=\zeta$ to the absorbing state $m_1+m_2+1$.
Since $q_{i,m_1+m_2+1}$ are identical for all $i\in\mathcal{X}_1$, the (potential) transition time, denoted as $T_{13}^{(\lambda)}$, follows an exponential distribution with rate $\zeta$, same as $T_{13}$.
By the memoryless property of exponential distributions, we see that the failure time of ${X}^{(\lambda)}$ has the same distribution as $T_{13}^{(\lambda)}$ when $T_{13}^{(\lambda)}\leq T_{12}^{(\lambda)}$.
On the other hand, we derive the system failure time when $T_{13}^{(\lambda)}>T_{12}^{(\lambda)}$ as follows.
When $T_{13}^{(\lambda)}>T_{12}^{(\lambda)}$, the system is at state $m_1+1$ at time $T_{12}^{(\lambda)}$.
From this time on, the probability that ${X}^{(\lambda)}(t)$ enters the absorbing state $m_1+m_2+1$ with additional $l$ transitions is
$$
p_l(\lambda)\prod_{i=1}^{l-1}(1-p_i(\lambda))=\frac{F(l/\lambda)-F((l-1)/\lambda)}{1-F((l-1)/\lambda)}\prod_{i=1}^{l-1}\frac{1-F(i/\lambda)}{1-F((i-1)/\lambda)}=F\left(\frac{l}{\lambda}\right)-F\left(\frac{l-1}{\lambda}\right)
$$
for $l\in[m_2-1]$,
and 
$$
\prod_{i=1}^{m_2-1}(1-p_i(\lambda))=\prod_{i=1}^{m_2-1}\dfrac{1-F(i/\lambda)}{1-F((i-1)/\lambda)}=1-F\left(\frac{m_2-1}{\lambda}\right)
$$
for $l=m_2$, where we use the fact that $F(0)=0$.
When ${X}^{(\lambda)}$ uses $l$ steps to transition from state $m_1+1$ to state $m_1+m_2+1$, the elapsed time follows an Erlang distribution with shape $l$ and rate $\lambda>0$ by our construction.
This fact, in conjunction with the probability of $l$ transitions to the absorbing state specified above, shows that the duration from time $T_{12}^{(\lambda)}$ to the absorbing state follows a mixed Erlang distribution.
It is easy to check that this mixed Erlang distribution is exactly the one specified in Equation~\eqref{eq:herlang} for $T_{23}^{(\lambda)}$.
Combining the above analysis to see that the absorption time of ${X}^{(\lambda)}$ has the same distribution as $T_{13}^{(\lambda)}\mathbbm{1}{\{T_{13}^{(\lambda)}\leq T_{12}^{(\lambda)}\}}+(T_{12}^{(\lambda)}+T_{23}^{(\lambda)})\mathbbm{1}{\{T_{13}^{(\lambda)}> T_{12}^{(\lambda)}\}}$,
whence the proposition holds.
\Halmos

\subsection{Proof of Theorem~\ref{thm:convergence}}

\iffalse
Let $\mu$ be the dominating measure for the joint distribution of $\{T_{12},T_{13},T_{23},Y_{1:N}\}$ and $\{T_{12}^{(\lambda)},T_{13}^{(\lambda)},T_{23}^{(\lambda)},Y_{1:N}^{(\lambda)}\}$, and let $p()$ and $p^{(\lambda)}()$ be the density of $\{T_{12},T_{13},T_{23},Y_{1:N}\}$ and $\{T_{12}^{(\lambda)},T_{13}^{(\lambda)},T_{23}^{(\lambda)},Y_{1:N}^{(\lambda)}\}$ under $\mu$.
Before the proof, we present a lemma that will be useful later on.
\begin{lemma}
	
\end{lemma}
Idea for the proof.
Show $p^{(\lambda)}(t_{12}, t_{13},t_{23}) \to p(t_{12}, t_{13},t_{23})$.
\fi

Given that the realization of the condition-monitoring signals $Y_{1:n}=y_{1:n}\in\{0,\ldots,K\}^n$, we let $\mathbb{C}_{\text{ab}}(n, y_{1:n})$ and $\mathbb{C}_{\text{c}}(n, y_{1:n})$ denote the minimum costs-to-go for $X(t)$ when the mission is aborted or continued at time $n\delta$, respectively.
Similarly, given the same information, let $\mathbb{C}_{\text{ab},\lambda}(n, y_{1:n})\triangleq V^{(\lambda)}_\text{ab}(n,\bm\pi_n^{(\lambda)})$ and $\mathbb{C}_{\text{c},\lambda}(n, y_{1:n})\triangleq V^{(\lambda)}_\text{c}(n,\bm\pi_n^{(\lambda)})$ be the minimum costs-to-go for $X^{(\lambda)}(t)$, where $V_\text{ab}$ and $V_\text{c}$ are defined in Section~\ref{subsec:pomdp_model}, and here we explicitly indicate their dependence on $\lambda$.
The optimal decision rules $a_n^*$ and $a_{\lambda,n}^*$ can be determined by comparing the two costs-to-go, i.e., $a_n^*(y_{1:n})=\mathbbm{1}{\left\{\mathbb C_{\text{ab}}(n,y_{1:n}) < \mathbb C_{\text{c}}(n,y_{1:n})\right\}}$ and $a_{\lambda,n}^*(y_{1:n})=\mathbbm{1}{\left\{\mathbb C_{\text{ab},\lambda}(n,y_{1:n}) < \mathbb C_{\text{c},\lambda}(n,y_{1:n})\right\}}$.
To show $a_{\lambda,n}^*\to a_n^*$ $\mathbb{P}$-almost surely as $\lambda\to\infty$, it suffices to show $\lim_{\lambda\to\infty}\mathbb C_{\text{ab},\lambda}(n,y_{1:n})=\mathbb C_{\text{ab}}(n,y_{1:n})$ and $\lim_{\lambda\to\infty}\mathbb C_{\text{c},\lambda}(n,y_{1:n})=\mathbb C_{\text{c}}(n,y_{1:n})$ for all $y_{1:n}\in\{0,\ldots,K\}^n$ and $n\in[N]$.

Recall that $T_{12}^{(\lambda)},T_{23}^{(\lambda)},T_{13}^{(\lambda)}$ are the independent (potential) transition times from state clusters $\mathcal{X}_1$ to $\mathcal{X}_2$, from $\mathcal{X}_2$ to $\mathcal{X}_3$, and from $\mathcal{X}_1$ to $\mathcal{X}_3$, respectively.
We then show $(T_{12}^{(\lambda)},T_{13}^{(\lambda)},T_{23}^{(\lambda)},Y^{(\lambda)}_{1:N})$ jointly converges to $(T_{12},T_{13},T_{23},Y_{1:N})$ 
in distribution as $\lambda\to\infty$.
Let 
$\boldsymbol{T}^{(\lambda)}=(T_{12}^{(\lambda)},T_{13}^{(\lambda)},T_{23}^{(\lambda)})$, 
$\boldsymbol{t}=(t_{12},t_{13},t_{23})$, 
$\boldsymbol{u}=(u_{12},u_{13},u_{23})$, and 
$\boldsymbol{F}^{(\lambda)}(\boldsymbol{t})=P^{(\lambda)}(T_{12}^{(\lambda)}\leq t_{12},T_{13}^{(\lambda)}\leq t_{13},T_{23}^{(\lambda)}\leq t_{23})$.
By the independence of $Y_1^{(\lambda)},\ldots,Y_n^{(\lambda)}$ conditional on the failure process $X^{(\lambda)}(t)$, we have
\begin{align*}
P^{(\lambda)}(\boldsymbol{T}^{(\lambda)}\leq \boldsymbol{t},~ Y^{(\lambda)}_{1:N}=y_{1:N})
=&\int_{\bm u\leq \boldsymbol{t}} P^{(\lambda)}(Y^{(\lambda)}_{1:N}=y_{1:N} \mid \boldsymbol{T}^{(\lambda)} = \boldsymbol{u})\td \boldsymbol{F}^{(\lambda)}(\boldsymbol{u})\\
=&\int_{\bm u \leq \boldsymbol{t}}\prod_{n=1}^N P^{(\lambda)}(Y^{(\lambda)}_{n}=y_n \mid \boldsymbol{T}^{(\lambda)} = \boldsymbol{u})\td \boldsymbol{F}^{(\lambda)}(\boldsymbol{u}).
\end{align*}
Similarly,
\begin{align*}
\mathbb P(\boldsymbol{T}\leq \boldsymbol{t},~ Y_{1:N}=y_{1:N})
=\int_{\bm u \leq \boldsymbol{t}} \prod_{n=1}^N \mathbb P (Y_{n}=y_n \mid \boldsymbol{T} = \boldsymbol{u})\td \boldsymbol{F}(\boldsymbol{u}),
\end{align*}
where $\boldsymbol{T}=(T_{12},T_{13},T_{23})$ and $\boldsymbol{F}(\boldsymbol{t})=\mathbb P(T_{12}\leq t_{12},T_{13}\leq t_{13},T_{23}\leq t_{23})$.
By our construction of the surrogate signals $Y^{(\lambda)}_{1:N}$, $P^{(\lambda)}(Y^{(\lambda)}_{n}=y_n \mid \boldsymbol{T}^{(\lambda)} = \boldsymbol{u})$ is equal to $\mathbb P(Y_{n}=y_n \mid \boldsymbol{T} = \boldsymbol{u})$.
If we can show that $g_n(\boldsymbol{u}) \triangleq \mathbb P(Y_{n}=y_n \mid \boldsymbol{T} = \boldsymbol{u})$ is continuous almost everywhere with respect to the distribution of $\boldsymbol{T}=(T_{12},T_{13},T_{23})$, then we can use the fact that $\boldsymbol{F}^{(\lambda)}\leadsto \boldsymbol{F}$ and apply Portmanteau's theorem to the function $g(\boldsymbol{u})\triangleq \prod_{n=1}^N g_n(\boldsymbol{u})\mathbbm{1}\{\bm u \leq \boldsymbol{t}\}$
to conclude that $P^{(\lambda)}(\boldsymbol{T}^{(\lambda)}\leq \boldsymbol{t}, ~Y^{(\lambda)}_{1:N}=y_{1:N})\to \mathbb P(\boldsymbol{T}\leq \boldsymbol{t}, ~Y_{1:N}=y_{1:N})$ as $\lambda\to\infty$.
Note that
\begin{align*}
\mathbb{P}(Y_n= y_n \mid \bm T= \bm u) &=
d_{1,y_n}\mathbbm{1}\{n\delta < \min\{u_1,u_2\}\}+
	d_{2,y_n}\mathbbm{1}\{u_1<n\delta \leq u_1+u_3,~u_1 < u_2\} \\
&\quad +	\mathbbm{1}\{y_n=0, ~n\delta \geq u_2\mathbbm{1}\{u_1\geq u_2\}+(u_1+u_3)\mathbbm{1}\{u_1< u_2\}\}.
\end{align*}
This function has finite number of discontinuous points.
It is continuous almost everywhere by the continuity assumption on $T_{12},T_{13},T_{23}$, whence Portmanteau's theorem applies. 
This implies that
$P^{(\lambda)}(\boldsymbol{T}^{(\lambda)}\leq \boldsymbol{t},~ Y^{(\lambda)}_{1:N}\leq y_{1:N})\to P(\boldsymbol{T}\leq \boldsymbol{t},~ Y_{1:N}\leq y_{1:N})$ as $\lambda\to\infty$.
That is, $(T_{12}^{(\lambda)},T_{13}^{(\lambda)},T_{23}^{(\lambda)},Y^{(\lambda)}_{1:N})$ converges to $(T_{12},T_{13},T_{23},Y_{1:N})$ in distribution as $\lambda\to\infty$.

With the weak convergence established above, we now show $\lim_{\lambda\to\infty}\mathbb C_{\text{ab},\lambda}(n,y_{1:n})=\mathbb C_{\text{ab}}(n,y_{1:n})$ and $\lim_{\lambda\to\infty}\mathbb C_{\text{c},\lambda}(n,y_{1:n})=\mathbb C_{\text{c}}(n,y_{1:n})$ for all $n\in[N]$ by induction.
For $n=N-1$, assume the system is working and the realization of conditioning monitoring signals is $Y_{1:N-1}=Y_{1:N-1}^{(\lambda)}=y_{1:N-1}$.
For the failure process $X^{(\lambda)}(t)$,
the cost-to-go of aborting the mission is 
$$
\mathbb{C}_{\text{ab},\lambda}(N-1, y_{1:N-1})= C_\text{m}+C_\text{s}{P}^{(\lambda)}(\xi^{(\lambda)}\leq (N-1)\delta+w_{N-1} \mid Y^{(\lambda)}_{1:N-1}=y_{1:N-1}).
$$
For the failure process $X(t)$,
the cost-to-go of aborting the mission is 
$$
\mathbb{C}_{\text{ab}}(N-1, y_{1:N-1}) = C_\text{m}+C_\text{s}\mathbb{P}(\xi\leq (N-1)\delta+w_{N-1} \mid Y_{1:N-1}=y_{1:N-1}).
$$
Since $(T_{12}^{(\lambda)},T_{13}^{(\lambda)},T_{23}^{(\lambda)},Y_{1:N}^{(\lambda)}) \rightsquigarrow (T_{12},T_{13},T_{23},Y_{1:N})$ as $\lambda\to\infty$ and all transition times in the same probability space are independent, we have
$$
\begin{aligned} 			
	& P^{(\lambda)}(\xi^{(\lambda)}\leq t\mid Y^{(\lambda)}_{1:N-1}=y_{1:N-1}) \\
	=~&P^{(\lambda)}(T_{13}^{(\lambda)}\mathbbm{1}{\{T_{13}^{(\lambda)}\leq T_{12}^{(\lambda)}\}}+(T_{12}^{(\lambda)}+T_{23}^{(\lambda)})\mathbbm{1}{\{T_{13}^{(\lambda)}> T_{12}^{(\lambda)}\}}\leq t\mid Y^{(\lambda)}_{1:N-1}=y_{1:N-1}) \\
	=~& \frac{P^{(\lambda)}(T_{13}^{(\lambda)}\mathbbm{1}{\{T_{13}^{(\lambda)}\leq T_{12}^{(\lambda)}\}}+(T_{12}^{(\lambda)}+T_{23}^{(\lambda)})\mathbbm{1}{\{T_{13}^{(\lambda)}> T_{12}^{(\lambda)}\}}\leq t,~Y^{(\lambda)}_{1:N-1}=y_{1:N-1})}
	{P^{(\lambda)}(Y^{(\lambda)}_{1:N-1}=y_{1:N-1})} \\
	\to ~& \frac{\mathbb P(T_{13}\mathbbm{1}{\{T_{13}\leq T_{12}\}}+(T_{12}+T_{23})\mathbbm{1}{\{T_{13}> T_{12}\}}\leq t,~Y_{1:N-1}=y_{1:N-1})}
	{\mathbb P(Y_{1:N-1}=y_{1:N-1})} \\
	=~& \mathbb{P}(\xi\leq t\mid Y_{1:N-1}=y_{1:N-1})
\end{aligned}
$$
as $\lambda\to \infty$ for all $t> 0$ by the continuous mapping theorem.
Hence, the cost-to-go $\mathbb{C}_{\text{ab},\lambda}(N-1, y_{1:N-1})$ of aborting the mission using the stochastic process ${X}^{(\lambda)}(\cdot)$ converges to $\mathbb{C}_{\text{ab}}(N-1, y_{1:N-1})$ based on $X(t)$ as $\lambda\to\infty$.
By the same reasoning, the costs-to-go $\mathbb{C}_{\text{c},\lambda}(N-1, y_{1:N-1})$ of continuing the mission under ${X}^{(\lambda)}(\cdot)$ also converges to $\mathbb{C}_{\text{c}}(N-1, y_{1:N-1})$.
By the assumption of unique set of optimal solutions of the original problem in \eqref{eq:stopping}, $\mathbb{C}_{\text{ab}}(N-1, y_{1:N-1})\neq \mathbb{C}_{\text{c}}(N-1, y_{1:N-1})$.
The uniqueness, in conjunction with the convergence above, implies that when $\lambda$ is large enough, the sign of $\mathbb{C}_{\text{ab}, \lambda}(N-1, y_{1:N-1})- \mathbb{C}_{\text{c}, \lambda}(N-1, y_{1:N-1})$ is the same as that of $\mathbb{C}_{\text{ab}}(N-1, y_{1:N-1})- \mathbb{C}_{\text{c}}(N-1, y_{1:N-1})$.
Consequently, $\lim_{\lambda\to\infty}\mathbb C_{\lambda}(N-1,y_{1:N-1})=\mathbb C(N-1,y_{1:N-1})$,
where $\mathbb C_{\lambda,n-1}(y_{1:n-1})\triangleq\min\{\mathbb C_{\textnormal{ab},\lambda}(n-1,y_{1:n-1}),\mathbb C_{\textnormal{c},\lambda}(n-1,y_{1:n-1})\}$ and $\mathbb C(n-1,y_{1:n-1})\triangleq\min\{\mathbb C_{\textnormal{ab}}(n-1,y_{1:n-1}),\mathbb C_{\textnormal{c}}(n-1,y_{1:n-1})\}$ for all $n\in[N]$.

Now assume $a_{\lambda,n}^*(Y_{1:n})\to a_n^*(Y_{1:n})$ $\mathbb{P}$-almost surely and $\mathbb C_{\lambda}(n,y_{1:n})\to \mathbb C(n,y_{1:n})$ as $\lambda\to\infty$ for all $y_{1:n}\in[K]^n$ and $n=m,\ldots,N-1$.
Then for $n=m-1$, we can repeat the above argument for $n=N-1$ to show that the costs-to-go for aborting the mission are the asymptotically same for the two stochastic processes $X(t)$ and ${X}^{(\lambda)}(t)$ as $\lambda\to\infty$.
Moreover, given $Y_{1:m-1}=y_{1:m-1}\in[K]^{m-1}$, the cost-to-go of continuing the mission for ${X}(t)$ is
\begin{align*}
\mathbb{C}_{\text{c}}(m-1, y_{1:m-1}) = \mathbb{E}\big[&(C_\text{m}+C_\text{s})\mathbbm{1}{\{\xi\leq m\delta\}}+\mathbb C_{\text{ab}}(m,Y_{1:m})\mathbbm{1}{\{\xi>m\delta,~a_m^*(Y_{1:m})=1 \}}\\
	&+\mathbb C_{\text{c}}(m,Y_{1:m})\mathbbm{1}{\{\xi>m\delta,~a_m^*(Y_{1:m})=0\}}\mid Y_{1:m-1}=y_{1:m-1}\big].
\end{align*}
Similarly, the cost-to-go of continuing the mission for ${X}^{(\lambda)}(t)$ is
\begin{align*}
\mathbb{C}_{\text{c}, \lambda}(m-1, y_{1:m-1}) = \mathbb{E}^{(\lambda)}
\big[&(C_\text{m}+C_\text{s})\mathbbm{1}{\{\xi^{(\lambda)}\leq m\delta\}}+ \mathbb C_{\text{ab},\lambda}(m,Y_{1:m}^{(\lambda)})\mathbbm{1}{\{\xi^{(\lambda)}>m\delta,~ a_{\lambda,m}^*(Y_{1:m}^{(\lambda)})=1 \}}\\
	&+\mathbb C_{\text{c},\lambda}(m,Y_{1:m}^{(\lambda)})\mathbbm{1}{\{\xi^{(\lambda)}>m\delta,~ a_{\lambda,m}^*(Y_{1:m}^{(\lambda)})=0\}}\mid Y_{1:m-1}^{(\lambda)}=y_{1:m-1}\big].
\end{align*}
The induction hypothesis implies that $\lim_{\lambda\to\infty}a_{\lambda,m}^*(y_{1:m}) = a_{m}^*(y_{1:m})$.
By the weak convergence of $(T_{12}^{(\lambda)},T_{13}^{(\lambda)},T_{23}^{(\lambda)},Y_{1:N}^{(\lambda)})$ to $(T_{12},T_{13},T_{23},Y_{1:N})$, it can be routinely verified that the conditional distributions of $\xi^{(\lambda)}$ and $Y_m^{(\lambda)}$ given $Y_{1:m-1}^{(\lambda)}=y_{1:m-1}$ converge to those of $\xi$ and $Y_m$ given $Y_{1:m-1}=y_{1:m-1}$, as $\lambda\to\infty$.
The random variables within the conditional expectations of $\mathbb{C}_{\text{c}, \lambda}(m-1, y_{1:m-1})$ above are bounded by $C_\text{m}+C_\text{s}$, and hence are uniformly integrable.
As such, weak convergence implies convergence of moments, and $\mathbb{C}_{\text{c}, \lambda}(m-1, y_{1:m-1})$ converges to $\mathbb{C}_{\text{c}}(m-1, y_{1:m-1})$ as $\lambda\to\infty$.
Again, by the assumption of unique set of optimal solutions of the original problem in \eqref{eq:stopping}, $\mathbb{C}_{\text{ab}}(m-1, y_{1:m-1})\neq \mathbb{C}_{\text{c}}(m-1, y_{1:m-1})$.
Therefore, the sign of $\mathbb{C}_{\text{ab}, \lambda}(m-1, y_{1:m-1})- \mathbb{C}_{\text{c}, \lambda}(m-1, y_{1:m-1})$ will eventually be the same as that of $\mathbb{C}_{\text{ab}}(m-1, y_{1:m-1})- \mathbb{C}_{\text{c}}(m-1, y_{1:m-1})$,
implying $a_{m-1}^*$ and $a_{\lambda,m-1}^*$ are asymptotically the same as $\lambda\to\infty$.
Consequently, $\lim_{\lambda\to\infty}\mathbb C_{\lambda}(m-1,y_{1:m-1})=\mathbb C(m-1,y_{1:m-1})$.
Whence, if the induction hypothesis holds for $n=m,\ldots,N-1$, we have shown that it also holds for $n=m-1$.
Therefore, the induction hypothesis is true for all $n\in[N-1]$.

Moreover, given $Y_{n}=Y_{n}^{(\lambda)}=0$, the corresponding actions are also the same by the last set of constraints in Problems~\eqref{eq:stopping} and \eqref{eq:stopping_CTMC}, i.e., $ a_{\lambda,n}^*(y_{1:n})=a_{n}^*(y_{1:n})=1$ when $y_{n}=0$ for all $n\in[N-1]$.
Combining all the above arguments shows the almost sure convergence of $a_{\lambda,n}^*(Y_{1:n})$ to $a_{n}^*(Y_{1:n})$.

Since for any realization of the condition-monitoring signals $Y_{1:N}$, the actions corresponding to $(a_{n}^*)_{n\in[N-1]}$ and $(a_{\lambda,n}^*)_{n\in[N-1]}$ are the same as $\lambda\to\infty$, it means that $\mathfrak{T}^{(\lambda)}\to\mathfrak{T}^*$ $\mathbb{P}$-almost surely as $\lambda\to \infty$.
Because $\mathfrak{T}^{(\lambda)}$ and $\mathfrak{T}^*$ have a finite support $\{n\delta:n\in[N]\}$, we have $X(\mathfrak{T}^{(\lambda)}+w_{\mathfrak{T}^{(\lambda)}/\delta})\to X(\mathfrak{T}^*+w_{\mathfrak{T}^*/\delta})$ $\mathbb{P}$-almost surely as $\lambda\to\infty$.
Since $X(\mathfrak{T}^{(\lambda)}+w_{\mathfrak{T}^{(\lambda)}/\delta})$ and $X(\mathfrak{T}^*+w_{\mathfrak{T}^*/\delta})$ also only have a finite support $\mathbb{S}=\{1,2,3\}$, 
we have the corresponding cost $C(X(\mathfrak{T}^{(\lambda)}+w_{\mathfrak{T}^{(\lambda)}/\delta}),\mathfrak{T}^{(\lambda)})\rightarrow C(X(\mathfrak{T}^*+w_{\mathfrak{T}^*/\delta}),\mathfrak{T}^*)$ $\mathbb{P}$-almost surely as $\lambda\to\infty$.
Then, because the cost during a mission is bounded from above by $C_\text{m}+C_\text{s}$, it follows from the dominated convergence theorem to conclude $\lim_{\lambda\to\infty}\mathbb E[C(X(\mathfrak{T}^{(\lambda)}+w_{\mathfrak{T}^{(\lambda)}/\delta}),\mathfrak{T}^{(\lambda)})]=\mathbb{E}[ C(X(\mathfrak{T}^*+w_{\mathfrak{T}^*/\delta}),\mathfrak{T}^*)]$.
\Halmos

\subsection{Proof of Lemma~\ref{lemma:increase_rate}}\label{proof:increase_rate}
We prove the monotonicity of $h^{(\lambda)}(t)$ by first showing that $q_{i,m_1+m_2+1}$ is nondecreasing in $i\in\mathcal{X}_2$.
With a slight abuse of notation, let $\{X^{(\lambda)}(t):t>0\}$ be the same CTMC defined in Section~\ref{sec:pomdp} but starting from state $m_1+1$ with probability one at time $t=0$.
Then the absorption time of $X^{(\lambda)}$ is $T_{23}^{(\lambda)}$,
% When the system is at state $i=1,\ldots,m_1$, the system failure rate is constant $\zeta$.
so that the system failure rate at time $t=0$ is 
\begin{align*}
	h^{(\lambda)}(0)&\triangleq \lim_{\Delta t\downarrow 0} \frac{P^{(\lambda)}(X^{(\lambda)}(\Delta t)=m_1+m_2+1\mid X^{(\lambda)}(0)= m_1+1)}{\Delta t}\\
	&=q_{m_1+1,m_1+m_2+1}
	=\lambda p_1(\lambda)
	=\lambda F(1/\lambda),
\end{align*}
where the second equality follows from the definition of transition rates of CTMCs.
By definition of $F(\cdot)$, we have $\lim_{\lambda\to\infty}\lambda F(1/\lambda)=h(0)$, which is larger than $\zeta$ by Assumption~\ref{assump:nondecreasing}.
When the system is at state $i=m_1+1,\ldots,m_1+m_2-1$, the transition rate to the failure state $m_1+m_2+1$ is $q_{i,m_1+m_2+1}=\lambda p_{i-m_1}(\lambda)$.
By \citetappendix[Lemma~1]{khaleghei2021optimal}, there exists $\bar\lambda$ such that $p_{i-m_1}(\lambda)$ is nondecreasing in $i$ for all $\lambda>\bar\lambda$, so is $q_{i,m_1+m_2+1}$.
Moreover, we have $q_{m_1+m_2,m_1+m_2+1}=\lambda\geq q_{i,m_1+m_2+1}$ for all $i\in\mathcal{X}_2$ because $p_{j}(\lambda)\leq 1$ for all $j\in[m_2-1]$ by definition.

Let {${\tilde{\bm p}}(t)\triangleq \left(\frac{p_{m_1+1,j}(t)}{1-p_{m_1+1,m_1+m_2+1}(t)}\right)_{j\in\mathcal{X}_2}$ be the vector of probabilities that $X^{(\lambda)}(t)=j\in\mathcal{X}_2$, given that $X^{(\lambda)}(0)=m_1+1$ and $X^{(\lambda)}(t)<m_1+m_2+1$.}
We then show ${\tilde{\bm p}}(t_1)\succeq_\text{LR}{\tilde{\bm p}}(t_2)$ for any $t_1\geq t_2>0$.
%where ``$\succeq$'' denotes the usual stochastic order, i.e., $\sum_{j=k}^{m_2}\tilde p_{j}(t_1)\geq \sum_{j=k}^{m_2}\tilde p_{j}(t_2)$ for all $k\in[m_2]$, where $\tilde p_{k}(t)$ is the $k$th entry of $\tilde{\bm p}(t)$.
Indeed, 
since our CTMC has a finite number of states, the transition between states follows Kolmogorov's forward equations, i.e., for any $t>0$,
\[
p_{ii}'(t)=-\lambda p_{ii}(t) 
\quad \text{and} \quad
p_{ij}'(t)=\lambda [1-p_{j-m_1-1}(\lambda)] p_{i,j-1}(t)-\lambda p_{ij}(t),
\]
with $i,j\in\mathcal{X}_2$ and $i<j$.
Solving these equations yields
\[
p_{ii}(t)=\exp(-\lambda t)
\quad \text{and} \quad
p_{ij}(t)=\lambda [1-p_{j-m_1-1}(\lambda)]e^{-\lambda t}\int_{0}^t e^{\lambda s}p_{i,j-1}(s) \td s.
\]
Fixing $i=m_1+1$ and recursively using the above identities for $j=m_1+2,\ldots,m_1+m_2$ give
\[
p_{m_1+1,j}(t)=\frac{\prod_{k=m_1+1}^{j-1}[1-p_{k-m_1}(\lambda)]}{(j-m_1-1)!}(\lambda t)^{j-m_1-1}\exp(-\lambda t), \quad t>0,
\] 
for all $j\in\mathcal{X}_2$.
We then have $p_{m_1+1,j+1}(t)/p_{m_1+1,j}(t)=[1-p_{j-m_1}(\lambda)]\lambda t/(j-m_1)$, which increases in $t>0$.
By definition of the likelihood ratio order, ${\tilde{\bm p}}(t_1)\succeq_\text{LR}{\tilde{\bm p}}(t_2)$ for any $t_1\geq t_2>0$ is equivalent to
\[
\frac{p_{m_1+1,j+1}(t_1)}{1-p_{m_1+1,m_1+m_2+1}(t_1)}\cdot
\frac{1-p_{m_1+1,m_1+m_2+1}(t_2)}{p_{m_1+1,j+1}(t_2)}\geq \frac{p_{m_1+1,j}(t_1)}{1-p_{m_1+1,m_1+m_2+1}(t_1)}\cdot
\frac{1-p_{m_1+1,m_1+m_2+1}(t_2)}{p_{m_1+1,j}(t_2)}
\]
\[
\Leftrightarrow \quad 
\frac{p_{m_1+1,j+1}(t_1)}{p_{m_1+1,j}(t_1)}\geq \frac{p_{m_1+1,j+1}(t_2)}{p_{m_1+1,j}(t_2)}
\]
for all $j\in\{m_1+1,\ldots,m_1+m_2-1\}$, which is true based on the above monotone property of $t\mapsto p_{m_1+1,j+1}(t)/p_{m_1+1,j}(t)$.

%For any $i,j\in\mathcal{X}_2$, $a$
%{\color{red}[i don't see why here. could you elaborate a bit more here?] Hence, we have $(p_{m_1+1,j}(t_1))_{j=m_1+1}^{m_1+m_2+1}
%\succeq_\textnormal{LR}(p_{m_1+1,j}(t_2))_{j=m_1+1}^{m_1+m_2+1}$ for any $t_1\geq t_2\geq 0$ by definition.}

The failure rate at time $t>0$ is then given by
\begin{align*}
	h^{(\lambda)}(t)
	&=\lim_{\Delta t\downarrow 0}\frac{P^{(\lambda)}(X^{(\lambda)}(t+\Delta t)=m_1+m_2+1 |X^{(\lambda)}(0)= m_1+1,~X^{(\lambda)}(t)\neq m_1+m_2+1)}{\Delta t} \\
	&=\lim_{\Delta t\downarrow 0}\sum_{j=m_1+1}^{m_1+m_2}
	\frac{P^{(\lambda)}(X^{(\lambda)}(t+\Delta t)=m_1+m_2+1 | X^{(\lambda)}(t)=j)P^{(\lambda)}(X^{(\lambda)}(t)=j| X^{(\lambda)}(0)= m_1+1)}{\Delta t P^{(\lambda)}(X^{(\lambda)}(t)\neq m_1+m_2+1| X^{(\lambda)}(0)= m_1+1)} \\
	&=\frac{\sum_{j=m_1+1}^{m_1+m_2}p_{m_1+1,j}(t)q_{j,m_1+m_2+1}}{1-p_{m_1+1,m_1+m_2+1}(t)}
	=\bm q'\tilde{\bm p}(t),
\end{align*}
where $\bm q\triangleq (q_{j,m_1+m_2+1})_{j=m_1+1}^{m_1+m_2+1}$.
{%\color{red}
%[i don't see why here. could you elaborate a bit more here?] 
Since $\bm q$ is a nondecreasing sequence and ${\tilde{\bm p}}(t_1)\succeq_\textnormal{LR}\tilde{\bm p}(t_2)$ for any $t_1\geq t_2\geq 0$, it then follows from \citetappendix[(1.A.7)]{shaked2007stochastic} that their inner product satisfies $h^{(\lambda)}(t_1)=\bm q'\tilde{\bm p}(t_1)\geq \bm q'\tilde{\bm p}(t_2)= h^{(\lambda)}(t_2)$ for any $0\leq t_2\leq t_1\leq T$.}
\Halmos

\noindent\textbf{Remark.}
Based on the above proof, it suffices to ensure $q_{i,m_1+m_2+1}$ is nondecreasing in $i\in[m_1+m_2]$ to guarantee the monotonicity of ${h}^{(\lambda)}(t)$.
To select a finite value of $\lambda$ for optimization under Assumption~\ref{assump:nondecreasing}, we can compute $q_{i,m_1+m_2+1}$ for all $i\in\{m_1+1,\dots,m_1+m_2+1\}$ by Equations~\eqref{eq:trans_prob} and \eqref{eq:jump_prob}, and check if $q_{i,m_1+m_2+1}$ is nondecreasing in $i$.

\subsection{Proof of Lemma \ref{lemma:piecewise_and_concave}}\label{proof:piecewise and concave}
We prove the lemma by induction.
By Equation~\eqref{3}, $\kappa(w_n,\bm\pi)$ is a linear function of $\bm\pi$ for all $n=0,\ldots,N$.
When $n=N$, the value function $V(N,\bm\pi)$ is also a linear function in $\bm\pi$ by \eqref{eq:bellman} and hence is concave in $\bm\pi$.

We then assume $V(n,\bm\pi)=\min\{V_{\text{ab}}(n,\bm\pi),V_{\text{c}}(n,\bm\pi)\}$ is concave in $\bm\pi$ for all $n=l+1,l+2,\ldots,N$.
We show that both $V_{\text{ab}}(n,\pi)$ and $V_{\text{c}}(n,\pi)$ are concave in $\pi$ when $n=l$.
We first consider $V_{\text{ab}}(n,\bm\pi)$.
By linearity of $\kappa(w_n,\bm\pi)$ in $\bm\pi$ and Equation~\eqref{eq:v_ab}, $V_{\text{ab}}(n,\bm\pi)$ is linear and hence concave in $\bm\pi$ for all $n=0,\ldots,N-1$.

We next consider $V_{\text{c}}(n,\bm\pi)$.
By Equation~\eqref{eq:v_c} and linearity of $\kappa(w_n,\bm\pi)$ in $\bm\pi$, it suffices to show that $U(n,\bm\pi)=\mathbb{E}^{(\lambda)}[V(n+1,\bm\Pi_{n+1})\mid \bm\Pi_n=\bm\pi]$ is concave in $\bm\pi$ when $n=l$.
Let $\bm\pi=\alpha\bm\pi^{(1)}+(1-\alpha)\bm\pi^{(2)}$, for some $\bm\pi^{(1)}, \bm\pi^{(2)}\in\mathcal{S}$ and $\alpha\in(0,1)$.
Then by \eqref{eq:bayes_update}, we can write 
$$
\bm\pi_{n}(\bm\pi,k)
=\dfrac{\bm \pi'\tilde{\mathbf{P}}(\delta)\mathrm{diag}(\tilde{\bm d}_k)}
{\bm \pi'\tilde{\mathbf{P}}(\delta)\tilde{\bm d}_k}
\triangleq\dfrac{\bm \pi ' \bm A(k)}{\bm\pi'\bm b(k)},
\quad \forall n=0,\ldots,N-1,
$$
where $\bm A(k)\in \mathbb{R}^{(m_1+m_2)\times(m_1+m_2)}$ with the ($i$-$j$)th entry being $\bm A(k)_{ij}\triangleq p_{ij}(\delta)\tilde d_{jk}$, and $\bm b(k)\in \mathbb{R}^{m_1+m_2}$ with the $i$th component being $\bm b(k)_i\triangleq \sum_{j=1}^{m_1+m_2}p_{ij}(\delta)\tilde d_{jk}$.
Let $\tilde\alpha \triangleq \alpha\frac{(\bm\pi^{(1)})'\bm b(k)}{\bm \pi'\bm b(k)}$.
Routine algebra then shows that 
$$
\bm\pi_{n}(\bm\pi,k)=\tilde\alpha\dfrac{(\bm \pi^{(1)})'\bm A(k)}{(\bm\pi^{(1)})'\bm b(k)}+(1-\tilde{\alpha})\dfrac{(\bm \pi^{(2)})'\bm A(k)}{(\bm\pi^{(2)})'\bm b(k)}, \quad \forall n=1,\ldots,N-1.
$$
Note that the map from $\alpha$ to $\tilde\alpha$ is one-to-one from $[0,1]$ to $[0,1]$.
Therefore, we see that
\begin{align*}
	V(l+1,\bm\pi_{l+1}(\bm\pi,k)) &= V\left(l+1,\dfrac{\bm \pi' \bm A(k)}{\bm\pi'\bm b(k)}\right)
	= V\left(l+1,\tilde\alpha\dfrac{(\bm\pi^{(1)})'\bm A(k)}{(\bm\pi^{(1)})'\bm b(k)}+(1-\tilde{\alpha})\dfrac{(\bm\pi^{(2)})'\bm A(k)}{(\bm\pi^{(2)})'\bm b(k)}\right) \\
	& \ge \tilde\alpha V\left(l+1,\dfrac{ (\bm\pi^{(1)})'\bm A(k)}{(\bm\pi^{(1)})'\bm b(k)}\right)+
	(1-\tilde\alpha)V\left(l+1,\dfrac{(\bm{\pi}^{(2)})'\bm A(k)}{(\bm\pi^{(2)})'\bm b(k)}\right),
\end{align*}
where the inequality follows from the induction hypothesis that $V(l+1,\bm\pi)$ is concave in $\bm\pi$.

Plugging the above inequality in \eqref{eq:U_function} gives
\begin{align*}
	&~U(l,\pi)=\sum_{k=1}^K\bm\pi'\bm b(k)V\left(l+1,\dfrac{\bm \pi' \bm A(k)}{\bm\pi'\bm b(k)}\right)\\
	\ge &~\sum_{k=1}^K \bm\pi'\bm b(k)\left[
	\tilde\alpha V\left(l+1,\dfrac{ (\bm\pi^{(1)})'\bm A(k)}{(\bm\pi^{(1)})'\bm b(k)}\right)
	+(1-\tilde\alpha)V\left(l+1,\dfrac{(\bm \pi^{(2)})'\bm A(k)}{(\bm\pi^{(2)})'\bm b(k)}\right)\right] \\
	= &~\sum_{k=1}^K\left[ \alpha (\bm\pi^{(1)}) '\bm b(k)V\left(l+1,\dfrac{(\bm \pi^{(1)})'\bm A(k)}{(\bm\pi^{(1)})' \bm b(k)}\right)
	+(1-\alpha)(\bm\pi^{(2)})'\bm b(k)V\left(l+1,\dfrac{ (\bm\pi^{(2)})'\bm A(k)}{(\bm\pi^{(2)})'\bm b(k)}\right)\right] \\
	= &~\alpha U(l,\bm\pi^{(1)}) + (1-\alpha) U(l,\bm\pi^{(2)}),
\end{align*}
which shows $U(l,\bm\pi)$ is concave in $\bm\pi$.

Since the minimum operator preserves concavity, $V(n,\bm\pi)$ is also concave in $\bm\pi$ when $n=l$.
This completes the induction.

We then show that $V(n,\bm\pi)$ is piecewise linear for all $n=0,\ldots,N$.Assume that the value function is piecewise linear in $\bm\pi$ for all $n=l+1,\ldots,N$, so $V(n,\bm\pi)=\max_{\bm\alpha\in \mathcal{A}}\bm\alpha'\bm\pi$ for some finite set $\mathcal{A}$, where $\bm\alpha=(\alpha_j)_{j=1}^{m_1+m_2}$.
Then for $n=l$, we already show that $V_\textnormal{ab}(l,\bm\pi)$ is linear in $\bm\pi$.
Moreover, routine algebra shows that
$$
V_\textnormal{c}(l,\bm\pi)=\kappa(w_l,\bm\pi)+ \sum_{k=1}^K \max_{\bm\alpha\in\mathcal{A}}\sum_{i=1}^{m_1+m_2}\sum_{j=1}^{m_1+m_2} \pi_i p_{ij}(\delta)d_{jk}\alpha_j.
$$
Since $\kappa(w_n,\bm\pi)$ is linear in $\bm\pi$,
the above display is clearly piecewise linear in $\bm\pi$ as well.
Since $V(n,\bm\pi)$ is a pointwise minimum of a linear function and a piecewise linear function, we have $V(n,\bm\pi)$ is also piecewise linear in $\bm\pi$ for any fixed $n$.
\Halmos

%\subsection{Proof of Lemma  \ref{lemma:upper_bound}}\label{proof:upper bound}
%If the system is functional at time $n\delta$, we consider a policy that does not conduct any abort from time $n\delta$ to the end of the mission.
%Routine algebra shows that the cost-to-go from time $n\delta$ is exactly $\overline{V}_{\text{c}}(n,\bm\pi)$.
%On the other hand, the cost-to-go $V_{\text{c}}(n,\bm\pi)$ can be understood as the cost that no abort is conducted at time $n\delta$, and then we follow the optimal action during the remaining of the mission period.
%It follows from the above relationship that the cost-to-go from time $(n+1)\delta$ of $V_{\text{c}}(n,\bm\pi)$ is no greater than that of $\overline{V}_{\text{c}}(n,\bm\pi)$.
%\Halmos

\subsection{Proof of Lemma~\ref{lemma:monotonous_pi_n}}\label{proof:monotonousin pi and n}
{
We first show the monotonicity of $p_{i,m_1+m_2+1}(t)$ in $i$ by induction.
First, $p_{m_1+m_2+1,m_1+m_2+1}(t)=1$, and by solving the Kolmogorov backward equation, $p_{m_1+m_2,m_1+m_2+1}(t)=1-\exp(-\lambda t)<1$.
For induction, assume $p_{i,m_1+m_2+1}(t)\leq p_{i+1,m_1+m_2+1}(t)\leq\cdots\leq p_{m_1+m_2+1,m_1+m_2+1}(t)$ for some $i=2,\ldots,m_1+m_2$.
We consider three cases to show $p_{i-1,m_1+m_2+1}(t)\leq p_{i,m_1+m_2+1}(t)$.

\textbf{Case 1: $i-1>m_1$.}
By the Kolmogorov backward equation and noting $p_{m_1+m_2+1,m_1+m_2+1}(t)=1$,
\[
p_{i-1,m_1+m_2+1}'(t)=\lambda(1-p_{i-m_1-1}(\lambda))p_{i,m_1+m_2+1}(t)+\lambda p_{i-m_1-1}(\lambda)-\lambda p_{i-1,m_1+m_2+1}(t),
\]
yielding
\[
p_{i-1,m_1+m_2+1}(t)=\lambda\int_{0}^{t}\exp(\lambda(u-t))[(1-p_{i-m_1-1}(\lambda))p_{i,m_1+m_2+1}(u)+p_{i-m_1-1}(\lambda)]\mathrm{d}u.
\]
Similarly, we can solve
\[
p_{i,m_1+m_2+1}(t)=\lambda\int_{0}^{t}\exp(\lambda(u-t))[(1-p_{i-m_1}(\lambda))p_{i+1,m_1+m_2+1}(u)+p_{i-m_1}(\lambda)]\mathrm{d}u.
\]
Then it suffices to show $(1-p_{i-m_1}(\lambda))p_{i+1,m_1+m_2+1}(u)+p_{i-m_1}(\lambda)\geq (1-p_{i-m_1-1}(\lambda))p_{i,m_1+m_2+1}(u)+p_{i-m_1-1}(\lambda)$, which holds because
\begin{align*}
	& (1-p_{i-m_1}(\lambda))p_{i+1,m_1+m_2+1}(u)+p_{i-m_1}(\lambda)
	\geq (1-p_{i-m_1}(\lambda))p_{i,m_1+m_2+1}(u)+p_{i-m_1}(\lambda) \\
	=~ & p_{i,m_1+m_2+1}(u)+p_{i-m_1}(\lambda)(1-p_{i,m_1+m_2+1}(u)) 
	\geq p_{i,m_1+m_2+1}(u)+p_{i-m_1-1}(\lambda)(1-p_{i,m_1+m_2+1}(u)).
\end{align*}
Here, the first inequality follows from the induction hypothesis, and the second inequality follows from Assumption~\ref{assump:nondecreasing}.

\textbf{Case 2: $i-1=m_1$.}
Similarly, we can solve the Kolmogorov backward equations for $p_{m_1,m_1+m_2+1}(t)$ and $p_{m_1+1,m_1+m_2+1}(t)$, where
\[
p_{m_1,m_1+m_2+1}(t)=\int_{0}^{t}\exp(\lambda(u-t))[\zeta+(\lambda-\zeta)p_{m_1+1,m_1+m_2+1}(u)]\mathrm{d}u,
\]
and $p_{m_1+1,m_1+m_2+1}(t)$ is given in Case~1 by letting $i=m_1+1$.
Then it suffices to show $\zeta+(\lambda-\zeta)p_{m_1+1,m_1+m_2+1}(u)\leq \lambda[p_1(\lambda)(1-p_{m_1+2,m_1+m_2+1}(u))+p_{m_1+2,m_1+m_2+1}(u)]$, which is true as 
\begin{align*}
	\zeta+(\lambda-\zeta)p_{m_1+1,m_1+m_2+1}(u)
	\leq~ &  \zeta+(\lambda-\zeta)p_{m_1+2,m_1+m_2+1}(u) \\
	=~ & \zeta(1-p_{m_1+2,m_1+m_2+1}(u))+\lambda p_{m_1+2,m_1+m_2+1}(u) \\
	\leq~ & \lambda[p_1(\lambda)(1-p_{m_1+2,m_1+m_2+1}(u))+p_{m_1+2,m_1+m_2+1}(u)].
\end{align*}
Here, the first inequality follows from the induction hypothesis.
The second inequality follows from Assumption~\ref{assump:nondecreasing} and Lemma~\ref{lemma:increase_rate} where for $\lambda> \bar{\lambda}$, we have $\zeta\leq\lambda F(1/\lambda)=\lambda p_1(\lambda)$.

\textbf{Case 3: $i-1<m_1$.}
We can again solve the Kolmogorov backward equation to obtain
\[
p_{i-1,m_1+m_2+1}(t)=\int_{0}^t \exp(\lambda(u-t))[\zeta+(\lambda-\zeta)p_{i,m_1+m_2+1}(u)]\mathrm{d}u;
\]
\[
p_{i,m_1+m_2+1}(t)=\int_{0}^t \exp(\lambda(u-t))[\zeta+(\lambda-\zeta)p_{i+1,m_1+m_2+1}(u)]\mathrm{d}u.
\]
It is then easy to see $p_{i-1,m_1+m_2+1}(t)\leq p_{i,m_1+m_2+1}(t)$ due to the induction hypothesis.}

We then show the monotonicity of $V_\text{ab}(n,\bm\pi)$ and $V_\text{c}(n,\bm\pi)$.
The value function for abort is $V_\textnormal{ab}(n,\bm\pi)=C_\textnormal{m}+C_\textnormal{s}\kappa(w_n,\bm\pi)$.
We show that $\kappa(w_n,\bm\pi)$ is monotone likelihood ratio (MLR) increasing in $\bm\pi$ for all $n=1,\ldots,N-1$.
Indeed, recall $\kappa(w_n,\bm\pi)=\sum_{i=1}^{m_1+m_2}\pi_i p_{i,m_1+m_2+1}(w_n)$.
By Lemma~\ref{lemma:monotonous_pi_n}, $p_{i,m_1+m_2+1}(\delta)$ is nondecreasing in $i$.
The MLR monotonicity of $\kappa(w_n,\bm\pi)$ in $\bm\pi$ then follows from the fact that $\sum_{j}\pi_{j}^{(1)}f(j)\geq\sum_{j}\pi_{j}^{(2)}f(j)$ for any nondecreasing function $f$ and $\bm\pi^{(1)},\bm\pi^{(2)}\in\mathcal{S}$ such that $\bm\pi^{(1)}\succeq_\textnormal{LR}\bm\pi^{(2)}$ \citepappendix{shaked2007stochastic}.
To show the monotone property of $V_{\text{ab}}(n,\bm\pi)$ in $n$, we note that $w_n$ is nondecreasing in $n$ by Assumption~\ref{assump:rescue}.
Hence, $\kappa(w_n,\bm\pi)=\sum_{i=1}^{m_1+m_2}\pi_ip_{i,m_1+m_2+1}(w_n)$ is nondecreasing with $n$ for any fixed $\bm\pi\in\mathcal{S}$,
whence $V_{\text{ab}}(n,\bm\pi)$ is also nondecreasing with $n$. 

For the function $\overline{V}_{\text{c}}(n,\bm\pi)$, we have 
\begin{align*}
	\overline{V}_{\text{c}}(n,\bm\pi)
	&=(C_{\text{m}}+C_{\text{s}})\kappa((N-n)\delta+w_N,\bm\pi) \\
	&=(C_{\text{m}}+C_{\text{s}})\sum_{i=1}^{m_1+m_2}\pi_ip_{i,m_1+m_2+1}((N-n)\delta+w_N).
\end{align*}
Note that the term $(N-n)\delta+w_N$ is nonincreasing with $n$. 
By the same way, we can prove $\overline{V}_{\text{c}}(n,\bm\pi)$ is MLR increasing with $\bm\pi$ and nonincreasing with $n$.
\Halmos

\noindent\textbf{Remark.}
{We remark that the construction of $G^{(\lambda)}$ using $\lambda-\zeta$ in \eqref{eq:herlang2} leads to identical transition rates out of any transient state, facilitating the analysis in Case 2 in the above proof.}

\subsection{Proof of Theorem~\ref{thm:never_abort}}\label{proof:never abort}
Based on Lemma \ref{lemma:monotonous_pi_n}, if 
\begin{equation}\label{eq:sufficient_cond}
	V_\textup{ab}(n,\bm\pi)\geq\overline{V}_\textup{c}(n,\bm\pi)\geq {V}_\textup{c}(n,\bm\pi), \quad \forall n=1,\ldots,N-1,~\bm\pi\in\mathcal{S},   
\end{equation}
then the mission is never aborted during the mission time $[0,H]$.
Note that $V_{\text{ab}}(0,\bm e_1)<V_{\text{ab}}(n,\bm\pi)$ for all $n>1$ and $\bm\pi\in\mathcal{S}$, where the inequality follows from Lemma~\ref{lemma:monotonous_pi_n}.
Similarly, we have $\overline V_{\text{c}}(0,\bm e_{m_1+m_2})>\overline V_{\text{c}}(n,\bm\pi)$ for all $n>1$ and $\bm\pi \in \mathcal{S}$ by monotonicity of $\overline V_{\text{c}}(\cdot)$.
Then a sufficient condition for \eqref{eq:sufficient_cond} to hold is $V_{\text{ab}}(0,\bm e_{1})\ge \overline V_{\text{c}}(0,\bm e_{m_1+m_2})$, i.e., $C_{\text{m}}\ge(C_{\text{s}}+C_{\text{m}})(1-\exp[-\lambda(T + w_N)])$, where we use the property $w_0=0$.
Routine algebra yields Eq.~\eqref{eq:cost_ratio}.
\Halmos

\subsection{Proof of Theorem~\ref{thm:threshold_n}}\label{proof:threshold n}
By Lemma~\ref{lemma:monotonous_pi_n}, $V_\text{ab}(n,\bm\pi)$ and $\overline V_\text{c}(n,\bm\pi)$ are, respectively, nondecreasing and nonincreasing in $n$ for any given $\bm\pi$.
Hence, for any given $\bm\pi\in\mathcal{S}$, there may exist a threshold $\hat n(\bm\pi)$ dependent on $\bm\pi$ such that $V_\text{ab}(n,\bm\pi)\geq\overline{V}_\text{c}(n,\bm\pi)$ for all $n\geq \hat n(\bm\pi)$ and $V_\text{ab}(n,\bm\pi)<\overline{V}_\text{c}(n,\bm\pi)$ for all $n< \hat n(\bm\pi)$.
Let $\hat n\triangleq \max_{\bm\pi\in\mathcal{S}}\hat{n}(\bm\pi)$.
Then we see $V_\text{ab}(n,\bm\pi)\geq\overline{V}_\text{c}(n,\bm\pi)\geq V_\text{c}(n,\bm\pi)$ for all $n\geq\hat n$ and $\bm\pi\in\mathcal{S}$.
This tallies with the theorem statement that the optimal action is to continue the mission when $n\geq\hat n$, no matter of the observed signals.
As such, it is easy to see $V(n,\bm\pi)=\overline{V}_\textnormal{c}(n,\bm\pi)$ for all $\bm\pi\in\mathcal{S}$, $n=\hat n,\ldots,N-1$, by the definition of $\overline{V}_\textnormal{c}(n,\bm\pi)$.
\Halmos

\subsection{Proof of Lemma ~\ref{lemma:MLR_montonicity}}\label{proof:MLR montonicity}
We first show the following lemma for the transition probability matrix $\tilde{\mathbf{P}}(\delta)$ defined in Section~\ref{subsec:pomdp_model}.

\begin{lemma}\label{lemma:TP2}
	{
	%[since this lemma is used to derive the next lemma only, i suggest that you move this lemma to the supp to save space.]
	The matrix $\tilde{\mathbf{P}}(t)$ is totally positive of order $2$ (TP2) for all $t>0$.}
\end{lemma}
\textit{Proof of Lemma~\ref{lemma:TP2}}:
{
% We follow the procedures in \citeappendix{khaleghei2021optimal} for proof.
Let $\Delta_{ij}(t)\triangleq p_{ik}(t)p_{jl}(t)-p_{il}(t)p_{jk}(t)$ for $i,j\in[m_1+m_2]$ and some fixed $k,l\in[m_1+m_2]$.
Note $\Delta_{ii}(t)=0$ for all $i\in[m_1+m_2]$ and $t>0$.
Without loss of generality, it suffices to show $\Delta_{ij}(t)\geq 0$ for all $t>0$ when $i<j$ and $k<l$.
If $k<j$, then $\Delta_{ij}(t)=p_{ik}(t)p_{jl}(t)\geq 0$.
Thus it suffices to prove $\Delta_{ij}(t)\geq 0$ for $i<j\leq k<l$.
We show this by following an exhaustive list of cases.
	
\textbf{Case 1: $i<j=k<l\leq m_1$.}
We can use the Kolmogorov backward equation to show
$\Delta_{ij}'(t)=-2\lambda\Delta_{ij}(t)+(\lambda-\zeta)[\Delta_{i+1,j}(t)+p_{j+1,l}(t)p_{ik}(t)]$, 	which yields
$\Delta_{ij}(t)=(\lambda-\zeta)\int_{0}^t \exp(2\lambda(u-t))(\Delta_{i+1,j}(u)+p_{j+1,l}(u)p_{ik}(u))\td u$.	
Due to the recursive formula, it suffices to show $\Delta_{ij}(t)\geq 0$ when $i=k-1$ and $j=k$ for any fixed $k$ and $l$.
This holds as $\Delta_{k-1,k}(t)=(\lambda-\zeta)\int_{0}^{t}\exp(2\lambda(u-t))(\Delta_{kk}(u)+p_{k+1,l}(u)p_{k-1,k}(u))\td u\geq 0$.
	
\textbf{Case 2: $i<j<k<l\leq m_1$.}
We can use the Kolmogorov backward equation to show
$\Delta_{ij}'(t)=-2\lambda\Delta_{ij}(t)+(\lambda-\zeta)[\Delta_{i+1,j}(t)+\Delta_{i,j+1}(t)]$, 
which yields
$\Delta_{ij}(t)=(\lambda-\zeta)\int_{0}^t \exp(2\lambda(u-t))(\Delta_{i+1,j}(u)+\Delta_{i,j+1}(u))\td u$.
Based on the recursive formula and $\Delta_{ii}=0$, it suffices to show $\Delta_{i,k-1}(t)\geq 0$ for all $i<k-1$.
As $\Delta_{k-1,k-1}(t)=0$ and we show $\Delta_{k-2,k}(t)\geq 0$ in Case~1, we can use the recursive formula to show $\Delta_{k-2,k-1}(t)\geq 0$.
We then use the recursive formula again to show $\Delta_{i,k-1}(t)\geq 0$ for $i<k-2$.
	
\textbf{Case 3: $i<j\leq k\leq m_1<l$.}
When $j=k$, the proof follows the same procedure as in Case~1.
When $j<k$, the proof follows the same procedure as in Case~2.

\textbf{Case 4: $m_1< i<j\leq k<l$.}
We can first use nearly the same proof as in Cases~1 for $j=k$ and then follow nearly the same proof for Case~2 to show $\Delta_{ij}(t)\geq 0$ for $j<k$.

\textbf{Case 5: $i\leq m_1<j\leq k<l$.}
We can solve the Kolmogorov backward equation to obtain $\Delta_{ij}(t)=\int_{0}^t \exp(2\lambda(u-t))[(\lambda-\zeta)\Delta_{i+1,j}(u)+\lambda(1-p_{j-m_1}(\lambda))p_{ik}(u)p_{j+1,l}(u)]\td u$ when $j=k$ and $\Delta_{ij}(t)=\int_{0}^t \exp(2\lambda(u-t))[(\lambda-\zeta)\Delta_{i+1,j}(u)+\lambda(1-p_{j-m_1}(\lambda))\Delta_{i,j+1}(u)]\td u$ when $j<k$.
Using the recursive formula, it suffices to show $\Delta_{m_1,k}(t)\geq 0$.
We can show this by using the recursive formula again and noting $\Delta_{m_1+1,k}(t)\geq 0$ as shown in Case~4.

\textbf{Case 6: $i<j\leq m_1<k<l$.}
Similar to Case~2, we have $\Delta_{ij}(t)=(\lambda-\zeta)\int_{0}^t \exp(2\lambda(u-t))(\Delta_{i+1,j}(u)+\Delta_{i,j+1}(u))\td u$.
Based on this recursion, $\Delta_{ij}(t)\geq 0$ if $\Delta_{i,m_1}(t)\geq 0$ for all $i<m_1$.
We can use the recursive formula to show $\Delta_{m_1-1,m_1}(t)\geq 0$, as $\Delta_{m_1,m_1}(t)=0$ and $\Delta_{m_1-1,m_1+1}(t)\geq 0$ as shown in Case~5.
Again using the recursive formula shows $\Delta_{ij}(t)\geq 0$ for all $i,j$ in this case.

Since the above proof is for all possible values of $k,l$, we conclude that $\tilde{\mathbf{P}}(\delta)$ is TP2.
We remark that the CTMC construction in Figure~\ref{fig:CTMC} significantly simplifies the Kolmogorov backward equations, hence facilitating the proof of the TP2 property for $\tilde{\mathbf{P}}(t)$.
\Halmos}

We now show the first statement of Lemma~\ref{lemma:MLR_montonicity} on the MLR ordering of $\bm\pi_{n}(\bm\pi,k)$ with respect to $\bm\pi$. 
We have $\bm\pi_n(\bm\pi,k)\propto \bm\pi'\tilde{\mathbf{P}}(\delta)\mathrm{diag}(\tilde {\bm d}_k)$ based on \eqref{eq:bayes_update}. 
By Lemma~\ref{lemma:TP2}, $\tilde{{\mathbf P}}(\delta)$ is TP2.
Moreover, $\mathrm{diag}(\tilde {\bm d}_k)$ is TP2 as any diagonal matrix is TP2 by definition.
By \citetappendix[Theorem 1.2]{karlin1980classes}, the product of two TP2 matrices is TP2.
Since the MLR ordering is preserved by right multiplying a TP2 matrix, we see $\bm\pi_n(\bm\pi^{(1)},k)\succeq_\text{LR}\bm\pi_n(\bm\pi^{(2)},k)$ for any $\bm\pi^{(1)}\succeq_\text{LR}\bm\pi^{(2)}$.

We then show the second statement on the MLR ordering of $\bm\pi_{n}(\bm\pi,k)$ with respect to $k$. 
By definition of MLR ordering, 
it suffices to show for any $k_1\geq k_2$,
$$
\dfrac{\bm\pi_{ni}(\bm\pi,k_1)}{\bm\pi_{ni}(\bm\pi,k_2)}\ge\dfrac{\bm\pi_{nj}(\bm\pi,k_1)}{\bm\pi_{nj}(\bm\pi,k_2)}, \quad \forall i>j.
$$
Based on Equation~\eqref{eq:bayes_update}, the above inequality is equivalent to 
$$
\dfrac{\tilde d_{ik_1}}{\tilde d_{ik_2}}\ge \dfrac{\tilde d_{jk_1}}{\tilde d_{jk_2}},\quad \forall i>j.
$$ 
This is satisfied by Assumption~\ref{assump:TP2} that the state-observation matrix $\mathbf{D}$ is TP2.
\Halmos

\noindent\textbf{Remark.} An alternative proof of Lemma~\ref{lemma:TP2} is to use the fact that a tridiagonal matrix $\mathbf{P}$ is TP2 if and only if $\mathbf{P}_{ii}\mathbf{P}_{i+1,i+1}\geq \mathbf{P}_{i,i+1}\mathbf{P}_{i+1,i}$ for all $i$ \citepappendix{krishnamurthy2011bayesian}.
By our construction of the CTMC in Figure~\ref{fig:CTMC}, $\tilde{\mathbf{P}}(t)$ trivially satisfies the above inequality.
This underscores the benefit of constructing the CTMC as in Figure~\ref{fig:CTMC}.

\subsection{Proof of Lemma~\ref{lemma:monotone}}\label{poof:monotone}
We first prove the MLR monotonicity of $V_\textnormal{c}(n,\bm\pi)$ in $\bm\pi\in\mathcal{S}$ by induction.
When $n=N-1$, $V_\textnormal{c}(N-1,\bm\pi)=(C_\textnormal{s}+C_\textnormal{m})\kappa(\delta+w_N,\bm\pi)$, and the MLR monotonicity holds since $\kappa(\delta+w_N,\bm\pi)$ is MLR increasing in $\bm\pi\in\mathcal{S}$.
Next assume that $V_\textnormal{c}(n,\bm\pi)$ is MLR increasing in $\bm\pi\in\mathcal{S}$ for all $n=l+1,\ldots,N-1$, so $V(n,\bm\pi)$ also has this MLR monotone property for $n=l+1,\ldots,N-1$.
Then for $n=l$, we have
$$
V_\textnormal{c}(l,\bm\pi)=(C_\textnormal{m}+C_\textnormal{s})\kappa(\delta,\bm\pi)+\sum_{k=1}^K V(l+1,\bm\pi_{l+1}(\bm\pi,k))\sum_{i=1}^{m_1+m_2}\pi_i\sum_{j=1}^{m_1+m_2} p_{ij}(\delta)\tilde d_{jk},
$$
where $\tilde d_{jk}\triangleq d_{1k}$ for $j\in\mathcal{X}_1$ and $\tilde d_{jk}\triangleq d_{2k}$ for $j\in\mathcal{X}_2$.
Let $\tilde h(k,\bm\pi)\triangleq\sum_{i=1}^{m_1+m_2}\pi_i\sum_{j=1}^{m_1+m_2} p_{ij}(\delta)\tilde d_{jk}$ be the conditional probability that the next observed signal is $k\in[K]$ given the current belief state $\bm\pi$.
Then we show that the likelihood ratio $\tilde h(k,\bm\pi^{(1)})/\tilde h(k,\bm\pi^{(2)})$ is increasing in $k$ for any $\bm\pi^{(1)}\succeq_\textnormal{LR}\bm\pi^{(2)}$.
Indeed, the transition matrix $\tilde{\mathbf{P}}(\delta)$ are TP2 by Lemma~\ref{lemma:TP2}.
Let $\tilde{\bm\pi}^{(i)}\triangleq \tilde{\mathbf{P}}(\delta)\bm\pi^{(i)}$, $i=1,2$.
Since $\bm\pi^{(1)}\succeq_{\textnormal{LR}}\bm\pi^{(2)}$ and $\tilde{\mathbf{P}}(\delta)$ is TP2,
we have $\tilde{\bm\pi}^{(1)}\succeq_{\textnormal{LR}}\tilde{\bm\pi}^{(2)}$ by \citetappendix[Lemma~1.3]{lovejoy1987some}.
For any $k_1>k_2$, we have
\begin{equation*}
	\frac{\sum_{i=1}^{m_1+m_2} \tilde\pi_{i}^{(1)}\tilde d_{ik_1}}
	{\sum_{i=1}^{m_1+m_2} \tilde\pi_{i}^{(2)}\tilde d_{ik_1}}\geq
	\frac{\sum_{i=1}^{m_1+m_2} \tilde\pi_{i}^{(1)}\tilde d_{ik_2}}
	{\sum_{i=1}^{m_1+m_2} \tilde\pi_{i}^{(2)}\tilde d_{ik_2}}
	~ \Leftrightarrow ~
	\sum_{i=1}^{m_1+m_2}\sum_{j=1}^{i-1}(\tilde\pi_{i}^{(1)}\tilde\pi_{j}^{(2)}-\tilde\pi_{j}^{(1)}\tilde\pi_{i}^{(2)})(\tilde d_{ik_1}\tilde d_{jk_2}-\tilde d_{ik_2}\tilde d_{jk_1})\geq 0.
\end{equation*}
The last inequality holds because (i) $\tilde\pi_{i}^{(1)}\tilde\pi_{j}^{(2)}-\tilde\pi_{j}^{(1)}\tilde\pi_{i}^{(2)}\geq 0$ for all $i>j$ due to $\tilde{\bm\pi}^{(1)}\succeq_\textnormal{LR} \tilde{\bm\pi}^{(2)}$ and (ii) $\tilde d_{ik_1}\tilde d_{jk_2}-\tilde d_{ik2}\tilde d_{jk_1}\geq 0$ for all $k_1>k_2$ due to the assumption that $\mathcal{D}$ is TP2 so $\tilde{\bm d}_{k_1}\succeq_\textnormal{LR} \tilde{\bm d}_{k_2}$.

For $\bm\pi^{(1)}\succeq_\textnormal{LR}\bm\pi^{(2)}$, recall that $\bm\pi_{n}(\bm\pi,k)$ is MLR increasing in $\bm\pi$ and $k$ for fixed $n$ by Lemma~\ref{lemma:MLR_montonicity}.
It then follows that
\begin{align*}
	\sum_{k=1}^K V(l+1,\bm\pi_{l+1}(\bm\pi^{(1)},k))\tilde h(k,\bm\pi^{(1)}) 
	& \geq \sum_{k=1}^K V(l+1,\bm\pi_{l+1}(\bm\pi^{(1)},k))\tilde h(k,\bm\pi^{(2)}) \\
	& \geq \sum_{k=1}^K V(l+1,\bm\pi_{l+1}(\bm\pi^{(2)},k))\tilde h(k,\bm\pi^{(2)}),
\end{align*}
where the first inequality follows from $(\tilde h(k,\bm\pi^{(1)}))_{k=1}^K\succeq_\textnormal{LR}(\tilde h(k,\bm\pi^{(2)})_{k=1}^K$, the induction hypothesis, and the MLR monotonicity of $\bm\pi_n(\bm\pi,k)$ in $k$, and the second inequality follows from the induction hypothesis and the MLR monotonicity of $\bm\pi_n(\bm\pi,k)$ in $\bm\pi$.
Together with the MLR monotonicity of $\kappa(\delta,\bm\pi)$ in $\bm\pi\in\mathcal{S}$, we see that $V_\textnormal{c}(l,\bm\pi)$ is also MLR increasing in $\bm\pi\in\mathcal{S}$.

We have shown the MLR monotonicity of $V_\text{ab}(n,\bm\pi)$ with $\bm\pi\in\mathcal{S}$ for all $n=0,\ldots,N-1$ in the proof of Lemma~\ref{lemma:monotonous_pi_n}.
Since $V(n,\bm\pi)=\min\{V_\textnormal{ab}(n,\bm\pi),V_\textnormal{c}(n,\bm\pi)\}$, its MLR monotonicity in $\bm\pi\in\mathcal{S}$ follows directly from the monotone properties of $V_\textnormal{ab}(n,\bm\pi)$ and $V_\textnormal{c}(n,\bm\pi)$.
\Halmos

\subsection{Proof of Lemma~\ref{lemma:ab_opt}}\label{proof:ab_opt}
Assume on the opposite that there exists $n$ such that $V_\textnormal{c}(n-1,\bm e_{m_1+m_2})<V_\textnormal{ab}(n-1,\bm e_{m_1+m_2})$ and $V_\textnormal{c}(n,\bm e_{m_1+m_2})> V_\textnormal{ab}(n,\bm e_{m_1+m_2})$.
We then have 
\begin{align*}
	V_\textnormal{c}(n-1,\bm e_{m_1+m_2})&=C_\textnormal{m}\kappa(\delta,\bm e_{m_1+m_2})+V(n,\bm e_{m_1+m_2})(1-\kappa(\delta,\bm e_{m_1+m_2}))\\
	&=C_\textnormal{m}+C_\textnormal{s}(p_{m_1+m_2,m_1+m_2+1}(\delta)+p_{m_1+m_2,m_1+m_2}(\delta)p_{m_1+m_2,m_1+m_2+1}(w_{n})),
\end{align*}
where the second equality follows from $V(n,\bm e_{m_1+m_2})=V_\textnormal{ab}(n,\bm e_{m_1+m_2})$.
Moreover, we have 
$$
V_\textnormal{ab}(n-1,\bm e_{m_1+m_2})=C_\textnormal{m}+C_\textnormal{s}p_{m_1+m_2,m_1+m_2+1}(w_{n-1}).
$$
Then $V_\textnormal{c}(n-1,\bm e_{m_1+m_2})<V_\textnormal{ab}(n-1,\bm e_{m_1+m_2})$ is equivalent to 
$$
p_{m_1+m_2,m_1+m_2+1}(\delta)+p_{m_1+m_2,m_1+m_2}(\delta)p_{m_1+m_2,m_1+m_2+1}(w_{n})<p_{m_1+m_2,m_1+m_2+1}(w_{n-1}).
$$
We note that left and right hand side of the above inequality are the system failure probabilities in the next $\delta+w_n$ and $w_{n-1}$ time units, respectively.
Since $w_n$ is nondecreasing in $n$ by Assumption~\ref{assump:rescue}, we have $\delta+w_n>w_{n-1}$, so that $p_{m_1+m_2,m_1+m_2+1}(\delta)+p_{m_1+m_2,m_1+m_2}(\delta)p_{m_1+m_2,m_1+m_2+1}(w_{n})>p_{m_1+m_2,m_1+m_2+1}(w_{n-1})$.
This leads to a contradiction.
\Halmos

\subsection{Proof of Proposition~\ref{prop:order}}\label{proof:order}
For any $1\leq j<i<m_1+m_2$, we have $\pi_{1i}/\pi_{1j}=\pi_{2i}/\pi_{2j}$ based on the inverse transform of the Cartesian belief vector in Table \ref{tab:spherical}.
For any $1<j<i=m_1+m_2$ and $l\in\{1,2\}$, we have 
$$
\frac{\pi_{m_1+m_2}^{(l)}}{\pi_{j}^{(l)}}
=\frac
{1+r_l\cos\phi_{m_1+m_2-1}\prod_{k=1}^{m_1+m_2-2}\sin\phi_k}
{r_l\cos\phi_{j-1}\prod_{k=1}^{j-2}\sin\phi_k}
=\frac
{1/r_l+\cos\phi_{m_1+m_2-1}\prod_{k=1}^{m_1+m_2-2}\sin\phi_k}
{\cos\phi_{j-1}\prod_{k=1}^{j-2}\sin\phi_k}.
$$
We then have $\pi_{m_1+m_2}^{(1)}/\pi_{j}^{(1)}>\pi_{m_1+m_2}^{(2)}/\pi_{j}^{(2)}$ as $r_1<r_2$.
We can similarly prove such an inequality for $1=j<i=m_1+m_2$.
Based on all the above cases, we can conclude $\bm\pi^{(1)}\succeq_\textnormal{LR}\bm\pi^{(2)}$.
\Halmos

\subsection{Proof of Theorem~\ref{thm:struct_policy}}\label{proof: structure of the policy}
We first consider the case $n\leq \tilde n$, where $V_\textnormal{ab}(n,\bm e_{m_1+m_2})\leq V_\textnormal{c}(n,\bm e_{m_1+m_2})$.
Recall that $r(\bm\pi)$ and $\bm\Phi(\bm\pi)$ are the radius and spherical angle of $\bm\pi$ in the spherical coordinate system with origin $\bm e_{m_1+m_2}$, respectively.
Define the subsimplex $\mathcal{S}_{-(m_1+m_2)}\triangleq \{\bm\pi\in\mathcal{S}:\bm\pi_{m_1+m_2}=0 \}$.
For any $\bm\pi_0\in\mathcal{S}_{-(m_1+m_2)}$, we consider the line segment $\mathcal{S}(\bm\pi_0)=\{\bm\pi\in\mathcal{S}: \bm\Phi(\bm\pi)=\bm\Phi(\bm\pi_0),~0\leq r(\bm\pi)\leq r(\bm\pi_0) \}$.
We then consider the value function restricted on this line segment $\{V(n,\bm\pi)\}_{\bm\pi\in\mathcal{S}(\bm\pi_0)}$ and write this function in the spherical coordinates as $\tilde V(n,r,\bm\Phi(\bm\pi_0))$.
Similarly, we write $\tilde V_\textnormal{ab}(n,r,\bm\Phi(\bm\pi_0))$ and $\tilde V_\textnormal{c}(n,r,\bm\Phi(\bm\pi_0))$ for $V_\textnormal{ab}(n,\bm\pi)$ and $V_\textnormal{c}(n,\bm\pi)$ on this line segment, respectively.
For any fixed $\bm\pi_0$ and the spherical angle $\bm\Phi(\bm\pi_0)$, the function $\tilde V_\textnormal{ab}(n,r,\bm\Phi(\bm\pi_0))$ is linear in $r$ because $V_\textnormal{ab}(n,\bm\pi)$ is linear in $\bm\pi$ and $\pi_i$ is linear in $r$ for all $i\in[m_1+m_2]$ as shown in Table \ref{tab:spherical}.
Moreover, $\tilde V_\textnormal{c}(n,r,\bm\Phi(\bm\pi_0))$ is concave in $r$ given $\bm\Phi(\bm\pi_0)$, because $V_\textnormal{c}(n,\bm\pi)$ is concave in $\bm\pi$ by Lemma~\ref{lemma:piecewise_and_concave} and $\pi_i$ is linear in $r$ for all $i\in[m_1+m_2]$.
In the spherical coordinate system, $V_\textnormal{ab}(n,\bm e_{m_1+m_2})\leq V_\textnormal{c}(n,\bm e_{m_1+m_2})$ translates into $\tilde V_\textnormal{ab}(n,0, \bm\Phi(\bm\pi))\leq \tilde V_\textnormal{c}(n,0, \bm\Phi(\bm\pi))$ for all $\bm\pi$.
Then by the concavity of $\tilde V_\textnormal{c}(n,r,\bm\Phi)$ and linearity of $\tilde V_\textnormal{ab}(n,r,\bm\Phi)$ in $r\geq 0$, we can conclude that $\tilde V_\textnormal{c}(n,r,\bm\Phi)$ and $\tilde V_\textnormal{ab}(n,r,\bm\Phi)$ have at most one intersection in the interval $0\leq r(\bm\pi)<r(\bm\pi_0)$ for any $\bm\pi_0\in\mathcal{S}_{-(m_1+m_2)}$ and fixed $\bm\Phi\triangleq\bm\Phi(\bm\pi_0)$.
If such an intersection exists, it is the threshold $\bar r_n$ so that it is optimal to abort the mission when $r<\bar r_n$ and to continue the mission otherwise.
If such an intersection does not exist, the mission would never be successfully executed, and we may set $\bar r_n=\infty$ for any $\bm\pi_0\in\mathcal{S}_{-(m_1+m_2)}$.

We then consider the case $\tilde n<n<\hat n$. 
As in the previous case, we consider $\tilde V_\textnormal{c}(n,r,\bm \Phi)$ and $\tilde V_\textnormal{ab}(n,r,\bm\Phi)$ restricted on the line segment $\mathcal{S}(\bm\pi_0)$ for $0\le r\le r(\bm\pi_0)$ and $\bm\pi_0\in\mathcal{S}_{-(m_1+m_2)}$.
When $\tilde n<n<\hat n$, we have $V_\textnormal{ab}(n,0,\bm\Phi(\bm\pi_0))>\tilde V_\textnormal{c}(n,0,\bm \Phi(\bm\pi_0))$ for all $\bm\pi_0\in\mathcal{S}_{-(m_1+m_2)}$.
By the linearity of $\tilde V_\textnormal{ab}(n,r,\bm\Phi)$ in $r$ and concavity of $\tilde V_\textnormal{c}(n,r,\bm \Phi)$ in $r$, there could be two scenarios: (a) There is no intersection between $\tilde V_\textnormal{ab}(n,r,\bm\Phi)$ and $\tilde V_\textnormal{c}(n,r,\bm\Phi)$ for the fixed $\bm\Phi\triangleq \bm\Phi(\bm\pi_0)$, so the mission abort is never executed; (b) There are two intersections of $\tilde V_\textnormal{ab}(n,r,\bm\Phi)$ and $\tilde V_\textnormal{c}(n,r,\bm\Phi)$, and we denote the two intersections by $\bar r_{n1}$ and $\bar r_{n2}$ such that $0\leq \bar r_{n1}\le \bar r_{n2}\leq r(\bm\pi_0)$ and $\tilde V_\textnormal{ab}(n,r,\bm\Phi)\leq V_\textnormal{c}(n,r,\bm\Phi)$ only when $r\in[\bar r_{n1},\bar r_{n2}]$ for fixed $\bm\Phi\triangleq\bm\Phi(\bm\pi_0)$, which is the region for abort.

We last show that the region $\{\bm\pi\in\mathcal{S}:V_\textnormal{ab}(n,\bm\pi)\leq V_\textnormal{c}(n,\bm\pi)\}$ is convex for all $n=0,\ldots,N$.
Let $\bm\pi^{(1)},\bm\pi^{(2)}\in\mathcal{S}$ be two belief states such that the corresponding optimal action is to abort the mission.
Then we have $V(n,\bm\pi^{(l)})=V_\textnormal{ab}(n,\bm\pi^{(l)})=C_\textnormal{m}+C_\textnormal{s}\kappa(w_n,\bm\pi^{(l)})$ for $l=1,2$. 
For some $\alpha\in(0,1)$, let $\bm\pi'=\alpha\bm\pi^{(1)}+(1-\alpha)\bm\pi^{(2)}$.
By concavity of $V(n,\bm\pi)$ in $\bm\pi$, we have $V(n,\bm\pi')\geq\alpha V(n,\bm\pi^{(1)})+ (1-\alpha) V(n,\bm\pi^{(2)})=V_\textnormal{ab}(n,\bm\pi')$.
Moreover, since $V(n,\bm\pi')\leq V_\textnormal{ab}(n,\bm\pi')$ by definition, we can conclude $V(n,\bm\pi')= V_\textnormal{ab}(n,\bm\pi')$, so the region for abort is convex for all $n=0,\ldots,N-1$.
\Halmos

\subsection{Proof of Corollary~\ref{coro:opt_policy_wo_a1}}
We can follow the same proof of Theorem~\ref{thm:struct_policy} to show that the structural results in Theorem~\ref{thm:struct_policy}(i) and (iii) still hold.
Since we cannot find a time threshold for no abort, we need to replace $\hat n$ with $N$ and then follow the same proof to show that Theorem~\ref{thm:struct_policy}(ii) holds for $n=\tilde n+1,\ldots,N-1$.
\Halmos

\subsection{Proof of Proposition~\ref{prop:no_intermediate}}\label{proof: condition to remove intermediate region}
We first consider condition (i).
By definition of $\hat n$, there exits $\bm\pi\in\mathcal{S}$ such that $\overline{V}_\textnormal{c}(\hat n-1,\bm\pi)>V_\textnormal{ab}(\hat n-1,\bm\pi)$.
This inequality is equivalent to
\begin{align*}
	& (C_\textnormal{s}+C_\textnormal{m})\kappa((N-\hat n+1)\delta+w_N,\bm\pi)>
	C_\textnormal{m}+C_\textnormal{s}\kappa(w_{\hat n-1},\bm\pi),~\exists \bm\pi\in\mathcal{S} \\
	\Leftrightarrow~
	& C_\textnormal{s}\sum_{i=1}^{m_1+m_2}\pi_i\big[{p}_{i,m_1+m_2+1}((N-\hat n+1)\delta+w_N)-{p}_{i,m_1+m_2+1}(w_{\hat n-1})\big] \\
	& +C_\textnormal{m}\sum_{i=1}^{m_1+m_2}\pi_i{p}_{i,m_1+m_2+1}((N-\hat n+1)\delta+w_N)> C_\textnormal{m},~\exists \bm\pi\in\mathcal{S} \\
	\Leftrightarrow~
	& \max_{\bm\pi\in\mathcal{S}}~ C_\textnormal{s}\sum_{i=1}^{m_1+m_2}\pi_i\big[{p}_{i,m_1+m_2+1}((N-\hat n+1)\delta+w_N)-{p}_{i,m_1+m_2+1}(w_{\hat n-1})\big] \\
	& \qquad\quad
	+C_\textnormal{m}\sum_{i=1}^{m_1+m_2}\pi_i{p}_{i,m_1+m_2+1}((N-\hat n+1)\delta+w_N)> C_\textnormal{m}.
\end{align*}
Since $\mathcal{S}$ is a probability simplex, it is clear that the optimal solution to the linear program on the left side of the last inequality is $\bm\pi=\bm e_{l}$ with $l=\arg\max_{i\in[m_1+m_2]} (C_\textnormal{s}+C_\textnormal{m})p_{i,m_1+m_2+1}((N-\hat n+1)\delta+w_N)-C_\textnormal{s}p_{i,m_1+m_2+1}(w_{\hat n-1})$. 
The condition in the proposition statement ensures $l=m_1+m_2$.
Since abort is not executed at periods $n=\hat n,\ldots,N-1$, we have $\overline{V}_\textnormal{c}(\hat n-1,\bm\pi)={V}_\textnormal{c}(\hat n-1,\bm\pi)$ for all $\bm\pi\in\mathcal{S}$; in particular $\overline{V}_\textnormal{c}(\hat n-1,\bm e_{m_1+m_2})={V}_\textnormal{c}(\hat n-1,\bm e_{m_1+m_2})$.
Then $\tilde n= \hat n-1$ by definition of $\tilde n$.

We then prove the statement for condition (ii) by showing that this condition implies condition (i).
Indeed, $p_{i,m_1+m_2+1}(w_{\hat n-1})=0$ for all $i\in[m_1+m_2]$ when $w_{\hat n-1}=0$.
Since $p_{i,m_1+m_2+1}(t)$ is increasing in $i\in[m_1+m_2]$ for all $t\geq 0$ by Lemma~\ref{lemma:monotonous_pi_n}, it is clear that Problem~\eqref{eq:remove_counter} attains the maximum at $i^*=m_1+m_2$.
Condition (i) hence holds, and we have $\tilde n=\hat n-1$.
\Halmos

\subsection{Proof of Proposition~\ref{lemma:preservation}}\label{proof:monotonicity and concavity of PBVI}
%We write $\mathbb{E}^{(\lambda)}$ as $\mathbb{E}$ for convenience.
At each iteration, the function $\widehat V^{(\tau)}(n,\bm\pi)$ is the pointwise minimum of linear functions in $\bm\pi$, so it is concave in $\bm\pi\in\mathcal{S}$.
We then use induction to prove the MLR monotonicity of $\widehat{V}^{(\tau)}(n,\bm\pi)$ in $\bm\pi\in\mathcal{S}$.
At the last stage $N-1$, we have
$\widehat{V}^{(\tau)}_\textnormal{c}(N-1,\bm\pi)=\kappa(\bm\pi,\delta)(C_\textnormal{s}+C_\textnormal{m})+\mathbb{E}^{(\lambda)}\left[\widehat{V}^{(\tau)}(N,\bm\Pi_N)\mid\bm\Pi_{N-1}=\bm\pi \right]$.
Because $V(N,\bm\pi)$ is linear in $\bm \pi$, the point-based value iteration (PBVI) algorithm can give the exact estimate of the terminal cost, i.e., $\widehat{V}^{(\tau)}(N,\bm\pi)=V(N,\bm\pi)$ for all $\tau$, which is MLR increasing in $\bm\pi$ by Lemma~\ref{lemma:monotone}. 
Then following the proof of Lemma~\ref{lemma:monotone}, we can conclude $\mathbb{E}^{(\lambda)}\left[\widehat{V}^{(\tau)}(N,\bm\Pi_N)\mid\bm\Pi_{N-1}=\bm\pi \right]$ is MLR increasing in $\bm\pi$.
This shows $\widehat{V}^{(\tau)}_\textnormal{c}(N-1,\bm\pi)$ is MLR increasing in $\bm\pi$ for all $\tau$.
We can then replace $V(n,\bm\pi)$ with $\widehat{V}^{(\tau)}(n,\bm\pi)$ in the proof of Lemma~\ref{lemma:monotone} and follow the same induction procedure to conclude $\widehat{V}^{(\tau)}(n,\bm\pi)$ is MLR increasing in $\bm\pi$ for $n=0,\ldots,N-2$.
\Halmos

\subsection{Proof of Corollary~\ref{thm:control_limit_CTMC}}
Following a similar procedure as in the proof of Lemma~\ref{lemma:piecewise_and_concave}, we can show that $V_\textnormal{c}(n,\pi)$ is concave in $\pi$ for all $n=0,\ldots,N$.
When $n\leq \tilde n$, we have $V_\textnormal{ab}(n,1)<V_\textnormal{c}(n,1)$.
By linearity of $V_\textnormal{ab}(n,\pi)$ and concavity of $V_\textnormal{c}(n,\pi)$, there is at most one intersection between $V_\textnormal{ab}(n,\pi)$ and $V_\textnormal{c}(n,\pi)$ for $\pi\in[0,1]$, which is the control limit $\underline{\pi}_n$.
When $n> \tilde n$, we have $V_\textnormal{ab}(n,1)\geq V_\textnormal{c}(n,1)$.
Again by linearity of $V_\textnormal{ab}(n,\pi)$ and concavity of $V_\textnormal{c}(n,\pi)$, there could be two scenarios, i.e., (i) there are two intersections between $V_\textnormal{ab}(n,\pi)$ and $V_\textnormal{c}(n,\pi)$ for $\pi\in[0,1]$, which are the thresholds $\underline{\pi}_n$ and $\overline{\pi}_n$, and (ii) there is no intersection between $V_\textnormal{ab}(n,\pi)$ and $V_\textnormal{c}(n,\pi)$, so abort is never optimal.
\Halmos

\subsection{Proof of Proposition~\ref{prop:dim_red}}\label{appen:proof_dim_red}
%The proof largely follows the proof of Lemma~6 in \cite{wang2015multistate}.
%We provide the proof here for completeness.
Note that there are $m_1+m_2=m_1+1$ states for the CTMC when $m_2=1$.
Let $\bm\pi_n^-$ be the realization of the system state before observing the signal $Y_n$ at time $n\delta$.
Given the belief vector $\bm\Pi_{n-1}=\bm\pi_{n-1}=(\pi_{n-1,i})_{i\in\mathcal{X}_1\cup\mathcal{X}_2}$, the $i$th coordinate of $\bm\pi_n^-$ is then written as $\pi_{ni}^-=\sum_{j=1}^{m_1+1}\pi_{n-1,j}p_{ji}(\delta)$.
By noting $p_{m_1+1,i}(\delta)=0$ for all $i\in[m_1]$, 
we can derive $\pi_{ni}^-=r_{n-1}\varphi_i(\bm\Phi_{n-1})$ for all $i\in[m_1]$,
where $r_{n-1}=r_{n-1}(\bm\pi_{n-1})$ and $\bm\Phi_{n-1}=\bm\Phi_{n-1}(\bm\pi_{n-1})$.
We note $\sum_{i=1}^{m_1+1}\pi_{ni}^-< 1$ because it is likely that the system directly fails between two consecutive inspections.
Nevertheless, since it is impossible to observe the signal $Y_n\in[K]$ if the system has failed before $n\delta$, $\sum_{i=1}^{m_1+1}\pi_{ni}^-<1$ does not affect the constraint $\bm\pi_n'\bm 1_{m_1+1}=1$ after observing the realization of $Y_n=y_n\in[K]$.

Assume the realization of $Y_n$ is $y_n\in[K]$.
Then the angle vector $\bm\Phi_n$ after observing $Y_n=y_n$ is 
\begin{align}
	\phi_{ni}
	&=\arccos\frac{\pi_{n,i+1}}{\sqrt{\pi_{n,i+1}^2+\cdots+(\pi_{n,m_1+1}-1)^2+\pi_{n1}^2}} \label{eq:phi1} \\
	&=\arccos\left[
	\frac{\pi_{n,i+1}^-\tilde d_{i+1,y_n}}
	{\sum_{j=1}^{m_1+1}\pi_{nj}^-\tilde d_{j,y_n}}
	\left(
	\left(\dfrac{\pi_{n,i+1}^-\tilde d_{i+1,y_n}}{\sum_{j=1}^{m_1+1}\pi_{nj}^-\tilde d_{j,y_n}}\right)^2
	+\cdots\right.\right. \nonumber \\
	& \qquad\qquad\quad \left.\left.
	+\left(\dfrac{-\sum_{j=1}^{m_1}\pi_{nj}^-\tilde d_{j,y_n}}
	{\sum_{j=1}^{m_1+1}\pi_{nj}^-\tilde d_{j,y_n}}\right)^2
	+\left(\dfrac{\pi_{n1}^-\tilde d_{1,y_n}}{\sum_{j=1}^{m_1+1}\pi_{nj}^-\tilde d_{j,y_n}}\right)^2\right)^{-1/2}\right] \label{eq:phi2} \\
	&=\arccos\dfrac{\pi_{n,i+1}^-}{\sqrt{(\pi_{n,i+1}^-)^2+\cdots+(\sum_{j=1}^{m_1}\pi_{nj}^-)^2+(\pi_{n1}^-)^2}}, \quad i\in[m_1-1], \nonumber
\end{align}
which is irrelevant to $y_n$ and is Equation~\eqref{eq:angle} by routine algebra. 
The last equality follows from $\tilde d_{i1}=\cdots=\tilde d_{iK}$ for all $i\in[m_1]$.
A similar procedure proves the case of $\phi_{n,m_1}$.

Last, we show that $r_n$ can be computed recursively by 
\begin{align*}
	%\label{eq:radius}
	r_n & = \sqrt{\sum_{i=1}^{m_1}\pi_{ni}^2+(\pi_{n,m_1+1}-1)^2} = \left[\sum_{i=1}^{m_1}\left(\frac{\pi_{ni}^-\tilde d_{i,y_n}}{\sum_{j=1}^{m_1+1}\pi_{nj}^-\tilde d_{j,y_n}}\right)^2+
	\left(
	\frac{-\sum_{j=1}^{m_1}\pi_{nj}^-\tilde d_{j,y_n}}
	{\sum_{j=1}^{m_1+1}\pi_{nj}^-\tilde d_{j,y_n}}
	\right)^2\right]^{\frac{1}{2}} \\
	& = \left[\sum_{i=1}^{m_1}\left(\frac{\pi_{ni}^-}{\sum_{j=1}^{m_1}\pi_{nj}^-+\pi_{n,m_1+1}^-d_{2,y_n}/d_{1,y_n}}\right)^2+
	\left(
	\frac{-\sum_{j=1}^{m_1}\pi_{nj}^-}
	{\sum_{j=1}^{m_1}\pi_{nj}^-+\pi_{n,m_1+1}^-d_{2,y_n}/d_{1,y_n}}
	\right)^2\right]^{\frac{1}{2}} \\
	& = \frac{\sqrt{\sum_{i=1}^{m_1}(\pi_{ni}^-)^2+(\sum_{i=1}^{m_1}\pi_{ni}^-)^2}}{\sum_{i=1}^{m_1}\pi_{ni}^-+\pi_{n,m_1+1}^-d_{2,y_n}/d_{1,y_n}} \\
	& =\frac{r_{n-1}\sqrt{\sum_{i=1}^{m_1}\varphi_i^2(\bm\Phi_{n-1})+(\sum_{i=1}^{m_1}\varphi_i(\bm\Phi_{n-1}))^2}}
	{r_{n-1}\sum_{i=1}^{m_1}\varphi_{i}(\bm\Phi_{n-1})+(r_{n-1}\varphi_{m_1+1}(\bm\Phi_{n-1})+p_{m_1+1,m_1+1}(\delta))d_{2,y_n}d_{1,y_n}^{-1}} \\
	& =\frac{\sqrt{\sum_{i=1}^{m_1}\varphi_i^2(\bm\Phi_{n-1})+(\sum_{i=1}^{m_1}\varphi_i(\bm\Phi_{n-1}))^2}}
	{\sum_{i=1}^{m_1}\varphi_{i}(\bm\Phi_{n-1})+(\varphi_{m_1+1}(\bm\Phi_{n-1})+p_{m_1+1,m_1+1}(\delta)r_{n-1}^{-1})d_{2,y_n}d_{1,y_n}^{-1}} 
\end{align*}
given $Y_n=y_n\in[K]$.
\iffalse
where
\begin{multline*}
	\varphi_{m_1+1}(\bm\Phi)\triangleq\cos \phi_1(p_{2,m_1+1}(\delta)-p_{1,m_1+1}(\delta)) \\
	+\sum_{i=3}^{m_1+1}\cos\phi_{i-1}\prod_{k=1}^{i-2}\sin\phi_k(p_{i,m_1+1}(\delta)-p_{1,m_1+1}(\delta))-
	p_{m_1+1,m_1+2}(\delta).
\end{multline*}
\fi
The equation follows from $\pi_{ni}^-=\sum_{j=1}^{m_1+1}\pi_{n-1,j}p_{ji}(\delta)$, the transform formula in Table~\ref{tab:spherical}, the relationship between $\bm\pi_{n}^-$ and $\bm\pi_n$ in \eqref{eq:phi1} and \eqref{eq:phi2}, and 
$\pi_{n,m_1+1}^-
=1+(\pi_{n,m_1+1}^--1)
=1+\sum_{j=1}^{m_1+1}\pi_{n-1,j}(p_{j,m_1+1}(\delta)-1)
=1+r_{n-1}\varphi_{m_1+1}(\bm\Phi_{n-1})+p_{m_1+1,m_1+1}(\delta)-1
=r_{n-1}\varphi_{m_1+1}(\bm\Phi_{n-1})+p_{m_1+1,m_1+1}(\delta)$.
\Halmos

\subsection{Proof of Theorem~\ref{thm:multi_opt}}\label{appen:proof_thm:multi_opt}
{
Similar to the single-mission case, we let $\overline{V}_\text{c}(n,\bm\pi_n)$ be the cost-to-go if we continue the mission at time $n\delta$ with belief state $\bm\pi_n$ until the end of the mission.
For any $n$, let $\tilde l(n)=\min\{l\in[L]:\sum_{l'=1}^{l}N_{l'}>n\}$.
We can readily compute
\[
\overline{V}_\text{c}(n,\bm\pi_n)=
C_\text{s}\kappa((N-n)\delta+w_N,\bm\pi_n)+
\sum_{l=\tilde l(n)}^L C_\text{m}^{(l)}\kappa\left(\left(\sum_{l'=1}^l N_{l'}-n\right)\delta+w_N\mathbbm{1}\{l=L\},\bm\pi_n\right).
\]
We know $n\geq\sum_{l'=1}^lN_{l'}$ if mission $l\in[L-1]$ has been completed.
Since $\kappa(t,\bm\pi)$ is increasing in $t$ and MLR increasing $\bm\pi$, an upper bound of $\overline{V}_\text{c}(n,\bm\pi)$ for $n\geq \sum_{l'=1}^lN_{l'}$ is
\begin{align*}
	\overline{V}_\text{c}\left(\sum_{l'=1}^lN_{l'},\bm e_{m_1+m_2}\right)=&~
	C_\text{s}\left[1-e^{-\lambda(T-\delta\sum_{l'=1}^lN_{l'}+w_N)}\right] \\
	&~+
	\sum_{l'=l+1}^L C_\text{m}^{(l')}
	\left\{1-\exp\left[-\lambda\left(\sum_{l''=l+1}^{l'} N_{l''}\delta+w_N\mathbbm{1}\{l'=L\}\right)\right]\right\}.
\end{align*}
By the same way, a lower bound of $V_\text{ab}(n,\bm\pi)$ with $n\geq \sum_{l'=1}^l N_{l'}$ is $C_\text{s}\bm e'_{1}\exp(\mathbf{Q}w_{\sum_{l'=1}^lN_{l'}})\bm e_{m_1+m_2+1}+C_\text{m}^{(L)}$.
The mission is never aborted after completing mission $l$ if $V_\text{c}(\sum_{l'=1}^lN_l,\bm e_{m_1+m_2})\leq C_\text{m}^{(L)}+C_\text{s}\bm e'_{1}\exp(\mathbf{Q}w_{\sum_{l'=1}^lN_{l'}})\bm e_{m_1+m_2+1}$, which is 
\begin{multline}\label{eq:no_abort_after_l}
	\sum_{l'=l+1}^L C_\text{m}^{(l')}\left\{\mathbbm{1}\{l'<L\}-\exp\left[-\lambda\left(\sum_{l''=l+1}^{l'}N_{l''}\delta+w_N\mathbbm{1}\{l'=L\}\right)\right]\right\} \\
	\leq C_\text{s}\left[\bm e_1'\exp\left(\mathbf{Q}w_{\sum_{l'=1}^lN_{l'}}\right)\bm e_{m_1+m_2+1}+\exp\left(-\lambda\left(T-\delta\sum_{l'=1}^{l}N_{l'}+w_N\right)\right)-1\right].
\end{multline}
When $l=0$, $\sum_{l'=1}^0 N_l=0$. 
The mission is never aborted in this case, and \eqref{eq:no_abort_after_l} simplifies to \eqref{eq:no_abort_multi}.
	
Since the expected repair cost for aborting or stopping system at $n\delta$ is $C_\text{r}\bm\pi' \tilde{\mathbf{P}}(w_n)\cdot(\bm 0_{m_1}',\bm 1_{m_2}')'$, which is linear in $\bm\pi$, 
we can use the same proof for Lemma~\ref{lemma:piecewise_and_concave} to show that $V^{(\lambda)}(n,\bm\pi)$, $V_\text{ab}(n,\bm\pi)$, and $V_\text{c}(n,\bm\pi)$ are all concave in $\bm\pi$ for all $n$.
However, since $w_n$ may not be nondecreasing in $n$, we cannot guarantee existence of the thresholds $\hat n$ and $\tilde n$.
Nevertheless, we can examine the values of $V_\text{ab}(n,\bm e_{m_1+m_2})$ and $V_\text{c}(n,\bm e_{m_1+m_2})$ for each $n$.
If $V_\text{ab}(n,\bm e_{m_1+m_2})>V_\text{c}(n,\bm e_{m_1+m_2})$, the optimal abort policy is of the form as in Theorem~\ref{thm:struct_policy}(i);
otherwise, it is of the form as in Theorem~\ref{thm:struct_policy}(ii), both proved using the concavity of $V^{(\lambda)}(n,\bm\pi)$, $V_\text{ab}(n,\bm\pi)$, and $V_\text{c}(n,\bm\pi)$.
Again using the concavity of value functions, the region in which abort is optimal can be shown to be convex.}
\Halmos

\section{Implementation Details}\label{appen:implement}
We first illustrate the details of the PBVI algorithm.
At the $\tau$th iteration of the algorithm, the backup step first constructs $\widehat{\mathcal{A}}_{n}^{(\tau)}$ based on $\{\bm\pi^{(l)}_n\}_{l\in[L_\tau]}$ in a backward induction way.
Specifically, we have $V(n,\bm\pi)=\overline{V}_\textnormal{c}(n,\bm\pi)$ for $n=\hat n,\ldots,N$ by Theorem~\ref{thm:threshold_n}, which is a linear function in $\bm\pi$.
All the sets ${\mathcal{A}}_n$, $n=\hat n,\ldots,N$, are hence a singleton, and the only element therein can be determined by Equation~\eqref{eq:upper_V}.
We simply let $\widehat{\mathcal{A}}^{(\tau)}_n={\mathcal{A}}_n$ for all $n\geq \hat n$.
When $n=\hat n-1$, we have two choices, i.e., to abort or to continue the mission, and $V(n,\bm\pi)=\min\{V_\text{ab}(n,\bm\pi),V_\text{c}(n,\bm\pi)\}$.
The cost-to-go $V_\text{ab}(n,\bm\pi)$ can be exactly computed since the map $\bm\pi \mapsto V_\textnormal{ab}(\hat n-1,\bm\pi)$ in \eqref{eq:v_ab} is linear.
% This means that $\widehat{\mathcal{A}}_{\hat n-1}^{(\tau)}$ is a singleton, and the only element therein can be obtained by \eqref{eq:v_ab}.
On the other hand,
exactly computing $V_\textnormal{c}(\hat n-1,\bm\pi)$ needs to solve the Bellman equation in~\eqref{eq:v_c}, which is computationally intensive.
Hence, we approximately compute $V_\textnormal{c}(\hat n-1,\bm\pi)$ by replacing $V(n,\bm\pi)$ with $\widehat{V}^{(\tau)}(n,\bm\pi)$ in~\eqref{eq:v_c}, yielding an approximation
\begin{multline}\label{eq:approx_N-1}
	\widehat{V}^{(\tau)}_\textnormal{c}(\hat n-1,\bm\pi)=\kappa(\bm\pi,\delta)(C_\textnormal{s}+C_\textnormal{m}) \\
	+\sum_{k=1}^K\min_{\bm\alpha_{\hat n}\in\widehat{\mathcal{A}}^{(\tau)}_{\hat n}}\sum_{i=1}^{m_1+m_2}\sum_{j=1}^{m_1+m_2}\alpha_{\hat n j}\pi_i p_{ij}(\delta) (d_{1k}\mathbbm{1}\{j\leq m_1\}+d_{2k}\mathbbm{1}\{j>m_1\}),
\end{multline}
where $\alpha_{\hat n j}$ is the $j$th coordinate of $\bm\alpha_{\hat n}$.
Since $\bm\pi\mapsto \kappa(\bm\pi,\delta)$ is linear, 
Equation~\eqref{eq:approx_N-1} implies both $\widehat V_\text{c}(n,\bm\pi)$ and $\widehat{V}^{(\tau)}(n,\bm\pi)\triangleq \min\{V_\text{ab}(n,\bm\pi),\widehat V_\text{c}(n,\bm\pi)\}$ are piecewise linear and concave in $\bm\pi$.
Based on $V_\text{ab}(n,\bm\pi)$ in \eqref{eq:v_ab} and $\widehat V_\text{c}(n,\bm\pi)$ in \eqref{eq:approx_N-1}, the set $\widehat{\mathcal{A}}_{\hat n-1}$ can be readily obtained.
This finishes the backward induction step for $n=\hat n-1$.
We keep repeating this procedure from $n=\hat n-2$ to $n=0$ to complete the backup step.

The expansion step then expands the sets $\{\bm\pi_n^{(l)}\}_{l\in [L_\tau]}$ by adding new elements $\bm\pi_n^{(l)}$, $l=L_\tau+1,\ldots,L_{\tau+1}$, for each $n\in[\hat n-1]$.
The expanded sets $\{\bm\pi_n^{(l)}\}_{l\in[L_{\tau+1}]}$ will then be used in the backup step at the $(\tau+1)$th iteration.
To ensure a good approximation, we need to find states $\{\bm\pi_n^{(l)}\}_{l=L_\tau+1}^{L_{\tau+1}}$ that are more likely to reach for the partially observable Markov decision process (POMDP) \citepappendix{pineau2003point,liu2022machine}.
Hence, we resort to the CTMC model in Section~\ref{subsec:pomdp_phase} to simulate system states.
Specifically, for each $l\in[L_\tau]$ and $n=0,\ldots,\hat n-2$, we run a batch of simulations starting from the state $\bm\pi_{n}^{(l)}$ and find the belief state at the next decision epoch.
To this end, we first sample an underlying state from each belief state $\bm\pi_{n}^{(l)}$ and then sample the state at the next decision epoch and the associated observation.
The belief state at the next decision epoch is computed based on the Bayesian update.
For all states simulated from $\bm\pi_{n}^{(l)}$, we keep the one that is farthest away from $\mathcal{A}^{(\tau)}_{n+1}$ as with \citeappendix{pineau2003point} and denote this state by $\bm\pi_{n+1}^{(L_\tau+l)}$.
We expand the set of belief states by adding all the above kept states.
In this way, we double the size of $\{\bm\pi_n^{(l)}\}_{l\in[L_\tau]}$, $n\in[\hat n-1]$, after each iteration, i.e., $L_{\tau+1}=2L_{\tau}$.
When the iteration number $\tau$ is larger than $Z_2$, we discard all the belief states following Algorithm~\ref{algorithm:PBVI}.
This completes the expansion step.

We next illustrate the binary search to find the time threshold $\hat n$.
The detailed procedures are provided in Algorithm~\ref{algo:binary}.

\begin{algorithm}
	\caption{Binary search to find the threshold $\hat n$.}
	\begin{algorithmic}[1]
		{\renewcommand\baselinestretch{1}\selectfont
		\Require
		The functions $\overline V_\textnormal{c}(n,\bm\pi)$ and $V_\textnormal{ab}(n,\bm\pi)$
		\Ensure
		The time threshold $\hat n$ for abort
		\State Init $\underline{n}\leftarrow 0$, $\overline{n}\leftarrow N$, and $\hat n\leftarrow \lfloor (\overline{n}+\underline{n})/2 \rfloor$
		\While{$\overline{n}> \underline{n}+1$}
		\State Set $\hat n\leftarrow \lfloor (\overline{n}+\underline{n})/2 \rfloor$
		\For{$\bm\pi=\bm e_1,\ldots,\bm e_{m_1+m_2}$}
		\If{$\overline{V}_\textnormal{c}(\hat n,\bm\pi)-{V}_\textnormal{ab}(\hat n,\bm\pi)> 0$}
		\State Set $\underline{n}\leftarrow \hat n$\;
		Break the for-loop
		\EndIf
		\EndFor
		\State Set $\overline{n}\leftarrow\hat n$
		\EndWhile
		\State Set $\hat n\leftarrow\underline{n}$
		\For{$\bm\pi=\bm e_1,\ldots,\bm e_{m_1+m_2}$}
		\If{$\overline{V}_\textnormal{c}(\underline n,\bm\pi)-{V}_\textnormal{ab}(\underline n,\bm\pi)> 0$}
		\State Set $\hat{n}\leftarrow \overline n$
		\State Break the for-loop
		\EndIf
		\EndFor
		\par}
	\end{algorithmic}
	\label{algo:binary}
\end{algorithm}

We then illustrate the optimization for the special cases of our model.
When $m_1=m_2=1$, the failure process follows a CTMC.
We discretize the belief space $[0,1]$ into $\{0,\Delta,2\Delta,\ldots,1-\Delta,1\}$, where $\Delta$ is the granularity parameter. 
We then obtain the control limits $\underline{\pi}_n,\overline{\pi}_n$ at each decision epoch $n\delta$, $n=0,\ldots,N-1$ by the backward algorithm summarized in Algorithm~\ref{algo:backward}.
When $m_1>1$ and $m_2=1$, we need to first compute $(\bm\Phi_n)_{n=0}^{N}$ by \eqref{eq:angle}. 
We then discretize the interval $[0,-(\cos \phi_{n,m_1}\prod_{k=1}^{m_1-1}\sin \phi_{nk})^{-1}]$ with granularity $\Delta$ for each $n=0,\ldots,N$.
Nearly the same as Algorithm~\ref{algo:backward}, we can use standard backward induction to exactly solve the corresponding Bellman equation after discretization.
%When $m_1=1$ and $m_2>1$, we can also discretize the interval for the radius $r$ under the new spherical coordinate system in Table~\ref{tab:spherical_new}.
%The optimization procedure is in the same spirit as the case that $m_1>1$ and $m_2=1$, and thus details are omitted for brevity.

\begin{algorithm}
	\caption{Backward induction algorithm to calculate the control limits ${\underline{\pi}_n},\overline{\pi}_n$.}
	\begin{algorithmic}[1]
		{\renewcommand\baselinestretch{1}\selectfont
		\Require
		$\mathbf{Q},\mathbf{D},C_{\text{m}},C_{\text{s}},\delta,N,\Delta$
		\Ensure
		The control limits $\underline{\pi}_n,\overline{\pi}_n$ at each decision epoch $n\delta$, $n=0,1,\ldots,N-1$
		\For{$\pi_N\in\{0,\Delta,2\Delta,\ldots,1-\Delta,1\}$}
		\State Calculate $V(N,\pi_N)$ based on \eqref{eq:bellman}
		\EndFor
		\For{$n\leftarrow N-1,\ldots,0$}
		\For{$\pi_{n}\in\{0,\Delta,2\Delta,\ldots,1-\Delta,1\}$}
		\State Calculate $V(n,\pi_{n})$ based on \eqref{eq:bellman_base}
		\EndFor
		\State Find the thresholds ${\underline{\pi}_n},\overline{\pi}_n$ based on the value function
			\EndFor
			\par}
	\end{algorithmic}
	\label{algo:backward}
\end{algorithm}

We last demonstrate the optimization of the benchmark methods used in Section~\ref{sec:num}.
For the optimal $\mathcal{C}$-policy and $\mathcal{R}$-policy, {two brute-force searches are used to separately optimize $(\check{m},\check{N})$ and $p$.}
To obtain the optimal $\mathcal{M}$-policy, we can formulate a POMDP as in Section~\ref{subsec:special_CTMC} under the three-state CTMC approximation.
As shown in Theorem~\ref{thm:control_limit_CTMC}, the optimal abort policy has a control-limit structure on the belief state $[0,1]$.
Specifically, we let the transition rates be $q_{01}=\nu/m_1=4.01\times10^{-3}$, $q_{02}=\lambda=10^{-3}$, and $q_{12}=[\alpha\Gamma(1+1/\beta)]^{-1}=1.04\times10^{-3}$, where $\Gamma(\cdot)$ is the gamma function, to match the mean time of transition between states for the original failure process. 
We then discretize the state space with granularity $\Delta=0.01$ and use standard backward induction as shown in Algorithm~\ref{algo:backward} to solve the POMDP.
The threshold to abort at each period can be readily obtained for the optimal $\mathcal{M}$-policy.
{The benchmark with the one-phase approximation uses similar procedures to the $\mathcal{M}$-policy while keeping the first $m_1$ phases unchanged.}

\section{Additional Numerical Settings and Results}
\subsection{Numerical Settings}\label{appen:num_setting}
When comparing our modified PBVI algorithm with the classical one in numerical experiments, we randomly sample 100 belief states at the end of each iteration of our modified PBVI algorithm.
For each sampled belief state, we evaluate the approximate value function under the current iteration of the modified PBVI algorithm and compare it with that already obtained from the classical PBVI algorithm.
The modified PBVI terminates if for all the 100 sampled belief states, the absolute differences between the two approximate values are smaller than one.

In Section~\ref{subsec:num_setting}, we assume $T_{12}$ follows an Erlang distribution, with shape $m_1$ and rate $\nu$.
We set the shape parameter of $T_{12}$ as $m_1=2$ \citepappendix{khaleghei2021optimal}, which is commonly used for machinery systems.
Then we determine $\nu=8.01\times 10^{-3}$ to make $\mathbb{E}[T_{12}]$ the same as that in \citeappendix{yang2019designing}.
Since the Erlang distribution can capture the absorption time of a CTMC, we construct a CTMC with transition rate matrix $\mathbf{Q}=(q_{ij})_{i,j\in\mathcal{X}}$, where
\[
q_{ij}=
\begin{cases}
	\nu, & j=i+1,~i\in\mathcal{X}_1,\\
	\zeta, & i\in\mathcal{X}_1,~ j=m_1+m_2+1,\\
	[1-p_{i-m_1}(\lambda)]\lambda, & j=i+1,~i,j\in\mathcal{X}_2,\\
	p_{i-m_1}(\lambda)\lambda, & i\in\mathcal{X}_2\backslash\{m_1+m_2\},~ j=m_1+m_2+1, \\
	\lambda, & i=m_1+m_2,~j=m_1+m_2+1,
\end{cases}
\]
$q_{ij}=0$ for other $i\neq j$, and $q_{ii}=-\sum_{j\in\mathcal{X}:j\neq i}q_{ij}$.
The CTMC starts from the first state with probability one, i.e., $X^{(\lambda)}(0)=1$.
The optimal mission abort policy can be obtained by solving $V(0,\bm e_1)$ from the Bellman recursion.
Lemmas~1 and 2 in \citeappendix{khaleghei2021optimal} suggest that all structural results in Section~\ref{sec:property} still hold if we formulate a POMDP based on this CTMC.

{Moreover, Assumption~\ref{assump:nondecreasing} fails when $F$ is a mixture of two Weibull distributions.
In a real application, if a mixture of Erlang distribution $F^{(\lambda)}$ is used to approximate a general positive distribution $F$, the values of both $\lambda>0$ and $m_2\in\mathbb{N}_+$ need to be determined.
In our numerical experiment, for a fixed value of $m_2$, we use the commonly used moment-matching method to determine $\lambda$ such that $F^{(\lambda)}$ and $F$ have the same mean value.
As such, we only need to choose the value of $m_2$ for optimization. 
An intuitive way to choose a reasonable $m_2$ value is to plot the cumulative distribution functions (CDFs) of the true $F$ and the mixture of Erlang for approximation in the same graph and visually compare the two curves.
We do this comparison for both the Weibull distribution in \eqref{eq:weibull} and the mixture distribution in \eqref{eq:mixture}. 
For the Weibull distribution, we examine $m_2\in\{10, 20, 30, 40\}$ and the results are provided in Figure~\ref{fig:inspection1}.
For the mixture distribution, we examine $m_2\in\{10,30,50,70\}$, with the results provided in Figure~\ref{fig:inspection}.
In the first case, $m_2=20$ has provided an accurate approximation to $F$, and further increasing $m_2$ only leads to marginal improvement.
In the second case, the lack of fit mainly comes from the CDF at the neighborhood of time $t=50$. The choice of $m_2=50$ seems to provide a big improvement on the fit in this area compared to $m_2=30$, while increasing $m_2$ to 70 improves the approximation marginally.
The observations validate our argument of choosing $m_2$ by visually inspecting the true CDF and the approximating CDF.}

\begin{figure}
	\centering
	\subfigure[$m_2=10$]{\includegraphics[scale=0.215]{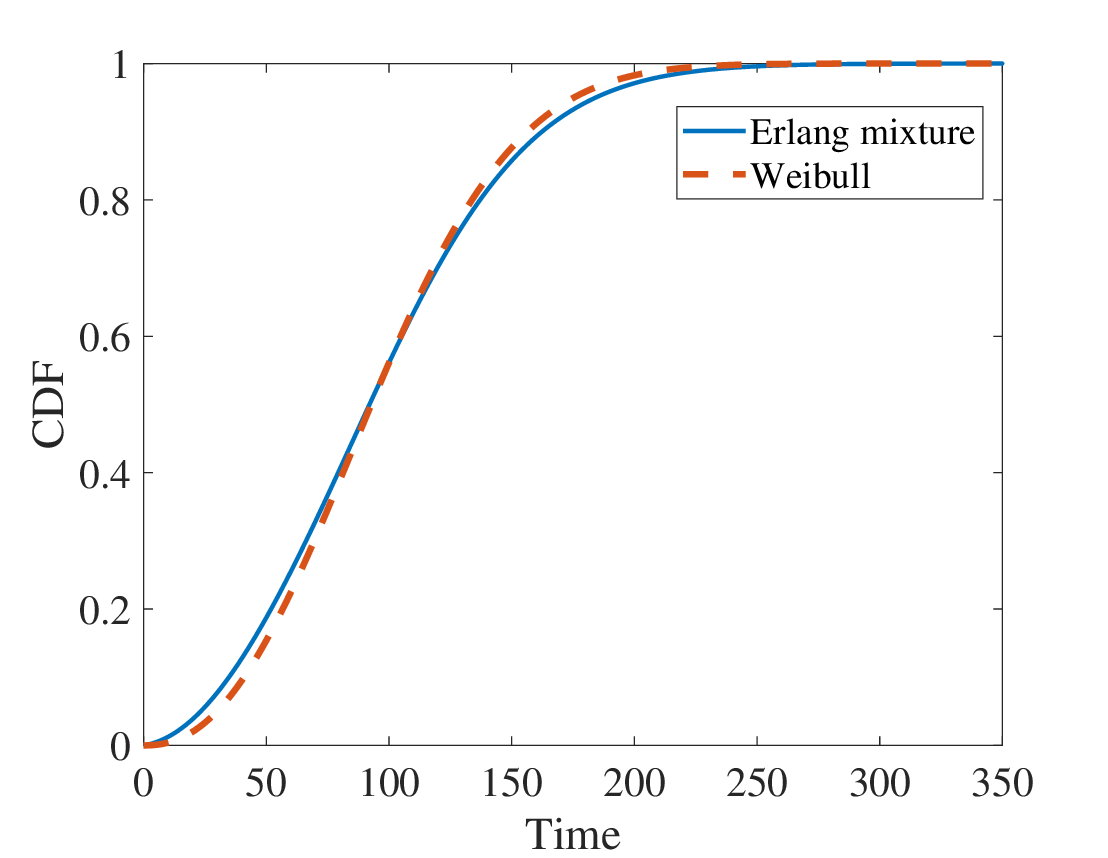}}
	\subfigure[$m_2=20$]{\includegraphics[scale=0.215]{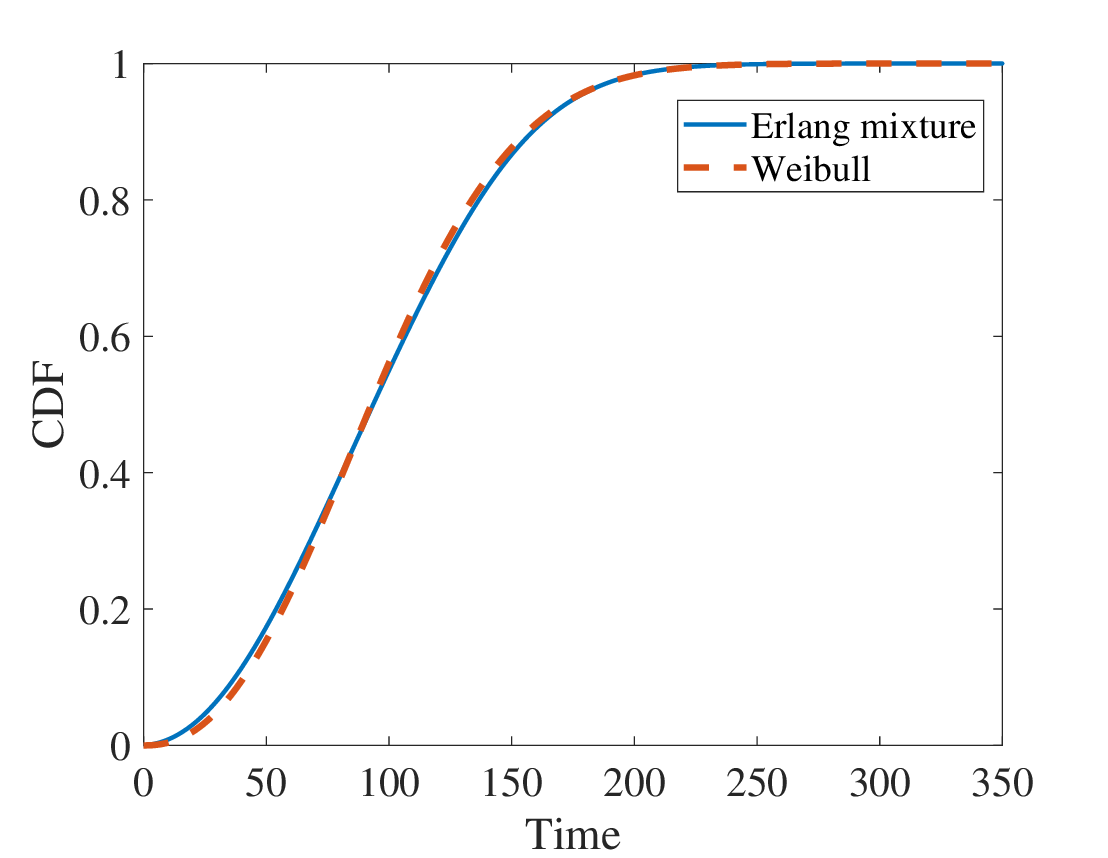}}
	\subfigure[$m_2=30$]{\includegraphics[scale=0.215]{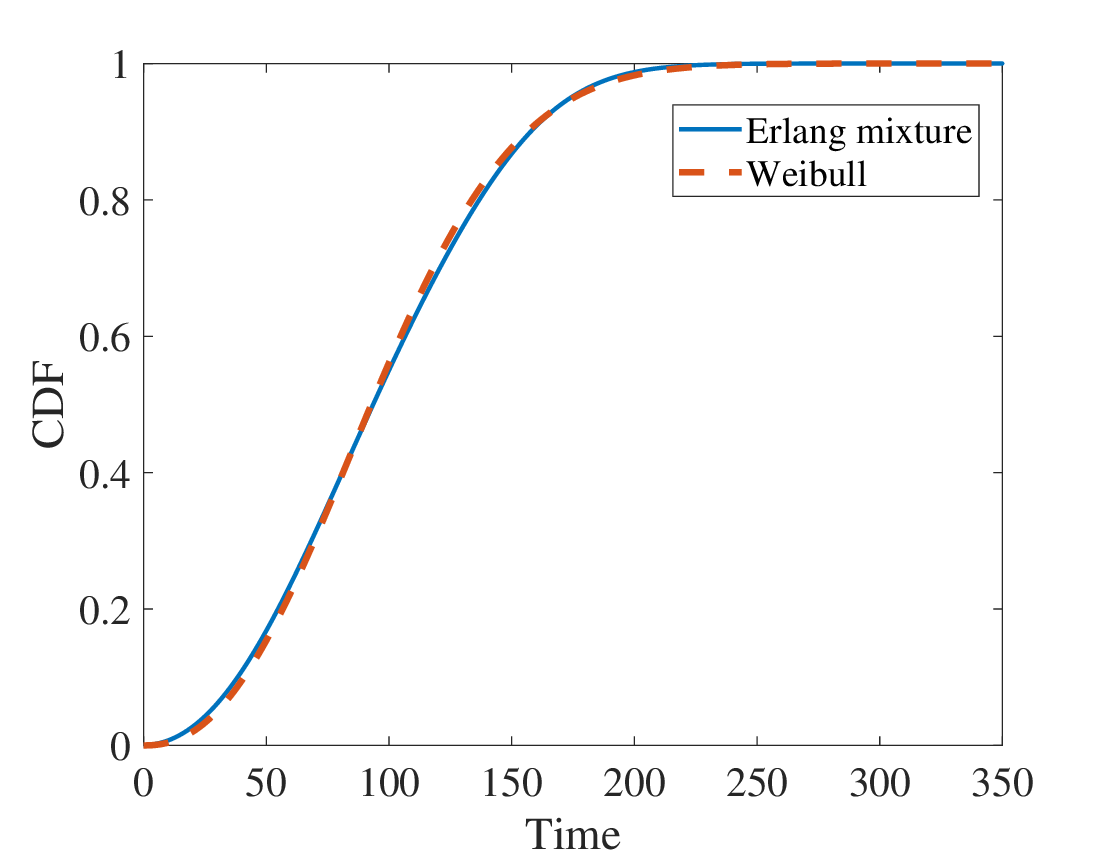}}
	\subfigure[$m_2=40$]{\includegraphics[scale=0.215]{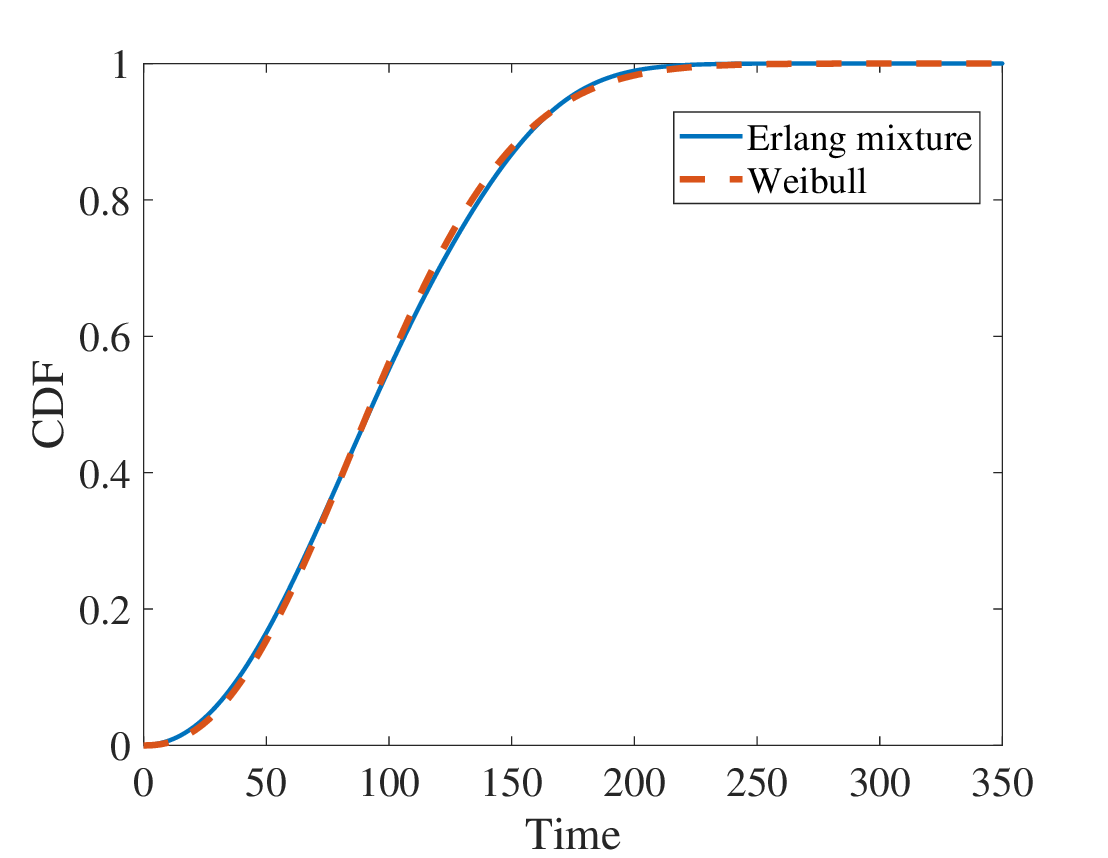}}
	\caption{Approximating the distribution $F$ of $T_{23}$ in \eqref{eq:weibull} by the Erlang mixture distribution $F^{(\lambda)}(\cdot)$ under various values of $m_2$.} 
	\label{fig:inspection1}
\end{figure}

\begin{figure}
	\centering
	\subfigure[$m_2=10$]{\includegraphics[scale=0.215]{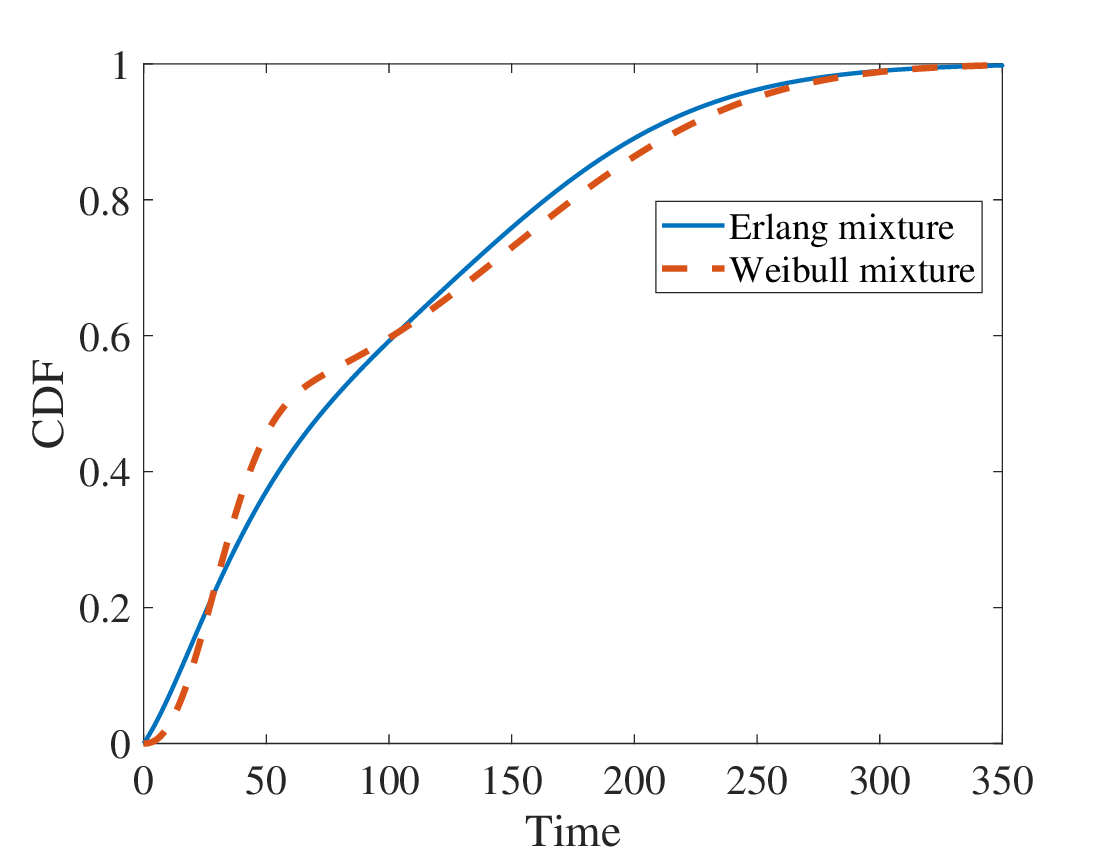}}
	\subfigure[$m_2=30$]{\includegraphics[scale=0.215]{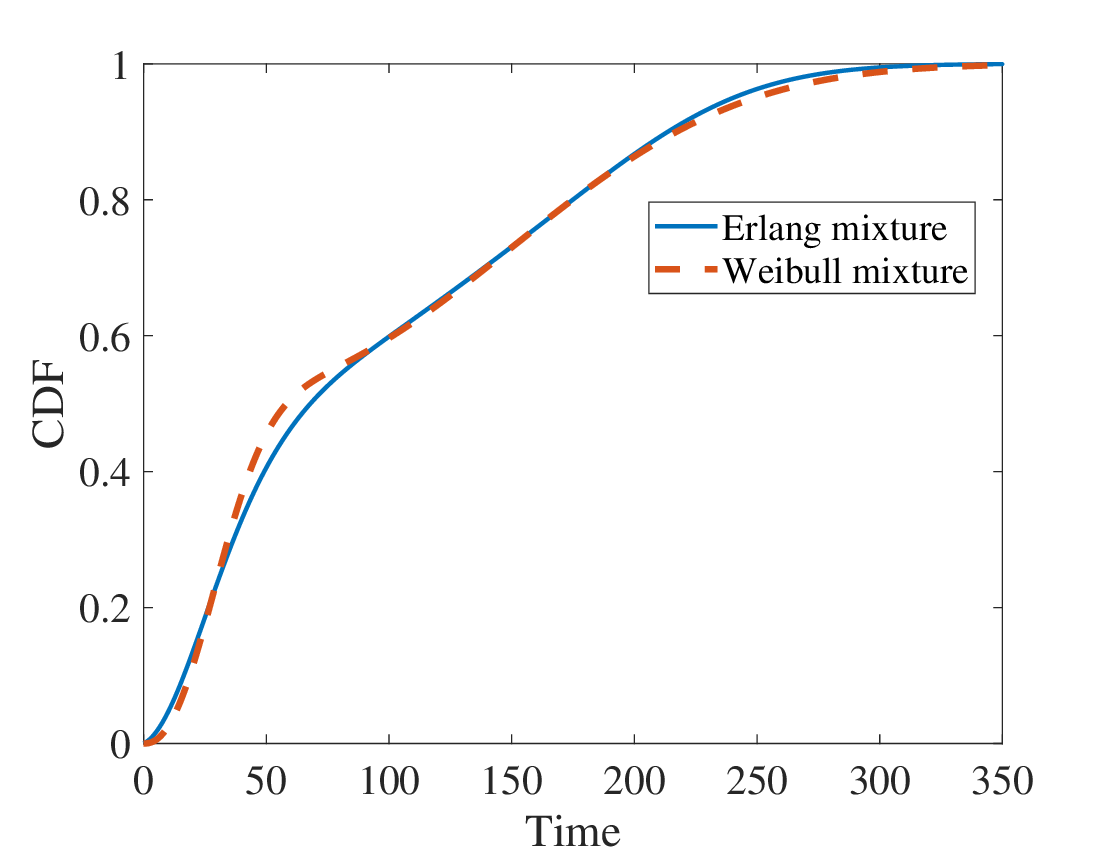}}
	\subfigure[$m_2=50$]{\includegraphics[scale=0.215]{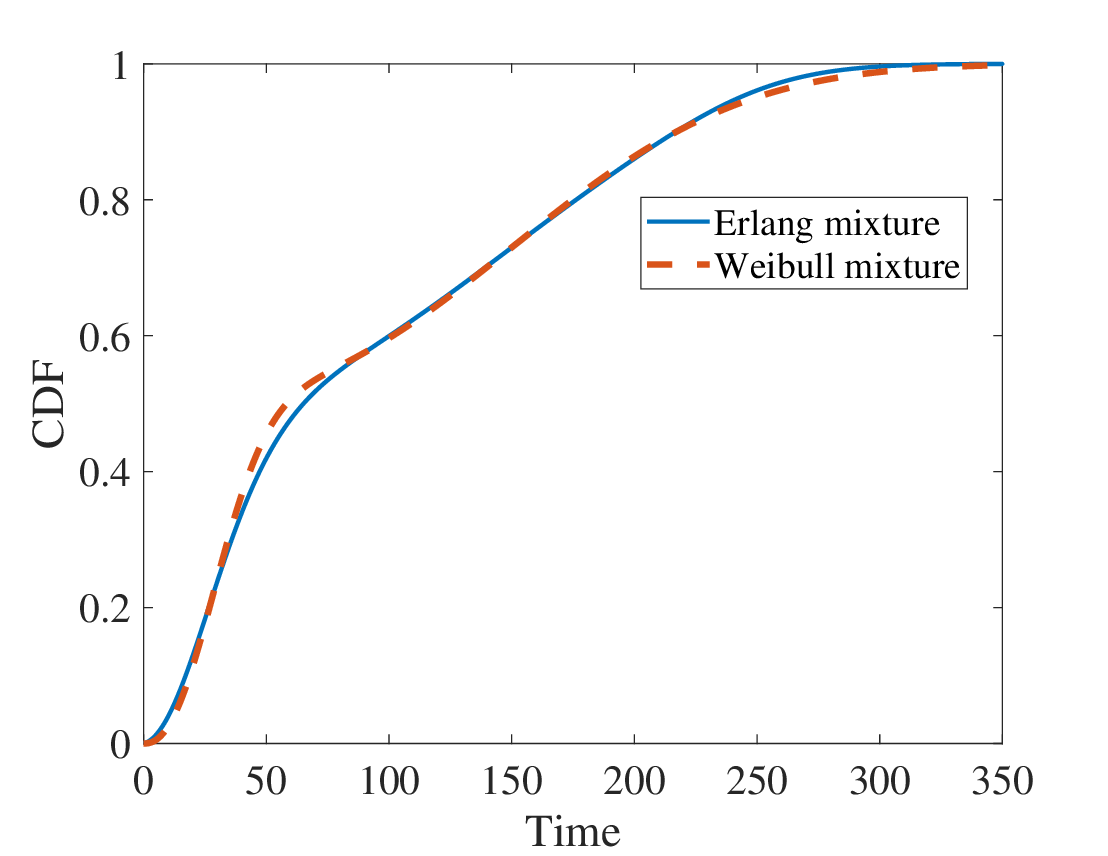}}
	\subfigure[$m_2=70$]{\includegraphics[scale=0.215]{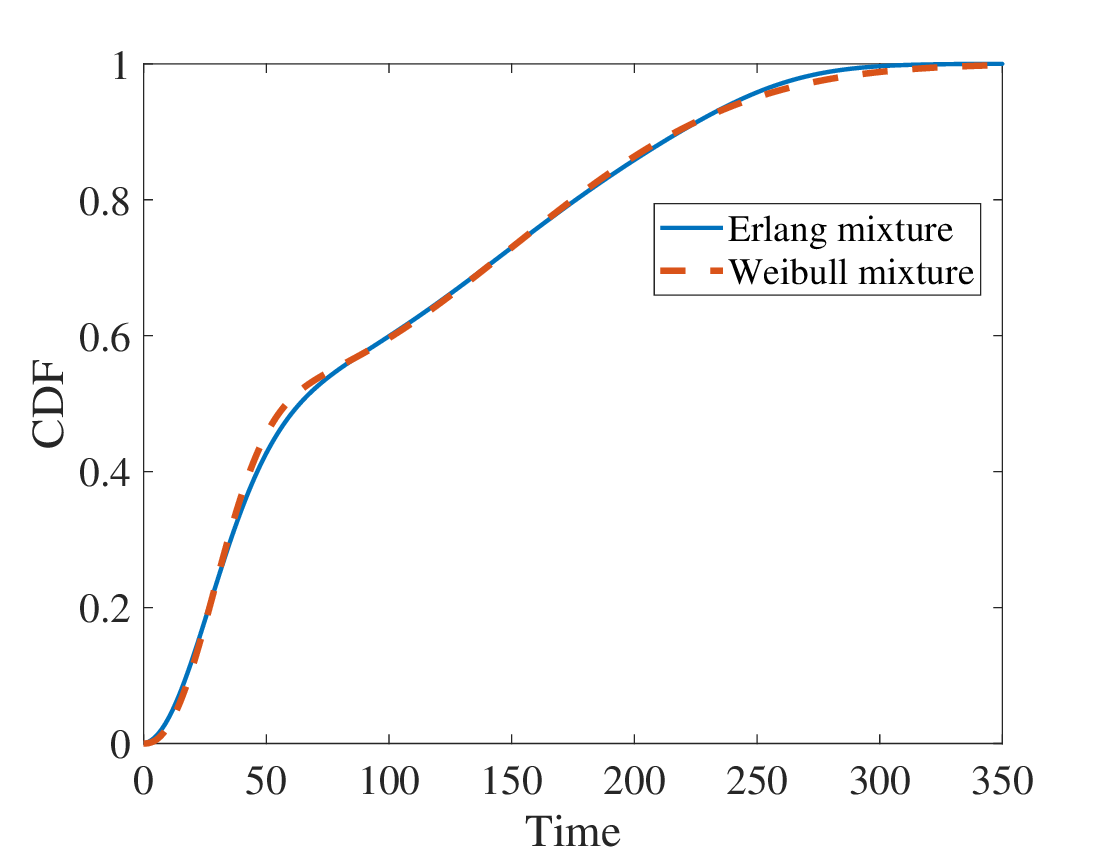}}
	\caption{Approximating the distribution $F$ of $T_{23}$ in \eqref{eq:mixture} by the Erlang mixture distribution $F^{(\lambda)}(\cdot)$ under various values of $m_2$.} 
	\label{fig:inspection}
\end{figure}

{To compare the two PBVI algorithms in terms of both running time and accuracy, we consider a problem instance with $m_1=1$ and $m_2=2$.
Specifically, the transition from the healthy to defective state follows an exponential distribution, and that from the defective to failure state follows a two-phase mixture of Erlang distribution.
This stochastic process can be represented by a four-state CTMC with the transition rate matrix given by
\[
\mathbf{Q}=
\left[
\begin{array}{cccc}
	-2.75\times10^{-3} & 2.29\times 10^{-3} & 0 & 4.59\times 10^{-4} \\
	0 & -1.04\times10^{-2} & 6.92\times 10^{-3} & 3.46\times 10^{-3} \\
	0 & 0 & -2.86\times10^{-2} & 2.86\times10^{-2} \\
	0 & 0 & 0 & 0 
\end{array}
\right].
\]
The state-observation matrix is the same as that in Section~\ref{subsec:num_setting}.
Since the failure state is observable, the belief state is three-dimensional, i.e., $\{\bm\pi\in\mathbb{R}^3_+:\pi_1+\pi_2+\pi_3= 1\}$.
We discretize $[0,1]$ for each dimension of $\bm\pi$ into $1,000$ states, resulting in $1,000^3=10^9$ states in total.
It takes around three days for a computing server to use backward induction to find the optimal solution for this problem instance.
%We fix the random seed so the two PBVI algorithms generate the same belief states in the expansion step in each simulation replication.
For our modified PBVI algorithm, we set $L_1=2$, $Z_1=30$, and $Z_2=10$.

Last, for the numerical study in Section~\ref{subsec:num_multi}, the rescue time $w_n$ can be computed by the law of cosines based on the locations of the depot and sites 1 to 3.
The rescue time is set to be
\[
w_n=
\begin{cases}
	n, & 0\leq n<25, \\
	25, & 25\leq n<60, \\
	\sqrt{625+(n-60)(n-85)}, & 60\leq n<85, \\
	25, & 85\leq n<110, \\
	\sqrt{625+(n-110)(n-135)}, & 110\leq n<135, \\
	25, & 135\leq n\leq 185, \\
\end{cases}
\]
to assess the model performance in a multi-task setting.}

\subsection{Sensitivity Analysis}\label{appen:sen}
For a comprehensive comparison, we conduct sensitivity analysis on several model parameters.
{We first consider the case when $F$ is Weibull.}
Since the mission failure cost is subjective, we fix $C_\text{s}$ and examine different values of $C_\text{m}$ for $C_{\text{m}}/C_{\text{s}}\in\{0.25,0.5,1,2,4\}$.
Table~\ref{tab:sens_cm} shows that our proposed model consistently achieves the smallest expected cost per mission under all the parameter settings.
% It is worth mentioning that the optimal policy of our model tends to continue the mission when the mission failure penalty $C_\text{m}$ is large, e.g., $C_\text{m}=8,000$.
% In this case, operational costs obtained from our policy and most benchmarks would be close.
{Among the benchmarks, the one-phase approximation performs the best when $C_\text{m}$ is small, and the $\mathcal{R}$-policy performs the best when $C_\text{m}$ is large.
Nevertheless, all benchmarks lead to a cost increase by around 2\%--13\% compared with our method.}

\begin{table}
	\centering
	\caption{Comparison of the operational costs per mission when the mission failure cost $C_{\text{m}}$ changes {and $F$ is Weibull}.
	The results are averaged under 10,000 simulation replications.}
	\label{tab:sens_cm}
	\begin{tabular}{cccccc}
		\hline	
		Policy & $C_\text{m}=500$ & $C_\text{m}=1,000$ & $C_\text{m}=2,000$& $C_\text{m}=4,000$& $C_\text{m}=8,000$ \\
		\hline
		$\mathcal{C}$-policy & 522.8 & 717.7 & 1063.0 & 1775.9 & 3063.6 \\
		$\mathcal{R}$-policy & 547.6 & 753.5 & 1089.6 & 1724.0 & 2952.7 \\
		%\midrule
		$\mathcal{M}$-policy & 503.8 & 721.5 & 1063.4 & 1730.1 & 2968.9 \\
		One-phase approximation & 500.5	& 717.3 & 1061.4 & 1741.2 & 3015.6 \\
		%\midrule
		Proposed model & 494.3 & 697.6 & 1013.4 & 1635.0 & 2783.3 \\
		\hline
	\end{tabular}%
\end{table}%

We then compare our model with the benchmarks under different values of the state-observation matrix $\mathbf{D}$.
We let the false alarm probability $d_{12}$ vary in $\{0.05,0.1,\ldots,0.3\}$, compute $d_{11}=1-d_{12}$, while keeping other parameters the same as those in Section \ref{subsec:num_setting}.
Table~\ref{tab:sens_false} compares the expected cost per mission of the {five} policies.
%Since $\mathcal{R}$-policy incorrectly overlooks the condition monitoring noises, the cost under this policy is very sensitive to $d_{12}$.
Although the $\mathcal{M}$-policy considers the noises in condition monitoring, the resulting cost is consistently more than 5\% higher than that of our optimal policy, because it overlooks the non-Markovian nature of the failure process.
We see that the operational cost per mission of our model and the $\mathcal{M}$-policy generally increases in $d_{12}$, because the signal-noise-ratio decreases in $d_{12}$.
%This phenomenon is more obvious for the $\mathcal{R}$-policy, where the cost increases rapidly in $d_{12}$.
%This is because the system is highly likely to receive a false alarm at an early stage of the mission when $d_{12}$ is large.
%As such, the rescue procedure follows and the chance of system failure is small.
{Moreover, the optimal $\mathcal{C}$-policy and one-phase approximation generally perform the best among all benchmarks.
The cost increase of the benchmarks are around 5\%--8\% compared with our method.}
%In the limiting case when $d_{12}$ approaches one, the mission is aborted whenever seeing consecutive warning signals,
%and thus the cost per mission would converge to the mission failure cost $C_\text{m}=1,000$.

\begin{table}
	\centering
	\caption{Comparison of the operational costs per mission on the false alarm probabilities $d_{12}$ {when $F$ is Weibull}.
	The results are averaged under 10,000 simulation replications.}
	\label{tab:sens_false}
	\begin{tabular}{ccccccc}
		\hline	
		Policy & $d_{12}=0.05$ & $d_{12}=0.1$ & $d_{12}=0.15$ & $d_{12}=0.2$ & $d_{12}=0.25$ & $d_{12}=0.3$ \\
		\hline
		$\mathcal{C}$-policy & 1048.5 & 1053.3 & 1055.5 & 1056.8 & 1062.2 & 1065.3 \\
		$\mathcal{R}$-policy & 1071.7 & 1075.1 & 1082.5 & 1087.9 & 1088.5 & 1092.2 \\
		%\midrule
		$\mathcal{M}$-policy & 1050.0 & 1055.8 & 1059.9 & 1058.8 &1062.6 & 1067.4 \\
		One-phase approximation & 1050.7 & 1053.9 & 1054.5 & 1059.6 & 1059.2 & 1061.4 \\
		%\midrule
		Proposed model & 996.5 & 1001.0 & 1005.3 & 1008.1 & 1011.2 & 1015.7 \\
		\hline
	\end{tabular}%
\end{table}%

Similar to the case examining $d_{12}$, we also conduct sensitivity analysis on the false negative rate $d_{21}$.
Table~\ref{t2} shows the comparison results, where the operational costs per mission as a function of $d_{21}$ for the benchmarks exhibit similar patterns to that when $d_{12}$ changes.
{The $\mathcal{C}$-policy performs the best when $d_{21}$ is small, and the one-phase approximation performs the best when $d_{21}$ increases.
Moreover, we see around 5\%--9\% increases of the operational cost for the benchmarks compared with our model, validating the effectiveness of the proposed method again.}
%We mention is passing that the expected cost under $\mathcal{R}$-policy is high even when the false negative rate $d_{21}=0$.
%aThis suggests that the high expected cost of $\mathcal{R}$-policy mainly comes from the ignorance of false alarms in our problem setting.

\begin{table}
	\centering
	\caption{Comparison of the operational costs per mission on the false negative rate $d_{21}$ {when $F$ is Weibull}.
	The results are averaged under 10,000 simulation replications.}
	\label{t2}
	\begin{tabular}{ccccccc}
		\hline	
		Policy & $d_{21}=0.05$ & $d_{21}=0.1$ & $d_{21}=0.15$ & $d_{21}=0.2$ & $d_{21}=0.25$ & $d_{21}=0.3$ \\
		\hline
		$\mathcal{C}$-policy & 1059.0 & 1063.0 & 1067.6 & 1071.2 & 1078.2 & 1086.6
		 \\
		$\mathcal{R}$-policy & 1081.7 & 1089.4 & 1094.2 & 1102.2 & 1104.2 & 1107.0 \\
		%\midrule
		$\mathcal{M}$-policy & 1062.4 & 1063.1 & 1065.9 & 1071.8 & 1074.6 & 1084.1 \\
		One-phase approximation & 1059.9 & 1061.3 & 1064.9 & 1070.8 & 1070.9 & 1080.3 \\
		%\midrule
		Proposed model & 1010.1 & 1013.4 & 1016.6 & 1020.1 & 1021.3 & 1024.8 \\
		\hline
	\end{tabular}%
\end{table}%

{We then examine the sensitivity of the our model performance compared to the benchmarks when $F$ follows a mixture of Weibull distributions.
We again fix $C_\text{s}$ and examine different values of $C_\text{m}$ for $C_{\text{m}}/C_{\text{s}}\in\{0.25,0.5,1,2,4\}$.
The results are summarized in Table~\ref{tab:sens_Cm_mixture}.
The $\mathcal{R}$-policy outperforms other benchmarks when $C_\text{m}\geq 1,000$.
However, the cost gaps between the benchmarks and our method can be very large when $F$ follows a mixture of Weibull distributions.
Nearly all benchmarks lead to a cost increase larger than 10\% in different parameter settings.
When $C_\text{m}=8,000$, the gaps between our method and the $\mathcal{C}$-policy, $\mathcal{M}$-policy, and one-phase approximation are all near 20\%.

We then examine the sensitivity of our model performance to $d_{12}$ and $d_{21}$.
The results are summarized in Tables~\ref{tab:sens_2peak_d12} and \ref{tab:sens_2peak_d21}, respectively.
We can again see that the optimal $\mathcal{R}$-policy outperforms other benchmarks in all parameter settings. 
Nevertheless, the optimal $\mathcal{R}$-policy still incurs a cost increase consistently around 12\% in all parameter settings compared with our model. 
Similar to the base setting, all other benchmarks lead to a cost increase more than 15\% under all parameter settings.
This demonstrates the effectiveness of our model when the time from the defective to failure state is bimodal.}

\begin{table}
	\centering
	\caption{{Comparison of the operational costs per mission when the mission failure cost $C_{\text{m}}$ changes and $F$ is a mixture distribution}.
	The results are averaged under 10,000 simulation replications.}
	\label{tab:sens_Cm_mixture}
	\begin{tabular}{cccccc}
		\hline	
		Policy & $C_\text{m}=500$ & $C_\text{m}=1,000$ & $C_\text{m}=2,000$& $C_\text{m}=4,000$& $C_\text{m}=8,000$ \\
		\hline
		$\mathcal{C}$-policy & 530.7 & 906.6 & 1293.1 & 1968.2 & 3290.4 \\
		$\mathcal{R}$-policy & 565.3 & 894.3 & 1248.6 & 1866.7 & 3099.3 \\
		%\midrule
		$\mathcal{M}$-policy & 506.8 & 908.2 & 1295.2 & 1970.7 & 3279.6 \\
		One-phase approximation & 524.3 & 909.2 & 1294.6 & 1874.5 & 3291.3 \\
		%\midrule
		Proposed model & 502.1 & 814.1 & 1116.4 & 1693.9 & 2797.1 \\
		\hline
	\end{tabular}%
\end{table}%

\begin{table}
	\centering
	\caption{The operational cost per mission when the parameter $d_{12}$ changes {and $F$ is a mixture distribution}.
	The results are averaged under 10,000 simulation replications.}
	\label{tab:sens_2peak_d12}
	\begin{tabular}{ccccccc}
		\hline	
		Policy & $d_{12}=0.05$ & $d_{12}=0.1$ & $d_{12}=0.15$ & $d_{12}=0.2$ & $d_{12}=0.25$ & $d_{12}=0.3$ \\
		\hline
		$\mathcal{C}$-policy & 1275.1 & 1278.2 & 1284.5 & 1289.5 & 1290.1 & 1301.6 \\
		$\mathcal{R}$-policy & 1235.2 & 1234.9 & 1242.6 & 1245.7 & 1247.5 & 1249.7 \\
		%\midrule
		$\mathcal{M}$-policy & 1273.6 & 1283.4 & 1283.5 & 1288.4 & 1290.0 & 1293.4 \\
		One-phase approximation & 1272.3 & 1284.0 & 1282.0 & 1289.2 & 1293.3 & 1294.2 \\
		%\midrule
		Proposed model & 1096.2 & 1103.6 & 1104.1 & 1106.9 & 1110.2 & 1114.3 \\
		\hline
	\end{tabular}%
\end{table}%

\begin{table}[htbp]
	\centering
	\caption{The operational cost per mission when the parameter $d_{21}$ changes {and $F$ is a mixture distribution}.
	The results are averaged under 10,000 simulation replications.}
	\label{tab:sens_2peak_d21}
	\begin{tabular}{ccccccc}
		\hline	
		Policy & $d_{21}=0.05$ & $d_{21}=0.1$ & $d_{21}=0.15$ & $d_{21}=0.2$ & $d_{21}=0.25$ & $d_{21}=0.3$ \\
		\hline
		$\mathcal{C}$-policy & 1290.2 & 1293.0 & 1300.2 & 1307.4 & 1308.4 & 1314.0 \\
		$\mathcal{R}$-policy & 1247.8 & 1248.6 & 1251.2 & 1252.0 & 1257.9 & 1260.9 \\
		%\midrule
		$\mathcal{M}$-policy & 1293.2 & 1294.9 & 1298.7 & 1302.3 & 1303.7 & 1309.7 \\
		One-phase approximation & 1284.0 & 1294.5 & 1297.4 & 1301.9 & 1303.8 & 1311.6 \\
		%\midrule
		Proposed model & 1115.1 & 1116.4 & 1120.3 & 1122.1 & 1126.5 & 1131.1 \\
		\hline
	\end{tabular}%
\end{table}%

{Finally, we examine the optimal operational cost under various values of $\lambda$ and $m_2$, the two parameters determining the Erlang mixture approximation.
This numerical experiment reveals how fast our surrogate policy converges to the optimal policy in the non-asymptotic regime.
For the Weibull distribution, we study the value of $m_2\in\{5,10,\ldots,35\}$ and compute the corresponding $\lambda$ by matching the mean with the Weibull distribution in Section~\ref{subsec:num_setting}.
When $F$ follows the mixture distribution, we set the value of $m_2\in\{10,\ldots,70\}$ as it needs more phases for the mixture of Erlang distribution to accurately approximate $F$ with a more irregular shape, as evidenced in Figure~\ref{fig:inspection}.
The results are summarized in Tables~\ref{tab:sens_m2} and \ref{tab:sens_m2_2}.
We can see that the $\lambda$ value obtained from moment-matching increases with $m_2$ in both tables.
This implies that $\lambda$ chosen in this way satisfies the conditions in Corollary~4.1 of \citet[]{wiens1979distributions},
and it ensures that $F^{(\lambda)}$ approximates $F$ well with a large enough $m_2$.
Moreover, when $F$ is Weibull, the expected cost per mission decreases significantly as $m_2$ increases from 5 to 20, and almost remains stable when $m_2$ increases from 20 to 35.
This implies that our surrogate policy converges quickly to the optimal policy for moderate $m_2$ values.
% We note that the costs in Table~\ref{tab:sens_m2} are not strictly decreasing in $m_2$, because of simulation noises and approximation errors of the PBVI algorithm.
On the other hand, when $F$ is a mixture distribution, it is expected that the expected cost by adopting the surrogate policy converges to the optimal cost more slowly. 
In this case, the mean cost per mission nearly keeps constant if we increase $m_2$ from $50$ and $\lambda$ from $0.209$, again validating our model performance in the non-asymptotic regime.}

\begin{table}
	\centering
	\caption{The operational cost per mission when the parameter $m_2$ changes {and $F$ is Weibull}.
		The results are averaged under 10,000 simulation replications.}
	\label{tab:sens_m2}
	\begin{tabular}{cccccccc}
		\hline	
		$m_{2}$& 5 & 10 & 15 & 20 & 25 & 30 & 35 \\
		\hline
		$\lambda$ & 0.041 & 0.074 & 0.105 & 0.134 & 0.163 & 0.191 & 0.218 \\
		Mean cost per mission & 1061.4 & 1044.8 & 1027.8 & 1013.4 & 1013.4 & 1013.1 & 1013.1 \\
		\hline
	\end{tabular}%
\end{table}%

\begin{table}
	\centering
	\caption{The operational cost per mission when the parameter $m_2$ changes {and $F$ is a mixture distribution}.
	The results are averaged under 10,000 simulation replications.}
	\label{tab:sens_m2_2}
	\begin{tabular}{cccccccc}
		\hline	
		$m_{2}$& 10 & 20 & 30 & 40 & 50 & 60 & 70 \\
		\hline
		$\lambda$ & 0.054 & 0.095 & 0.134 & 0.172 & 0.209 & 0.245 & 0.281 \\
		Mean cost per mission & 1256.6 & 1200.4 & 1155.2 & 1127.0 & 1116.4 & 1111.8 & 1110.8 \\
		\hline
	\end{tabular}%
\end{table}%

\bibliographystyleappendix{informs2014} % outcomment this and next line in Case 1
\bibliographyappendix{ref2} % if more than one, comma separated

\clearpage

\end{APPENDICES}

\end{document}